\documentclass[10pt]{article}
\usepackage{graphicx}
\usepackage{amsmath}
\usepackage{graphicx}
\usepackage{esint}
\usepackage{xcolor}
\usepackage{amsfonts}
\usepackage{amssymb}
\usepackage{amsmath, amssymb ,amsthm, amsfonts, amsgen}
\usepackage{mathrsfs}
\usepackage{hyperref}
\usepackage{mathtools}
\usepackage{subcaption}
\usepackage{verbatim} 
\usepackage{enumitem}

\DeclareGraphicsExtensions{.eps,.bmp,.jpg,.pdf,.mps,.png,.gif}
\numberwithin{equation}{section} 
\makeatletter
\newsavebox{\@brx}
\newcommand{\llangle}[1][]{\savebox{\@brx}{\(\m@th{#1\langle}\)}%
	\mathopen{\copy\@brx\kern-0.5\wd\@brx\usebox{\@brx}}}
\newcommand{\rrangle}[1][]{\savebox{\@brx}{\(\m@th{#1\rangle}\)}%
	\mathclose{\copy\@brx\kern-0.5\wd\@brx\usebox{\@brx}}}
\makeatother

\let\pa\partial  
  
\let\eps\varepsilon  
\newcommand{\N}{{\mathbb N}} 
\newcommand{\R}{{\mathbb R}} 
\newcommand{\Z}{{\mathbb Z}} 
\newcommand{\diver}{\operatorname{div}}
\newcommand{\trace}{\operatorname{tr}}
\newcommand{\dist}{\operatorname{dist}}

\newcommand{\dd}{{\mathrm{d}}}
\newcommand{\sym}{\operatorname{sym}}

\newcommand{\mat}[1]{\boldsymbol #1}
\newcommand{\vect}[1]{\boldsymbol #1}

\newcommand{\NN}{\mathbb{N}}

\newcommand{\RR}{\mathbb{R}}
\newcommand{\ZZ}{\mathbb{Z}}

\newcommand{\ee}{\varepsilon}

\theoremstyle{plain}
\newtheorem{teo}{Theorem}[section]

\theoremstyle{definition}
\newtheorem{definition}[teo]{Definition}

\newtheorem{assumption}[teo]{Assumption}
\newtheorem{remarkica}[teo]{Remark}
\newtheorem{propozicija}[teo]{Proposition} 
\newtheorem{lema}[teo]{Lemma} 
\newtheorem{corollary}[teo]{Corollary} 
\newtheorem{teorem}[teo]{Theorem}

\usepackage[titletoc]{appendix}

\newcommand{\BBB}{\color{black}}
\newcommand{\CCC}{\color{black}}

\mathtoolsset{showonlyrefs}

\begin{document}
	\title{Poroelastic plate model obtained by simultaneous homogenization and dimension reduction}
	
	\author{Marin Bu\v{z}an\v{c}i\'{c} \footnote{Faculty of Chemical Engineering and Technology, University of Zagreb, Marulićev trg 19, 10000 Zagreb, Croatia. Email: buzancic@fkit.hr} \and Pedro Hern\'andez-Llanos\footnote{Instituto de Ciencias de la Ingenier\'ia, Universidad de O'Higgins, Avenida Libertador Bernardo O'Higgins 611, 2841935 Rancagua, Chile. E-mail: pedro.hernandez@uoh.cl}\and Igor Vel\v ci\'c\footnote{Faculty of Electrical Engineering and Computing, 
			University of Zagreb, Unska 3,  10000 Zagreb, Croatia. E-mail: igor.velcic@fer.hr}\and Josip \v Zubrini\'c\footnote{Faculty of Electrical Engineering and Computing, 
			University of Zagreb, Unska 3,  10000 Zagreb, Croatia. E-mail: josip.zubrinic@fer.hr}}
	\date{\today}

	\maketitle{\textit{Dedicated to the memory of Andro Mikeli\'c.}}

\begin{abstract}
In this paper, the starting point of our analysis is a coupled system of linear elasticity and Stokes equation. We consider two small parameters: the thickness $h$ of the thin plate and the pore scale $\varepsilon(h)$ which depends on $h$. We will focus specifically on the case when the pore size is comparatively small relative to the thickness of the plate. The main goal here is derive a model of a poroelastic plate, starting from the $3D$ problem as $h$ goes to zero, using simultaneous homogenization and dimension reduction techniques. The obtained model generalizes the poroelastic plate model derived by A. Mikeli\'c et. al. in 2015 using dimension reduction techniques from $3D$ Biot's equations in the sense that it also covers the case of contacts of poroelastic and (poro)elastic plate as well as the evolution equation with inertial term. 



\end{abstract}

\section{Introduction}

The derivation of Biot’s equations using the homogenization approach is a relatively old topic (see \cite{gilbmikelic2000},\cite{clopeaumikelic2001}). Recently, there has been growing attention devoted to deriving and studying poroelastic equations in various contexts. We will mention just a few of them. In \cite{gahn1}, the authors investigate a thin layer of an elastic body immersed in a fluid modeled by Stokes equations. In a recent work \cite{gahn2024}, the author derives the poroelastic plate model in a regime where the size of the pores is on the same scale as the thickness of the body. In \cite{gurvich},   the quasi-static and evolutionary equations for poroelastic plates are analyzed. In \cite{Avalos},  the existence result for the interaction between fluid and poroelastic structures is discussed. In \cite{bociu}, the authors analyze the quasi-static Biot system of poroelasticity for both compressible and incompressible constituents. In \cite{Mikelic2012}, the model for the contact between poroelastic and elastic bodies is derived. Finally, in \cite{Mikelic2015}, the authors derive the model of quasi-static Biot’s plates by performing a dimension reduction of $3D$ Biot’s equations.

The current paper focuses on deriving Biot's plate equations under the condition that the size of the pores, denoted as $\varepsilon(h)$, is significantly smaller than the thickness of the plate, denoted as $h$, that is: $\varepsilon(h) \ll h$. 
The starting point is the coupled system consisting of the evolution Stokes equation for incompressible fluid and the linearized equations of elasticity. Utilizing modifications of Griso's decomposition (referenced in \cite{Grisosdecomp2005} and \cite{Vel14a}, \cite{mbukalvelcic2017},   \cite{Marinvelciczubrinic2022}), we perform simultaneous homogenization and dimension reduction. Griso's decomposition has previously been employed to obtain the model of homogenized plates within the contexts of both linear and non-linear elasticity (see \cite{Vel14a,mbukalvelcic2017} for the case of linear and non-linear elasticity in moderate contrast and \cite{Marinvelciczubrinic2022} for the case of the elastic plate in high contrast). 

In the limit (resulting in the quasi-static case), we derive a generalized version of Biot's plate model obtained in \cite{Mikelic2015}. Our derivation also encompasses scenarios involving the contact between elastic and poroelastic plates, as well as interactions between two distinct poroelastic plates. Consequently, the resulting limit equations are more intricate and do not exhibit the same decoupling observed in the equations derived in \cite{Mikelic2015}, and also analyzed in \cite{gurvich}.

The derivation of the poroelastic plate model via simultaneous homogenization and dimension reduction also furnishes additional insights into the limit Darcy's law and appropriate boundary conditions, as well as interface conditions in the case of contacts, establishing their connection with the microscopic model.

There is a notable difference between the quasi-static model obtained in this paper and the one derived in \cite{gahn2024}, where  the macroscopic pressure depends only on macroscopic in-plane variables, leading to different limit equations. In this paper (akin to \cite{Mikelic2015}), the limit equations solely involve the derivative of pressure in the vertical direction. More precisely,  the microscopic fluid velocity in the effective Darcy's law is driven only by the derivative of the pressure in the vertical direction.

One of the challenges encountered in deriving the model was to establish the correct scalings of the constants  (i.e. the relation of the elasticity constants with respect to viscosity) which gives Biot's plate model in the limit, since we have to scale them with respect to two small parameters $\eps$ and $h$.

In our model derivation, we also address the scenario resulting in the plate model with inertial terms, thereby providing justification for the model analyzed in \cite{gurvich}. When considering models with inertial terms, it's well recognized that the full $3D$ macroscopic Biot's equations exhibit memory effects (see \cite{clopeaumikelic2001}).
In this context, there exists a certain analogy between poroelastic equations and linear elastic equations with high-contrast inclusions. Namely, to obtain Biot's equations, viscosity is put in high contrast with elasticity constants. However, the limit plate equations derived here do not exhibit memory effects (see also \cite[Section 3.4.1]{Marinvelciczubrinic2022}). This absence can be attributed to the fact that we derive the model in the bending regime, which is valid for long times (we implicitly scale time for the derivation of evolution equations with inertial terms). 
The techniques outlined in this paper could also enable us to derive models in the membrane regime (referenced in \cite{Marinvelciczubrinic2022}), where it is expected that memory effects for macroscopic equations would be present. However, these equations are degenerate in the vertical component of deformation, thus their significance for the applications is not clear (see \cite[Section 3.4.2]{Marinvelciczubrinic2022}).

For analyzing the limit equations, we adopt the approach of \cite{gurvich}. Additional effort is required due to the more generalized nature of the equations compared to those analyzed in \cite{gurvich}. This is particularly true for equations incorporating inertial terms.

Next, we briefly outline the structure of the paper. In Section \ref{secsetting}, we discuss the problem's setting and provide some\emph{ a priori }estimates for fluid-elastic structure interaction. Section \ref{secquasi-static} is dedicated to deriving and analyzing the quasi-static Biot's plate model, while Section \ref{secinertia} focuses on deriving and analyzing the Biot's plate model with inertial term in the bending regime. The Appendix contains auxiliary claims essential for dimension reduction, such as Griso's decomposition and its consequences, along with auxiliary definitions and claims about two-scale convergence.

The main results of the paper are summarized as follows. In the quasi-static case
Theorem \ref{gammazeromay11} provides the compactness result; Theorem \ref{effectiveequationsmay11} provides the convergence result, identifying the limit model;
Theorem \ref{existeo} provides uniqueness and existence results and
Theorem \ref{proplenka1} and Theorem \ref{twoscaleconvergencesmay16} provide the strong convergence result.
Similarly in the case with inertial term
Theorem \ref{compinertia} provides the compactness result and the limit Biot's equations; 
Theorem \ref{uniqueinertia} provides the uniqueness result;
Theorem \ref{existinertia} provides the existence result and
Theorem \ref{stronginertia} provides the strong convergence result. \CCC We provide more details about the results in the next section. \BBB

\CCC \subsection{Techniques and novelties}
 
 The limit models for Biot's plate are derived  using simultaneous homogenization and dimension reduction (see Theorem \ref{effectiveequationsmay11} and Theorem \ref{compinertia}). In the quasi-static case our model extends the model obtained in \cite{Mikelic2015} (see Remark \ref{usporedba} below), obtained by doing dimension reduction from $3D$ macroscopic Biot's equations (adapting the techniques from \cite{ciarlet2}). In the case with inertial term this is, to the best of our knowledge, the first justification of such a model. 
We emphasize the fact that, starting from the full $3D$ Biot's model with inertial term and doing dimension reduction,  it is not a priori clear how the limit plate model would look like, since the full $3D$ Biot's model  has memory effects (see \cite{clopeaumikelic2001} and Remark \ref{interpret1} below). 
Doing simultaneous homogenization and dimension reduction additionally enables us to discuss influence of microscopic boundary conditions on the limit model and to obtain limit Darcy's law (see Remark \ref{remdarcy} and Section \ref{secgen} below). 

 To obtain the limit model, we adapt the method used in \cite{Vel14a,mbukalvelcic2017,Marinvelciczubrinic2022}. The method consists in using modification of Griso's decomposition to obtain the compactness statements in Theorem \ref{gammazeromay11} and Theorem \ref{compinertia}. It relies on the decomposition of sequences with bounded symmetrized scaled gradients and obtaining the strongest compactness statement for such sequences (the results are recalled and adapted in Appendix \ref{griso}). Unlike the approach of \cite{ciarlet2} which focuses more on obtaining limit displacement and justifying Kirchoff-Love ansatz, the method we use here focuses more on obtaining appropriate compactness results for sequences with bounded symmetrized scaled gradients, while Kirchoff-Love ansatz is incorporated in these results.  In order to obtain such results, we need to use certain estimates valid on $h$ level, the information on the limit displacement is not sufficient (see Theorem \ref{aux:thm.griso} below).  

 We also obtain the limit model for contacts between different types of poroelastic cells and identify two main types of contacts under which we are able to obtain the limit boundary condition for the pressure at the interface (one of them is Neumann condition  and the other is continuity of the pressure). They correspond to two different situations, depending whether there is a flow at the interface or not. 
To the best of our knowledge this is the first treatment of such  problems in the context of poroelasticity.

 The limit models are derived under the minimal regularity assumptions on the loads (the results in \cite{clopeaumikelic2001,Mikelic2012,gahn2024} require more regularity in time for the loads). 
On the microscopic level, we assume that the forces possess $L^2$ or $H^1$ regularity in time, depending on whether we want to obtain in the limit quasi-static case or case with inertial term (which is, anyhow, quasi-static in in-plane components). 
More precisely, to obtain the quasi-static limit equations, we require boundedness in $H^1$ in time for the loads, while in the evolution case, boundedness in $H^1$ in time is only enforced for the in-plane components of the loads and for the vertical component, we impose only $L^2$ boundedness in time. It's worth noting that this requirement is less restrictive than the ones in \cite{Mikelic2012} and \cite{gahn2024}, which require higher regularity. This minimal requirement is exactly the one that is necessary for the existence result of the limit equations.
Consequently, on the microscopic level we have to deal with weak solution defined in \cite{DUGun}, which doesn't guarantee that the pressure posses $L^2$ regularity in time (unless the forces have $H^1$ regularity in time). 
Also, the pressure converges weakly only in $H^{-1}$ in time, and an additional regularity of the limit pressure is attained from the limit equations themselves (see the proof of Theorem \eqref{effectiveequationsmay11}).
We also don't use the condition that the loads at zero moment are zero (cf. \cite{clopeaumikelic2001,Mikelic2012,gahn2024}) except for the strong convergence results in Theorem \ref{proplenka1}, Theorem \ref{twoscaleconvergencesmay16} and Theorem \ref{stronginertia}, where it is necessary.

  To obtain the existence and uniqueness result for the limit problems we use the approach from \cite{gurvich} (see Theorem \ref{existeo} and Theorem \ref{existinertia}). This approach required additional considerations since our limit equations are more complex than the one analyzed in \cite{gurvich} (in particular in-plane and vertical components of displacement do not decouple).  Moreover, in the case with inertial term we prove uniqueness directly (see Theorem \ref{uniqueinertia}, the proof is not given in \cite{gurvich}) and prove existence for less regular loads by using Galerkin approximation (and relying on semigroup approach as suggested in \cite{gurvich}), see Remark \ref{obsshane} for comparison between semigroup approach and approach via Galerkin approximation.

 The strong convergence result is proved in Theorem \ref{proplenka1}, Theorem \ref{twoscaleconvergencesmay16} and Theorem \ref{stronginertia} (as in \cite{clopeaumikelic2001,Mikelic2012} in the context of $3D$ Biot's equations). This requires to prove certain energy-type equality (see Proposition \ref{propen} and Proposition \ref{propen2} below), which has to be done by the approximation argument for the solutions that don't posses enough regularity (as a consequence of less regular loads than the one considered e.g., in \cite{clopeaumikelic2001,Mikelic2012}). This approximation is done via semigroup approach.

\BBB

\subsection{Notation}
\label{notation}

Throughout this paper, we denote by $\NN$, $\ZZ$, $\RR$, $\RR^{m\times n}$ and $\RR^{n\times n}_{\rm sym}$ the sets of natural, integer, real,  the space of real $m\times n$ matrices and the space of real $n\times n$ symmetric matrices respectively, this two last spaces endowed with the usual Euclidean norm $|\vect{F}|=\sqrt{\text{tr}\, \vect{F}^T\vect{F}}$. If $\mathbf{A}, \mathbf{B} \in \RR^{n \times m}$, we denote by $\mathbf{A}:\mathbf{B}$ the scalar product 
$$ \mathbf{A}:\mathbf{B}=\text{tr} (\mathbf{A}^T \mathbf{B}). $$
For $\mathbf{A} \in \RR^{n \times n}$, $\sym \mathbf{A}$ denotes the symmetric part of $\mathbf{A}$. 
We denote the coordinate vectors by $\vect{e}^i$, $i=1,2,3$.  Unless otherwise stated, the  Greek  indices $\alpha,\beta$ take values in the set $\{1,2\}$. For  the  vectors  $x=(x_1,x_2,x_3) \in \RR^3$, we denote by $\widehat{x}=(x_1,x_2)$  the vector containing the first two variables. If the vector valued function $\vect{u}$ is taking values in $\RR^3$,  we denote the first two components by $\vect{u}_*=(u_1,u_2)$. 
If $\mathbf{a}, \mathbf{b} \in \RR^n$, we denote by $\mathbf{a} \odot \mathbf{b}$ the symmetric part of the tensor product  $\mathbf{a} \otimes \mathbf{b}$, i.e. $\mathbf{a} \odot \mathbf{b} := \sym (\mathbf{a} \otimes \mathbf{b} )$. If $\vect{a},\vect{b} \in \R^3$, we denote by $\vect{a} \wedge \vect{b}$ their cross product. 

With $\iota: \RR^{2 \times 2} \to \RR^{3 \times 3}$, we denote the natural inclusion 
$$ \iota(\mathbf{A})_{\alpha\beta}=\mathbf{A}_{\alpha \beta}, \quad \iota(\mathbf{A})_{i3}= \iota(\mathbf{A})_{3i}=0,\quad \alpha,\beta \in\{1,2\},\ i \in \{1,2,3\}. $$

For $A \subset \RR^n$, we denote by $\bar{A}$ its closure and  by $|A|$ its Lebesgue measure. Furthermore, $\chi_A$ denotes the characteristic function. The open ball of radius $r > 0$ centered at $x\in \RR^k$ is denoted by  $B(x,r)$.

We denote the unit interval by $ I = \left(-\frac{1}{2}, \frac{1}{2}\right)$, the unit cell by $ Y = [0, 1)^3 $, and the flat torus with quotient topology by $ \mathcal{Y} = \mathbb{R}^3 / \mathbb{Z}^3 $. 
For a subset $\mathcal{Y}'\subset \mathcal{Y}$, we denote by $Y'$ the associated subset of $Y$ through the identification map. 
\CCC This means that $\mathcal{Y}'=\{y+\Z^3 :  y \in Y'\}$. \BBB
\CCC  Analogously we define $\widehat{Y}=[0,1]^2$, $\widehat{\mathcal{Y}}=\RR^2/\ZZ^2$. 
Again, for $y \in \mathcal{Y}$, $\hat{y}$ denotes $\hat{y}=(y_1,y_2)$. 
\BBB

For  an  open set $A$ and  $k \in \NN\cup \{\infty\}$,   $C^k(A)$ denotes the set of  $k$-times continuously  differentiable functions on $A$.
If $A$ is a closed set, then $C^k(A)$ is defined as the space of functions that are restrictions of functions in $C^k(A')$ for some open set $A' \supset A$.
The space $C_c^k(A)$ denotes the set of $k$-times differentiable functions with compact support in $A$. The space $C^k(\mathcal{Y})$ denotes the set of $k$-times differentiable functions on the torus, while $H^k(\mathcal{Y})$ denotes the closure of $C^k(\mathcal{Y})$ in the $H^k$ Sobolev norm.
If $\omega \subset \RR^2$ is a rectangle  (i.e. $\omega = (a,b)\times (c,d)$), then the set $C^k_{\#}(\omega)$ is the set of  $k$-times continuously differentiable functions with bounded derivatives that can be extended to  periodic functions in $C^k(\RR^2)$, with period $\omega$.
In addition, for an open subset 
$B \subset \RR$, we denote by $C^k_{\#}(\omega \times B)$ the set of $k$-times  differentiable functions on $\omega \times B$ that can be extended to  functions in $C^k(\RR^2 \times B)$ which are periodic in $(x_1,x_2)$-variables with period $\omega$.  The spaces $H^k_{\#}(\omega)$ and  $H^k_{\#}(\omega \times B)$ are the closures of sets  $C^k_{\#}(\omega)$ and $C^k_{\#}(\omega \times B)$ in $H^k$ norm, respectively.  
The set $\dot{C}^k(A)$ denotes the subspace of functions in $C^k(A)$ with zero mean. 
In an analogous way we define the spaces $\dot{C}_c^k(A)$, $\dot{H}^k(A)$, $\dot{H}_{\#}^k(\omega)$ $\dot{H}^k_{\#}(\omega\times B)$, $B \subset \R$. These definitions of functional spaces are naturally extended for spaces of functions taking values in $\RR^k$. 
In general, $\rightarrow$ denotes the strong convergence and $\rightharpoonup$ denotes the weak (or weak-*) convergence. 

 For $A \subset \R^n$ and $f \in L^1(A)$, $\fint_A f$ denotes the mean value of $f$ on the set $A$.   
For $f \in L^1(A \times I)$, $A \subset \RR^n$, we denote by $\bar{f} \in L^1(A)$ the function\begin{equation} \label{defbar} \bar{f}=\int_I f \,dx_3, \end{equation} 
and if $f \in L^1(Y)$, we denote by 
\begin{equation} \label{deflangle}
\langle f \rangle_Y:=\int_Y f \,dy.
\end{equation}

For a smooth function $f$,  its support is denoted by $\textrm{supp} f$. 

For a Banach space $V$, we denote its dual by $V'$.
If $f \in V'$ and $v \in V$, then $_{V'}\langle f, v \rangle_V$ denotes the value $f(v)$.
The set of continuous linear operators from vector space $V$ to vector space $W$ is denoted by $L(V,W)$.
Given an operator $\mathcal{A}$, we denote by $\mathcal{D}(\mathcal{A})$ and $\textrm{Ker}(\mathcal{A})$ its domain and kernel.
For  $h>0$,  we denote by $\pi_h:\RR^3 \to \RR^3$ the mapping 
\begin{equation}
\label{pimapping}
    \pi_h(\mathbf{a})=(ha_1,ha_2, a_3).
\end{equation}
For a function $\vect{u}:\Omega \subset \RR^3 \to \RR^3$, we denote  its symmetric gradient by $e(\vect{u}):=\sym \nabla \vect{u}$.  In the analogous  way, we define $e_{\widehat{x}} (\vect{u})$ and $e_y(\vect{u})$  for the mappings $\vect{u}:\omega \subset \RR^2 \to \RR^2$ and $\vect{u}: Y(\mathcal{Y})\to \RR^3$, respectively. 
For functions $\vect{u}:\Omega \subset \RR^3 \to \RR^k$, $k \in \N$, and $h>0$  we denote by $\nabla_h$ the scaled gradient 
\begin{equation} \label{defscaledg} 
\nabla_h \vect{u}:=\left(\partial_1 \vect{u},\partial_2 \vect{u}, \frac{1}{h} \partial_3 \vect{u} \right). 
\end{equation}
If $\vect{u}$ is taking values in $\RR^3$, we denote the scaled divergence and the symmetric scaled gradient by 

\begin{equation}\label{defscalede} 
\diver_h \vect{u}:=\text{tr} \,\nabla_h \vect{u}, \quad e_h(\vect{u}):=\sym \nabla_h \vect{u}. 
\end{equation} 

Throughout Section \ref{secsetting}, Section \ref{secquasi-static} and Section \ref{secinertia} we consider $\omega \subset \RR^2$ a rectangle with vertices  having  integer  coordinates, and assume that there exist $k\in \NN$ and $l\in \NN$ such that $h=\frac{1}{k}$ and $\eps=\frac{1}{kl}$. 
We will write $\eps=\eps(h)$ and assume in Section \ref{secquasi-static} and Section \ref{secinertia} that $\lim_{h \to 0} \frac{\eps(h)}{h}=0$. 

We define $\Omega^h=\omega \times (hI)$ for $h > 0$, and $\Omega =\omega \times I$.

\section{Setting of the problem}
\label{secsetting} 

In this section, we present the microscopic problem. In Section \ref{geometry}, we discuss the geometry of the microscopic problem. In Section \ref{microscopic}, we provide the equations of the microscopic problem and establish the existence result. Finally, in Section \ref{apriori}, we derive some\emph{ a priori }estimates necessary to obtain the limit equations. 

\subsection{Geometry and material assumptions} \label{geometry} 

\medskip
\subsubsection{ The poroelastic cells and interface conditions}

\begin{figure}[ht]
\centering\includegraphics[draft=false,width=0.9\textwidth]{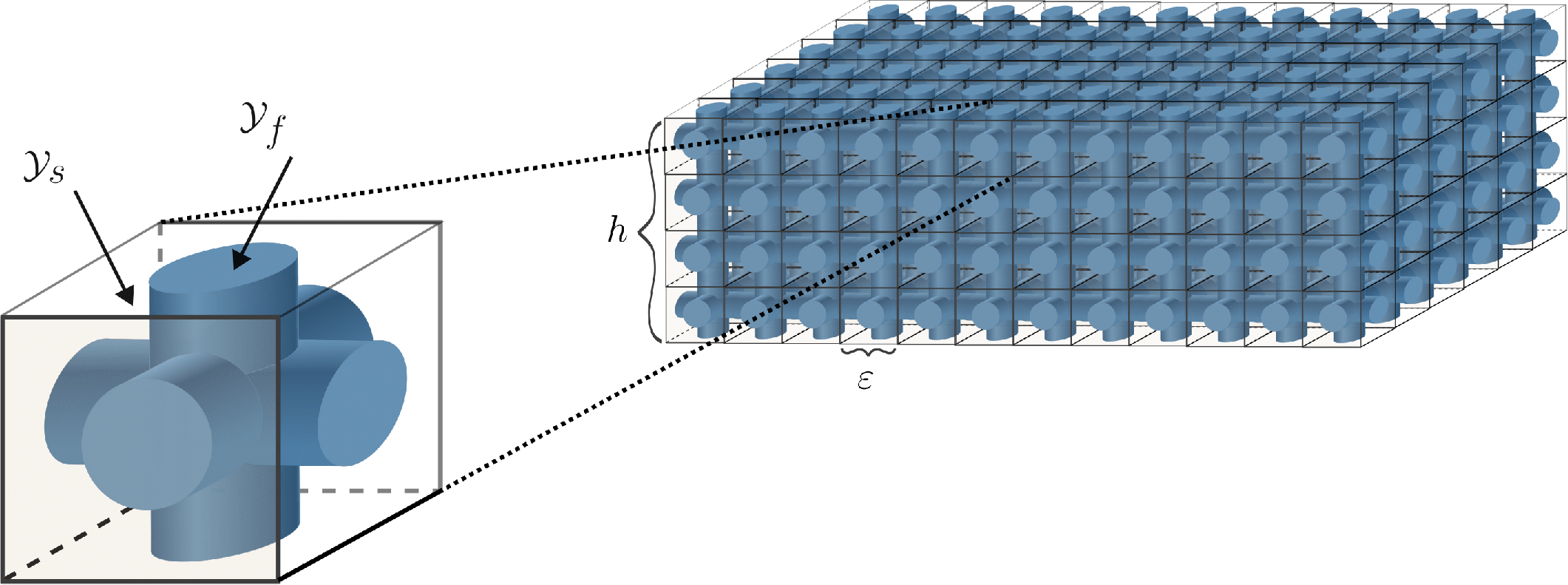}
	\caption{Illustration of the poroelastic plate, i.e. the fluid and the solid phases}
 \label{figure1} 
\end{figure}

\begin{figure}[h!]
	\centering
	\begin{subfigure}{\linewidth}
		\centering
		\begin{subfigure}{0.3\linewidth}
			\centering
			\includegraphics[draft=false,trim={10cm 4cm 10cm 3cm},clip,width=0.7\linewidth]{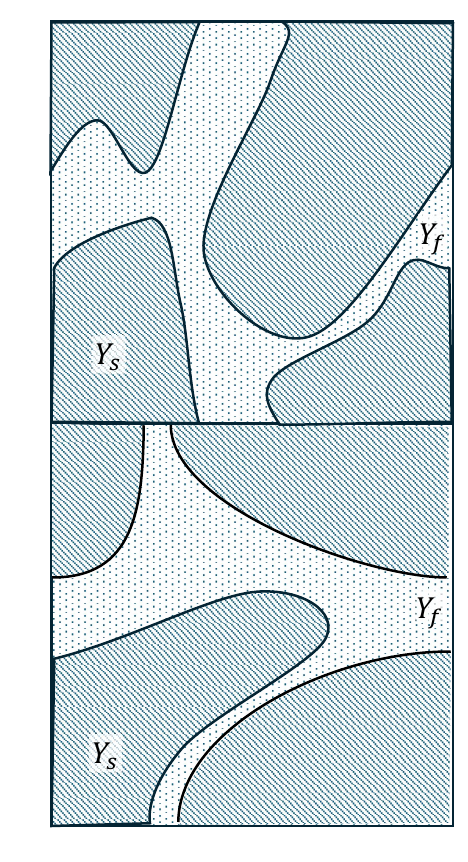}
			\caption*{(a)}
		\end{subfigure}
		\begin{subfigure}{0.3\linewidth}
			\centering
			\includegraphics[draft=false,trim={10cm 4cm 10cm 3cm},clip,width=0.7\linewidth]{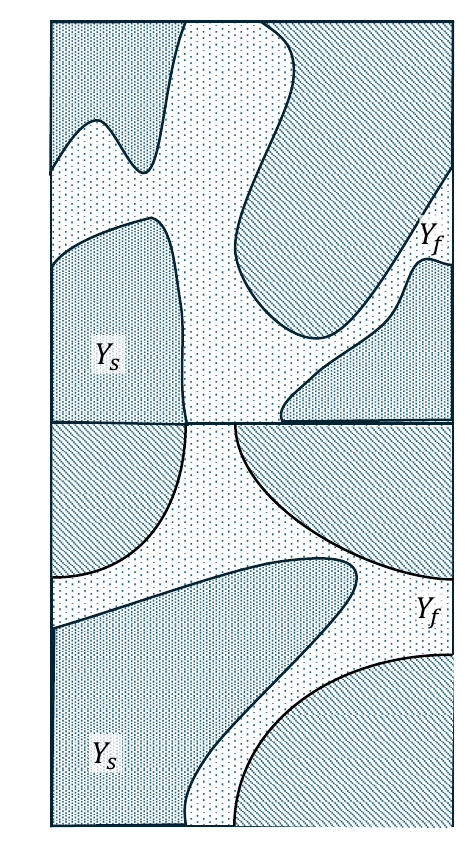}
			\caption*{(b)}
		\end{subfigure}
	\end{subfigure}
	\caption{(a) is an example of interface where there is no flow, while (b) is an example of interface such that there is a flow. Note also that the cells in (b) don't have the opposite boundaries matching.  
		\BBB}
	\label{figurecf}
\end{figure}

To extend our analysis to encompass interactions between various types of poroelastic plates, including interactions between poroelastic and elastic plates, we permit our material to consist of diverse poroelastic materials.
To this end, we  assume the existence of  a finite  number of typical cells consisting of fluid and solid  part. In other words, we assume that  there exist  $m \in \N$  and pairs $(\mathcal{Y}_f^i,\mathcal{Y}_s^i)$, $i=1,\dots,m$ of open  subsets of $\mathcal{Y}$ with Lipschitz boundary such that 
(see Figure \ref{figure1})
$$ \mathcal{Y}_f^i \cap \mathcal{Y}_s^i=\emptyset, \quad \overline{\mathcal{Y}_f^i} \cup \overline{\mathcal{Y}_s^i}=\mathcal{Y}, \quad i=1,\dots, m, $$
i.e. there exist pairs $({Y}_f^i,{Y}_s^i)$\footnote{Recall, we assume that $Y_f^i$, $Y_s^i$ are associated sets of $\mathcal{Y}_f^i$, $\mathcal{Y}_s^i$ through the identification map (\CCC recall Section \ref{notation}). \CCC Note that the closures of $Y_f^i$ i.e., $Y_s^i$ (in the topology of $\R^3$) don't need to have opposite boundaries matching. This happens if $\mathcal{Y}_f^i$ (i.e. $\mathcal{Y}_s^i$) has part of its boundary (in the torus topology) on $y_j=0$, $j=1,2,3$, see Figure \ref{figurecf} (b).  \BBB}, $i=1,\dots,m$, of open sets of $Y$ (in the relative topology) with Lipschitz boundary, such that
$$ {Y}_f^i \cap {Y}_s^i=\emptyset, \quad \overline{{Y}_f^i} \cup \overline{{Y}_s^i}={Y},\quad  i=1,\dots, m.  $$ 
We also assume that both $Y_f^i$ and $Y_s^i$ (consequently, both $\mathcal{Y}_f^i$ and $\mathcal{Y}_s^i$ \CCC) \BBB are connected sets. 

The set $\Omega^h$ is  partitioned into   $\widetilde{\Omega}_s^h$ and $\widetilde{\Omega}_f^h$,  representing the regions occupied by solid and fluid, respectively.  For each  $h > 0$, these sets are defined in the following way.  Let $i_h:\mathbb{Z}^3 \to \{1,\dots,m\}$ be a given map, which selects and associates a typical cell type to every index $k \in \Z^3$. 
We set
\begin{equation}
\widetilde{\Omega}^h_{f/s}= \bigcup_{k} \varepsilon\,(k+Y_{f/s}^{i_h(k)}),
\end{equation}
where the union is taken over all $k \in \mathbb{Z}^3$ such that $\varepsilon(k+Y) \subset  \Omega$. 
Note that the dependence on $\eps$ is concealed as a dependence on $h$ (since $\eps = \eps(h)$), and will remain so throughout the paper.
We assume  that  both $\widetilde{\Omega}^h_{s}$ and $\widetilde{\Omega}^h_{f}$  are Lipschitz and that $\widetilde{\Omega}_s^h$ is also connected. 
We  denote the interface between the two phases with  $\widetilde{\Gamma}^h=\partial \widetilde{\Omega}^h_f \cap \partial \widetilde{\Omega}^h_s$. In addition, we will assume that the fluid domain does not intersect the upper and lower boundaries, i.e.
\begin{equation} \label{nodownup} 
\partial \widetilde{\Omega}^h_{f} \cap \left(\{x_3=-h/2\}\cup\{x_3=h/2\}\right)=\emptyset.	
\end{equation} 
We discuss the consequences of dropping this condition in Section \ref{secintersection} below. 
Furthermore, we introduce sets  $\Omega^h_{f/s}\subset  \Omega$ as the images of $\widetilde{\Omega}^h_{f/s}$ through the rescaling 
\begin{equation} \label{rescaling} 
	(x_1, x_2, x_3)=(x_1^h, x_2^h, h^{-1}x_3^h),\quad (x^h_1, x^h_2, x^h_3)\in \Omega^h.
\end{equation}
Additionally we assume layer-like structure of our composite poroelastic material. 
Thus, we impose the following compatibility condition on the limiting structure: 
\begin{equation}
\left|\chi_{\Omega^h_{f/s}}(x)-\sum_{i=1}^m \chi_{U_i}(x) \chi_{\mathcal{Y}^i_{f/s}}\left(\frac{\widehat{x}}{\ee}, \frac{x_3}{\frac{\ee}{h}}\right)\right| \to 0 \textrm{ almost everywhere in } \Omega, 
\label{karakt} 
\end{equation}
where $U_1=\omega \times I_1,\dots,U_m=\omega \times I_m \subset \Omega$ 
denote the layers of the material. \BBB 
Here $I_1, \dots  I_m$
are disjoint, open, non-empty sub-intervals of $I$ of the form:
\begin{equation} \label{intervals} 
 I_i=(d_{i},d_{i+1}), \quad i=1,\dots,m, \quad -\frac{1}{2}= d_1< \dots< d_{m+1}=\frac{1}{2}.
 \end{equation} 

\begin{remarkica}
    In other words, the sets $U_1, \dots U_m$ model regions in $\Omega$ occupied by a poroelastic material with cell type $\mathcal{Y}^1_{f/s}, \dots, \mathcal{Y}^m_{f/s}$, respectively. 
    We allow that some of $\mathcal{Y}_f^i$, for $i=1,\dots,m$, are empty (meaning that the region $U_i$ is purely elastic).
    By the same analysis presented in this paper, it is also possible to consider more general types of regions than the ones defined here (see Section \ref{secgen}), which we don't analyze for the simplicity of the exposure.
\end{remarkica}
 However, for the derivation of the model we will need more restrictive assumption than \eqref{karakt}. 
 \CCC Namely, we will impose layer--like structure for every $h$. \BBB
 
Before stating the required assumption, we introduce the following notation. 
For a subset $S \subset Y$ we denote by $\widetilde{S}_h^{\#}$ the set 
$$ \widetilde{S}_h^{\#}:= \bigcup_{k \in \ZZ^3} \eps\,(k+S), $$
while by $S_h^{\#}$ we denote the image of $\widetilde{S}_h^{\#}$ through the rescaling \eqref{rescaling}.

\CCC
\begin{assumption} \label{assumption on regions}
For every $h$ there exists numbers $(d_i^h)_{i=1,\dots,m+1}$, such that 
\begin{equation} \label{intervals2} 
 -\frac{1}{2}= d_1^h< \dots< d_{m+1}^h=\frac{1}{2},  
 \end{equation}
and we have
$$\Omega^h_{f/s} \cap (\omega \times I_i^h) = \left(Y_{f/s}^{i}\right)_h^{\#}\cap (\omega \times I_i^h), \quad \forall h<h_0, \quad i = 1,\dots,m,$$
where 
$I_i^h=(d_{i}^h,d_{i+1}^h), \quad i=1,\dots,m.$
\end{assumption}
The above assumption tells us that whole region $\omega \times I_i^h$ consists only of poroelastic cells of type $i$.  	As a consequence of \eqref{karakt} we have that 
\begin{equation} \label{intervals3} 
 d_i^h \to d_i, \textrm{ as } h \to 0, \quad i=1,\dots,m+1. \end{equation} 
We distinguish two different type of contacts between the regions (see also Figure \ref{figurecf}): 
More precisely, if for $i\in {1,\dots,m-1}$ we have that $\mathcal{Y}_f^i \cap \mathcal{Y}_f^{i+1}\cap\{y_3=0\} =\emptyset$, then we say that there is \emph{no flow  at the interface}. \footnote{\CCC Strictly speaking the set $\{y_3 =0\}$, looked as the subset of $\mathcal{Y}$, is actually the set $\{y+\Z^3:y_3 =0\}.$\BBB} On the other hand, 
if for $i\in {1,\dots,m-1}$ we have that $\mathcal{Y}_f^i \cap \mathcal{Y}_f^{i+1}\cap\{y_3=0\} \neq \emptyset$, then we say that there is \emph{a flow at the interface}. \footnote{\CCC Note that, since both $\mathcal{Y}_f^i$ and $\mathcal{Y}_f^{i+1}$ are open, there is actually a ball contained in $\mathcal{Y}_f^i \cap \mathcal{Y}_f^{i+1}$ whose center belongs to the plane $y_3=0$.   \BBB}

\begin{figure}[ht]
	\begin{minipage}[c]{0.5\textwidth}
		\centering
		\begin{subfigure}{0.5\linewidth}
			\centering
			\includegraphics[draft=false,trim={13.1cm 3cm 10cm 2cm},clip,height= 10cm]{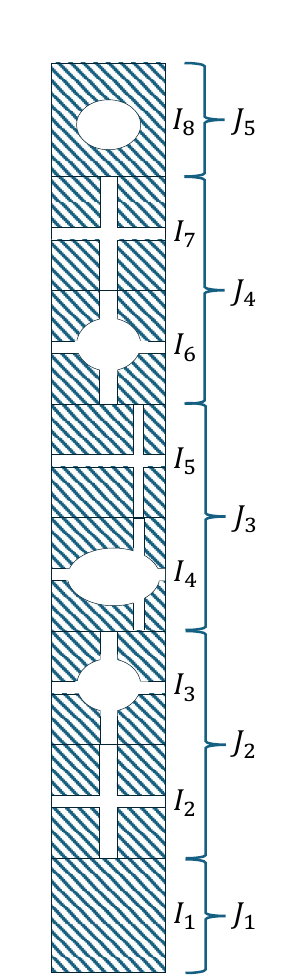}
		\end{subfigure}
	\end{minipage}\hfill
	\begin{minipage}[c]{0.5\textwidth}
	\caption{ \CCC In this figure we explain poroelastic region and the domain of vertical permeability. 
		For simplicity of the picture we drew that each  $I_i$, $i=1,\dots,8$ has only one (layer of) cell.  
		Here 
			$J_p=J_2 \cup J_3 \cup J_4\cup J_5$ and $J_K=J_2 \cup J_3 \cup J_4$ and $J_i$, for $i=1,\dots,5$ are assumed to be open intervals.  \BBB }
	\label{visokeslike}
	\end{minipage}
\end{figure}
\BBB
\subsubsection{The poroelastic region and the domain of vertical permeability}

Next we wish to emphasize the crucial part of the domain on which pressure is not zero.
This set (which we denote by $\Omega_p$) is obtained by excluding from $\Omega$ layers of purely elastic medium.
More precisely,
if there is a flow on the interface between $U_i$ and $U_{i+1}$, we join them into one region (including their interface). 
After doing this procedure for every $i \in \{1,\dots,m\}$ we obtain the regions $V_1,\dots,V_s$, $s\leq m$, of the form $V_i=\omega \times J_i$, for $j=1,\dots,s$. 
Here $J_i$ is an open, non-empty sub-interval of $I$ that is a union of two or more neighbouring intervals from the family $\{I_1,\dots,I_m\}$ including their endpoints (except the left endpoint of the first interval and the right endpoint of the last interval in the union).
We then define
\begin{equation} \label{defjp} 
    \Omega_p :=\omega \times J_p, \quad
    J_p := \bigcup_{\begin{matrix} \scriptstyle i \in\{1,\dots,s\}\\[-1ex] \scriptstyle V_i \textrm{ not purely elastic} \end{matrix}} J_i.
\end{equation} 
We refer to $\Omega_p$ as the poroelastic region (\CCC see Figure \ref{visokeslike}). \BBB


Furthermore, we will distinguish layers of $\Omega_p$ with cells which allow the fluid flow in $x_3$-direction.
We will denote this subset as $\Omega_K$, and refer to it as the domain of vertical permeability.
In order to do this, first we define the set $J_K \subset J_p$ by excluding those $I_i$, $i=1,\dots,m$, for which we have that  $\mathcal{Y}_f^i \cap \{y_3=0\}=\emptyset$, i.e. 
\begin{equation} \label{defjk} 
J_K:=J_p \backslash \bigcup_{\begin{matrix} \scriptstyle i \in\{1,\dots,m\}\\[-1ex] \scriptstyle \mathcal{Y}_f^i \cap \{y_3=0\}=\emptyset \end{matrix}} I_i
\end{equation} 
and we define 
$$ \Omega_K:=\omega \times J_K.$$

As we will see below, $\Omega_p$ is the region where the coefficient $M_0$, appearing in the limit equations, is uniformly strictly greater than zero, while $\Omega_K\subset \Omega_p$ is the region where the coefficient $\mathbb{K}_{33}$ is strictly greater than zero  (both coefficients $M_0$ and $\mathbb{K}_{33}$ are defined in Section \ref{effective} below).  On the region $\Omega\backslash \Omega_p$ both coefficients are equal to zero and furthermore, $\mathbb{K}_{33}$ is equal to zero on $\Omega_p \backslash \Omega_K$ (see Proposition \ref{tensorprop}). The importance of uniform positivity of these coefficients is seen in Section \ref{analysis}.
\begin{remarkica}\label{marindod1}
Note that $\Omega \backslash \Omega_p$ (i.e. $\Omega \backslash \Omega_K$) might have an isolated set of the form $\omega \times \{d_i\}$, for some $i=2,\dots,m$, as its subset \CCC(see again Figure \ref{visokeslike}). \BBB   
For this reason $\Omega_p$ might not be Lipschitz domain. 
Also, the functions belonging to $H^1_{\#}(\Omega_p)$ might have a jump in the trace on the part of the boundary of the form $\omega \times \{d_i\}$. Consequently, the space $C_{\#}^1(\bar{\Omega}_p)$ is not  dense in $H_{\#}^1(\Omega_p)$. However, the functions in $C_{\#}^1(\Omega_p)$, which together with their derivatives belong to $L^{\infty} (\Omega_p)$, are dense in $H_{\#}^1(\Omega_p)$.
The space $H^1_{\#}(\Omega_p)$ can be easily understood, since we have 
$$ H^1_{\#}(\Omega_p)=\bigoplus _{\begin{matrix} \scriptstyle i \in\{1,\dots,s\}\\[-1ex] \scriptstyle  V_i \textrm{ not purely elastic} \end{matrix}} H_{\#}^1(V_i). $$
The analogous observations are valid for the set $J_p$, i.e. $H^1(J_p)$ (also for $\Omega_K,J_K$ and $H_{\#}^1(\Omega_K),H^1(J_K)$). For the consequence of these observations on the limit problem see Remark \ref{remgotovo}.
\end{remarkica}

\subsubsection{The elasticity tensor}

The elastic properties of the material are modeled with  elasticity tensors $\mathbb{A}^1,\dots, \mathbb{A}^m$ that satisfy
\begin{equation}\label{tensorA}
 \exists \nu>0 \textrm{ such that }	\nu |\vect{\xi}|^2\leq \mathbb{A}^i(y)\vect{\xi}: \vect{\xi} \leq \nu^{-1}  |\vect{\xi}|^2,\quad \forall \vect{\xi}\in \RR^{3\times 3}_{\text{sym}}, \quad y \in Y,  i=1,\dots, m. 
\end{equation}
We  assume that the following standard symmetries hold 
\begin{equation}
	\mathbb{A}^s_{ijkl}(y)=\mathbb{A}^s_{jikl}(y)=\mathbb{A}^s(y)_{klij},\quad y \in Y,  i,j,k,l\in\{1,2,3\}, \quad s=1, \dots,m.
\end{equation}
 Finally, we define 
$$\widetilde{\mathbb{A}}^h(x^h)=\sum_{k \in \ZZ^3,\, \eps( k+Y) \subset \Omega^h}\mathbb{A}^{i_h(k)}(x^h/\eps-k)\chi_{\eps(k+Y_s^{i_h(k)})}(x^h), \quad x^h \in \Omega^h. $$ 
As a consequence of \eqref{karakt}, i.e. Assumption \ref{assumption on regions}, we have 
\begin{equation} \label{depra} 
	\left|\mathbb{A}^h(x)-\mathbb{A}\left(\CCC x_3 \BBB,\frac{\widehat{x}}{\ee}, \frac{x_3}{\frac{\ee}{h}}\right)\right|\to 0, \textrm{ almost everywhere in } \Omega,
\end{equation} 
where 
$  \mathbb{A}^h(x):=  \widetilde{\mathbb{A}}^h(x_1,x_2,hx_3)$ and
\CCC
\begin{equation} \label{defAAA} 
 \mathbb{A}(x_3,y):=  \chi_{I_1}(x_3) \mathbb{A}^1(y) \chi_{\mathcal{Y}_s^1}(y)+\cdots+ \chi_{I_m}(x_3)\mathbb{A}^m(y) \chi_{\mathcal{Y}_s^m}(y).
\end{equation} 
Note that 
For $x_3 \in I_i$, we  define  $\mathcal{Y}_{f/s}(x_3):=\mathcal{Y}_{f/s}^i$, and for $x_3=d_i$, $i=1,\dots,m+1$, $\mathcal{Y}_f(x_3):=\emptyset$, $\mathcal{Y}_s(x_3):=\mathcal{Y}$.  
We also define the set 
 \begin{equation}
 \label{yxspace}
    \Omega \times \mathcal{Y}_{f/s}^{x_3} :=\{(x,y) \in \Omega \times \mathcal{Y}: y \in \mathcal{Y}_{f/s}(x_3)\}.
 \end{equation}
 \BBB
All the functions on $L^2(\Omega \times \mathcal{Y}_{f/s}^{x_3})$ we can naturally extend by zero and consider them as functions in $L^2(\Omega \times \mathcal{Y})$ ($\equiv L^2(\Omega \times Y)$).
\CCC For given $x_3 \in I$ and $\mathcal{C} \subset \mathcal{Y}$ the space $H^1(\mathcal{C})$ is defined as the set of restrictions of functions  belonging to $H^1(\mathcal{Y})$ on $\mathcal{C}$.  
The set $H_0^1(\mathcal{C})$ is understood as the set of functions which, when extended by zero outside $\mathcal{C}$, belong to $H^1(\mathcal{Y})$. 
\BBB  
\subsection{The microscopic equations}

\label{microscopic} 
We consider the following fluid-solid-structure interaction problem on $\Omega^{h}$ which couples Stokes equation with linearized elasticity:
\begin{equation}\label{eq:fluid}
\displaystyle \eta \widetilde{\kappa}^h_{f}\partial_{tt}\widetilde{\vect{u}}^h+\frac{\nabla \widetilde{p}^h}{h^2}=\frac{\ee^{2}}{h^{4}} \Delta \partial_t \widetilde{\vect{u}}^h+\widetilde{\vect{F}}^{h},\quad   \displaystyle\text{div}\,\partial_t \widetilde{\vect{u}}^h=0\quad \text{in}\,\,\widetilde{\Omega}^{h}_f,
\end{equation}
\begin{equation}\label{eq:solid}
\displaystyle \eta\widetilde{\kappa}_{s}^{h}\partial_{tt}\widetilde{\vect{u}}^h=\frac{1}{h^2}\text{div}\,(\widetilde{\mathbb{A}}^he(\widetilde{\vect{u}}^h))+\widetilde{\vect{F}}^h\quad \text{in}\,\,\widetilde{\Omega}^{h}_s, 
\end{equation}
\begin{equation}\label{eq:tensors}
\widetilde{\sigma}^{f,h}:=-\frac{1}{h^2}\widetilde{p}^hI+2\frac{\ee^{2}}{h^4} e\left(\displaystyle\partial_t \widetilde{\vect{u}}^h\right)\quad \text{in}\,\,\widetilde{\Omega}^{h}_{f}, \quad \quad\widetilde{\sigma}^{s, h}:=\frac{1}{h^2}\widetilde{\mathbb{A}}^he(\widetilde{\vect{u}}^h)\quad\text{in}\,\,\widetilde{\Omega}_{s}^{h},
\end{equation}
\begin{equation}\label{eq:boundary}
 \displaystyle [\widetilde{\vect{u}}^h]=0, \quad \widetilde{\sigma}^{f,h}\vect{n}_f=\widetilde{\sigma}^{s,h} \vect{n}_s \quad\text{on}\,\,\widetilde{\Gamma}^{h}, \quad \quad \widetilde{\sigma}^{s,h}\vect{n}=0\quad\text{on}\,\,\{x_3=h/2\}\cup \{x_3=-h/2\},
\end{equation}
\begin{equation}\label{eq:initial}
\left.\widetilde{\vect{u}}^h\right|_{\{t=0\}}=\left.\partial_{t}\widetilde{\vect{u}}^h\right|_{\{t=0\}}=0\quad\text{on}\,\, \Omega^{h}, \quad \quad  \left\{\widetilde{\vect{u}}^h, \widetilde{p}^h\right\}\quad \text{is $\omega$-periodic in }\, (x_1,x_2).
\end{equation}
Here, $\widetilde{\vect{u}}^h|_{\Omega_s^{h}}$ represents the deformation in the solid, while $\partial_t \widetilde{\vect{u}}^h|_{\Omega_f^{h}}$ denotes the fluid velocity. The function $\widetilde{p}^h$ denotes pressure, and $\widetilde{\vect{F}}^h$ represents volume forces. The parameters $\widetilde{\kappa}_f^h$ and $\widetilde{\kappa}_s^h$ denote the fluid and solid densities, respectively, which are positive functions \CCC(for the analysis below they don't have to be constants) \BBB  bounded from below and above \CCC by a positive constant .\BBB

The parameter $\eta=\eta(h)>0$ serves as a  time-scaling parameter. One scenario we analyze is when $\lim_{h \to 0} \eta(h)=0$, which results in the quasi-static regime. The other scenario is $\eta(h)=1$, leading to the bending regime with inertial term.
\CCC In Section \ref{apriori} we will prove {\it a  priori} estimates and see how these two regimes influence them. Later, in Section \ref{secetanula} and Section \ref{secinertia} we will analyze these two regimes separately, In Section \ref{secetanula} we will analyze the regime $\lim_{h \to 0} \eta(h)=0$, while in Section \ref{secinertia} we will analyze the case $\eta(h)=1$.     \BBB

The tensors $\widetilde{\sigma}^{f,h}$ and $\widetilde{\sigma}^{s,h}$, given by \eqref{eq:tensors}, represent the stresses in the fluid and solid, respectively. The vector $\vect{n}$ in  the second and third equation of \eqref{eq:boundary} denotes the unit normal at the interface point, from different sides of the interface (in the second equation) or at the upper or lower boundary point (in the third equation).

The equations \eqref{eq:fluid} are Stokes equations for incompressible fluid, while \eqref{eq:solid} is the equation of linearized elasticity (with the appropriate scaling of the constants applied to all of them). The first equation in \eqref{eq:boundary} expresses the continuity of the deformation at the interface (here $[\widetilde{\vect{u}}^h]$ denotes the jump of the function, i.e. the difference of traces from two sides of the interface), while the second equation in \eqref{eq:boundary} represents the continuity of stresses in the normal direction. The third equation in  \eqref{eq:boundary} is the Neumann boundary condition at the upper and lower boundary (this can be modified by adding surface loads, see Section \ref{secgen}), while the first equation in  \eqref{eq:initial} represents the zero initial condition (our analysis also accommodates other initial conditions).  For simplicity, periodic boundary conditions are imposed on the transverse boundary by the second \CCC expression \BBB in \eqref{eq:initial}. Different initial and boundary conditions are discussed in Section \ref{secgen}.\BBB 

\begin{remarkica}
     It may not be immediately evident what is the motivation behind the particular scaling of coefficients in \eqref{eq:fluid}-\eqref{eq:tensors}, with respect to $\varepsilon$ and $h$. The main feature of this chosen scaling is that it leads to the Biot's plate model as $h \to 0$.  One way to interpret it is to consider it as a scaling of viscosity and time (or density) variables (to see this, one has to multiply equations \eqref{eq:fluid} and \eqref{eq:solid} by $h^2$).
\end{remarkica}

\subsubsection{ The weak formulation on the physical domain}

 The functional space for \eqref{eq:fluid}-\eqref{eq:initial} is the following:
\begin{equation} 
	\widetilde{V}^{h}:=\left\{\psi \in H^1_{\#}(\Omega^h;\RR^3) : \diver \psi=0 \textrm{ in } \widetilde{\Omega}_f^{h} \right\}. 
\end{equation}
The variational formulation which corresponds to (\ref{eq:fluid})-(\ref{eq:initial}) is given by:\\
For $\widetilde{\vect{F}}^h \in L^2(0,T; L^2(\Omega^h;\RR^3))$ find $\widetilde{\vect{u}}^h\in L^{\infty}(0,T;\widetilde{V}^h)$ such that  $\partial_t \widetilde{\vect{u}}^h \in  L^2(0,T;L^2(\Omega^h;\RR^3)) \cap L^2(0,T;H_{\#}^1(\widetilde{\Omega}^h_{f};\RR^3))$ and $\partial_{tt} \widetilde{\vect{u}}^h\in L^2(0,T; (\widetilde{V}^h)')$ and
\begin{equation}
\label{varformporousnonrescaled}
  \begin{split}    & \eta\widetilde{\kappa}^h \CCC_{(\widetilde{V}^h)'} \langle \partial_{tt}\widetilde{\vect{u}}^h(t),\vect{v}\rangle_{\widetilde{V}^h}\BBB+\displaystyle\frac{\ee^2}{h^{4}}\int_{\widetilde{\Omega}^{h}_{f}}2 e(\partial_t\widetilde{\vect{u}}^h(t)):e(\vect{v})\,dx^h+\frac{1}{h^{2}}\int_{\widetilde{\Omega}_{s}^{ h}}\widetilde{\mathbb{A}}^h(x^h) e(\widetilde{\vect{u}}^h(t)):e(\vect{v})\,dx^h\\
&=\int_{\Omega^h}\widetilde{\vect{F}}^{h}\vect{v}\,dx^h,\quad\forall \vect{v}\in \widetilde{V}^h,\quad\text{a.e.}\,\,\text{in}\,\,(0,T),
  \end{split}  
\end{equation}
where \CCC
\begin{equation} \label{defkappa}
\widetilde{\kappa}^{ h}=\widetilde{\kappa}^h_{f}\chi_{\widetilde{\Omega}^{ h}_{f}}+\widetilde{\kappa}^h_{s}\chi_{\widetilde{\Omega}^{h}_{s}}.
\end{equation} \BBB
\CCC For the analysis below it is enough to assume that $\widetilde{\kappa}^h$ is a function bounded from below and above by a positive constant. \BBB
The variational formulation is supplemented with initial conditions (\ref{eq:initial}).
The existence and uniqueness of the solution of the above problem is given in \cite[Theorem 2.4]{DUGun}. 
\begin{remarkica} \label{rempressure} 
Here we discuss the existence and uniqueness of the pressure variable \CCC for the problem \BBB  \eqref{varformporousnonrescaled} \CCC (cf. \eqref{eq:fluid}) \BBB. 
\CCC Classically, as in \eqref{varformporousnonrescaled}, weak formulation can be written without pressure with divergence free test functions.  Pressure  then can be introduced as a Lagrange multiplier and it is important to establish its existence, uniqueness and regularity (see, e.g. \cite[Section 5]{giraultraviart1986} for Stokes equation). 
From the mathematical point of view, introducing the pressure as a variable enables us to write weak formulation for arbitrary test functions and to write the strong formulation.	\BBB
To this end, 
we introduce: 
\begin{equation} 
	\widetilde{H_T^1L^2_f}:=\left\{\varphi \in H^1(0,T;L^2(\widetilde{\Omega}^h_f)) : \varphi(T)=0\right\}.
\end{equation}	
By integrating \eqref{varformporousnonrescaled} over $[0,t]$, for $t \in [0,T]$, we can \CCC introduce \BBB the function $\widetilde{P}^h(t)$: 
\begin{equation}
\label{eqandrey1}
    \begin{split}
        & -\frac{1}{h^2} \int_{\widetilde{\Omega}^h_f} \widetilde{P}^h(t) \diver \vect{v} \, dx^h= -\eta\int_{\Omega^h}\widetilde{\kappa}^h\partial_{t}\widetilde{\vect{u}}^h(t)\vect{v}\,dx^h-\displaystyle\frac{\ee^2}{h^{4}}\int_{\widetilde{\Omega}^{h}_{f}}2 e(\widetilde{\vect{u}}^h(t)):e(\vect{v})\,dx^h\\
&-\frac{1}{h^{2}}\int_{\widetilde{\Omega}_{s}^{ h}}\widetilde{\mathbb{A}}^h(x^h) e\left(\int_0^t \widetilde{\vect{u}}^h(\tau)\, d\tau\right):e(\vect{v})\,dx^h+\int_{\Omega^h}\int_0^t \widetilde{\vect{F}}^{h}(\tau)\, d\tau \vect{v}\,dx^h,  \quad \vect{v} \in H^1_{\#}(\Omega^h;\RR^3). 
    \end{split}
\end{equation}
By using the same approach as in \cite[Section 3.2]{DUGun}, we can conclude the existence of $\widetilde{P}^h \in L^2(0,T;L^2(\widetilde{\Omega}^h_f))$, which satisfies  \eqref{eqandrey1}, \CCC provided that we have the existence of the solution of \eqref{varformporousnonrescaled}. \BBB
\CCC Note that \eqref{eqandrey1} is simply the integrated (in time) variant  of \eqref{varformporousnonrescaled} with arbitrary test functions and with "pressure" variable (actually $\widetilde{P}^h$ is the integral in time of physical pressure, see below, and consequently has more regularity in time, namely $L^2$).  This integration in time was done to obtain the existence of $\widetilde{P}^h$, by using the results of  \cite[Section 3.2]{DUGun}.
In order to come to version of \eqref{varformporousnonrescaled} with pressure variable we need to do the derivation in time of \eqref{eqandrey1}. This is done in the following way: \BBB 
One tests \eqref{eqandrey1} with $\partial_t \vect{v}(t)$, where $\vect{v}\in H^1(0,T;H^1_{\#}(\Omega^h;\RR^3))$ is such that  $\vect{v}(T)=0$, and integrates over $[0,T]$. 
By doing integration by parts in every term on the right hand side except the first one, and considering \eqref{varformporousnonrescaled}, we have that for some $\widetilde{p}^h \in (\widetilde{H^1_TL^2_f})'$ the following is satisfied:  
\begin{align}
	&\nonumber\displaystyle -\eta\int_0^T\int_{\Omega^h}\widetilde{\kappa}^h \partial_{t}\widetilde{\vect{u}}^h\partial_t \vect{v}\, dt\,dx^h+\displaystyle\frac{\ee^2}{h^{4}}\int_0^T\int_{\widetilde{\Omega}^{h}_{f}}2 e(\partial_t\widetilde{\vect{u}}^h):e(\vect{v})\,dx^h\, dt\\
	&\nonumber+\frac{1}{h^{2}}\int_0^T\int_{\widetilde{\Omega}_{s}^{ h}}\widetilde{\mathbb{A}}^h(x^h) e(\widetilde{\vect{u}}^h):e(\vect{v})\,dx^h\,dt-\displaystyle\frac{1}{h^2} \CCC
	_{(\widetilde{ H_T^1L^2_f })'}\langle \widetilde{p}^h,\diver \, \vect{v} \rangle_{ \widetilde{H_T^1L^2_f}}\BBB
	\\
	&\label{varformporousnonrescaled11}=\int_{\Omega^h}\widetilde{\vect{F}}^{h}\vect{v}\,dx^h\, dt,\quad\forall \vect{v}\in H^1(0,T;H^1_{\#}(\Omega^h;\RR^3)) \textrm{ such that } \vect{v}(T)=0.
\end{align}
%
\CCC Here: \BBB 
\begin{equation} \label{defphfrom}  
_{(\widetilde{\CCC H_T^1L^2_f \BBB})'}\langle \widetilde{p}^h, \varphi\rangle_{\CCC \widetilde{H_T^1L^2_f}\BBB}:=-\int_0^T\int_{\widetilde{\Omega}^{h}_{f} } \widetilde{P}^h \, \partial_t \varphi\, dx^h \,dt, \quad \forall \varphi \in \widetilde{H_T^1L^2_f}.
\end{equation} 	
We write formally $\widetilde{P}^h(t)=\int_0^t \widetilde{p}^h(\tau) \, d\tau. $ 

\CCC This simple observation and the existence  of this irregular pressure, under the condition that loads have only $L^2$ regularity in time, is not discussed in \cite{DUGun}.  
The authors of \cite{DUGun} obtain the existence of the solution  of the problem \eqref{varformporous2gamma0may10} below, by using only weak formulation with divergence free (on fluid part) test functions. 
However, the existence of  $L^2$ regular in time pressure for the loads that have $H^1$ regularity in time is obtained in \cite{DUGun}. It is again done by eliminating the divergence free condition in the space of test functions (i.e., from geometrical point of view by introducing the Lagrange multiplier).     \BBB
\end{remarkica} 	

\subsubsection{ The weak formulation on the fixed domain}

In order to work in a fixed domain $\Omega$,  we apply the change of variables \eqref{rescaling} 
and we define $\vect{u}^h(x):=\widetilde{\vect{u}}^{h}(x^h)$. In the analogous way we define $p^h$, $\vect{F}^h$, $\kappa^h$ (\CCC recall \eqref{defkappa}\BBB).
After doing this transformation we obtain the following problem: For $\vect{F}^h \in L^2(0,T; L^2(\Omega))$    find ${\vect{u}}^h\in L^{\infty}(0,T;{V}^h)$ such that  $\partial_t \vect{u}^h \in  L^2(0,T;L^2(\Omega;\RR^3)) \cap L^2(0,T;H_{\#}^1({\Omega}^h_{f};\RR^3))$, $\partial_{tt} \vect{u}^h\in L^2(0,T; ({V}^h)') $ and $p^h \in (H_T^1L^2_f)'$ satisfy (\CCC recall \eqref{defscalede}\BBB)
\begin{equation}
\label{varformporous2gamma0may10}
    \begin{split}
        	& \displaystyle-\eta \int_0^T \int_{\Omega} \kappa^{h} \partial_t \vect{u}^h\CCC \partial_t \vect{v}\,dx\BBB+\displaystyle\frac{\ee^2}{h^{4}}\int_0^T\int_{\Omega^{h}_{f}}2 e_h(\partial_t \vect{u}^h):e_h(\vect{v})\,dx\, dt+\frac{1}{h^{2}}\int_0^T\int_{\Omega_{s}^{h}}\mathbb{A}^h(x) e_h(\vect{u}^h):e_h(\vect{v})\,dx\, dt \\ &-\displaystyle\frac{1}{h^{2}} \CCC
        	_{( H_T^1L^2_f )'}\langle p^h,\diver_h \, \vect{v} \rangle_{ H_T^1L^2_f}\BBB=\int_0^T\int_{\Omega}\vect{F}^h\vect{v}\,dx\, dt, \quad \forall \vect{v}\in H^1(0,T;H^1_{\#}(\Omega^h;\RR^3)) \textrm{ such that } \vect{v}(T)=0. 
    \end{split}
\end{equation}
 The initial condition is given with 
\begin{equation}\label{eq:51rescaledver}
\vect{u}^h(0)=0\quad\text{on}\,\, \Omega.
\end{equation}
Here we denoted by $\Omega^h_{f/s}$ the rescaled $\widetilde{\Omega}^h_{f/s}$. 
 The space $V^h$  is defined in the following way 
 \begin{equation}
 	V^{h}:=\left\{\psi \in H^1_{\#}(\Omega;\RR^3) : \diver_h \psi=0 \textrm{ in } \Omega_f^{h} \right\} 
 \end{equation}
 and the space $H_T^1L^2_f$ as
 \begin{equation} 
 H_T^1L^2_f:=\left\{\varphi \in H^1(0,T;L^2(\Omega_f^h)) : \varphi(T)=0\right\}.
 \end{equation}	
We also define $\Gamma^h$ as rescaled $\widetilde{\Gamma}^h$. 

Next we discuss the existence result for the problem \eqref{varformporous2gamma0may10} -\eqref{eq:51rescaledver}. 
The following result is the consequence of \cite[Theorem 2.4, Theorem 3.2, Theorem 3.4]{DUGun} and Remark \ref{rempressure}. 

\begin{propozicija}\label{newpropmay10}
\begin{enumerate} 
\item Let us suppose that  $\vect{F}\in L^2(0,T;L^2(\Omega;\RR^3))$.  Then,  $3D$ problem given in (\ref{varformporous2gamma0may10}) has a unique solution $(\vect{u}^h, p^h)$ such that  ${\vect{u}}^h\in L^{\infty}(0,T;{V}^h)$ and $\partial_t \vect{u}^h \in  L^2(0,T;L^2(\Omega;\RR^3)) \cap L^2(0,T;H_{\#}^1({\Omega}^h_{f};\RR^3))$, $\partial_{tt} \vect{u}^h\in L^2(0,T; ({V}^h)') $ and $p^h \in (H_T^1L^2_f)'$. It additionally satisfies \footnote{The expression $t\mapsto\int_0^tp^h(\tau)\, d\tau$ is, as explained in Remark \ref{rempressure}, understood in the symbolic sense.  }
\begin{align*} 
	& \|\partial_{tt} \vect{u}^h\|_{L^2(0,T;(V^h)')}+\|\partial_{t} \vect{u}^h\|_{L^2(0,T;H^1(\Omega_f^{h};\RR^3))}+ \sup_{t \in[0,T]}\|\partial_t \vect u^h(t)\|_{L^2(\Omega;\RR^3)}+ \sup_{t \in[0,T]}\|\vect{u}^h(t)\|_{H^1(\Omega^h_s;\RR^3)}\\ & \hspace{+8ex}+\|p^h\|_{(H_T^1L^2)'}+\left\lVert \int_0^{\cdot} p^h(\tau) \,d\tau\right\rVert_{L^2(0,T;L^2(\Omega^h_f))}\leq    C(h,\Omega) e^{C(h,\Omega)T}\|\vect F\|_{L^2(0,T;  L^2(\Omega;\RR^3))},
\end{align*} 
for some $C(h,\Omega)>0$.
\item 	Let us suppose that $\vect{F}\in H^1(0,T;L^2(\Omega;\RR^3))$. Then,  the unique solution of 1. satisfies $(\vect{u}^h, p^h)\in \left(H^1(0,T; V^h) \cap H^2(0,T; L^2(\Omega;\RR^3))\right) \times L^2(0,T; L^2(\Omega^{h}_f))$. It additionally satisfies: 
\begin{align*} 
	& \sup_{t \in[0,T]}\|\partial_{tt} \vect{u}^h(t)\|_{L^2(\Omega;\RR^3)}+\|\partial_{tt} \vect{u}^h\|_{L^2(0,T;H^1(\Omega_f^{h};\RR^3))}+  \sup_{t \in[0,T]}\|\partial_t \vect u(t)\|_{H^1(\Omega_s^{h};\RR^3)}\\ &\hspace{+10ex}+\|p^h\|_{L^2(0,T;L^2(\Omega_f^{h}))}\leq  C(h,\Omega) e^{C(h,\Omega)T}\|\vect F\|_{H^1(0,T;  L^2(\Omega;\RR^3))},
\end{align*}
for some $C(h,\Omega)>0$. 
\end{enumerate} 
\end{propozicija} 	

\begin{remarkica} \label{obsjosip1}
	In the case when $\vect{F} \in H^1(0,T;L^2(\Omega;\RR^3))$ we have that the solution of \eqref{varformporous2gamma0may10}  satisfies the following equality pointwise in time
 \begin{equation}
 \label{varformporous2gamma0may1000}
     \begin{split}
         		\displaystyle \eta\int_{\Omega} \kappa^{h}&\partial_{tt}\vect{u}^h(t)\vect{v}\,dx+\displaystyle\frac{\ee^2}{h^{4}}\int_{\Omega^{h}_{f}}2 e_h(\partial_t \vect{u}^h(t)):e_h(\vect{v})\,dx +\frac{1}{h^{2}}\int_{\Omega_{s}^{h}}\mathbb{A}^h(x) e_h(\vect{u}^h(t)):e_h(\vect{v})\,dx \\
		&-\displaystyle\frac{1}{h^{2}}\int_{\Omega^{h}_{f} }p^h \diver_h\,\vect{v}\,dx =\int_{\Omega}\vect{F}^h\vect{v}\,dx,\quad\forall \vect{v} \in H^1_{\#}(\Omega;\RR^3),\quad\text{a.e.}\,\,\text{in}\,\,(0,T),
     \end{split}
 \end{equation}
	with initial conditions 
	$ \vect{u}^h(0)=\partial_t \vect{u}^h(0)=0. $
 This is the consequence of additional regularity of the solution in time given by Proposition \ref{newpropmay10} 2.
\end{remarkica}

\subsection{ \textit{A priori} estimates} \label{apriori} 

Here we provide \textit{a priori} estimates for the solutions to the problem \eqref{varformporous2gamma0may10} regardless of the assumptions on the time-scaling parameter $\eta=\eta(h)$. We define $\kappa \in L^\infty(\Omega)$ as weak star limit of $\kappa^h$.   It is a positive function bounded from below and above. $\bar{\kappa}$ is then defined according to \eqref{defbar}.

\subsubsection{A priori estimates for the displacement} 
The main result in this section is given in Proposition \ref{prop 2may10}. However, before stating it, we give some auxiliary claims.
First, we introduce the following spaces: 
\begin{eqnarray*} 
V(\Omega^{h}_{f})&:=&\left\{ \vect \varphi\in H^{1}(\Omega_{f}^{h};\RR^3) : \varphi=0\,\text{ on }\,\, \Gamma^{h}\,\text{and}\, \vect \varphi\,\text{ is $\omega$-periodic in}\,(x_1, x_2)\right\}, \\	V(\Omega^{h}_{s})&:=&\left\{\vect \varphi\in H^{1}(\Omega_{s}^{h};\RR^3) :  \varphi \,\text{is $\omega$-periodic in}\,(x_1, x_2) \right\}.
\end{eqnarray*} 	
 We give two helpful lemmas.
\begin{lema} \label{lext} 
There exists $C>0$ such that for every $h>0$ there exists a linear extension operator from $V(\Omega^{h}_{s})$ to $H^1_{\#}(\Omega;\RR^3)$ ($\vect{\varphi} \mapsto \widehat{\vect{\varphi}} $) that satisfies 
\begin{eqnarray*} 
\|\widehat{\vect{\varphi}}\|_{H^1(\Omega;\RR^3)} &\leq & C  \|\vect{\varphi}\|_{H^1(\Omega^h_s;\RR^3)},\\
 \| e_h(\widehat{\vect{\varphi}}) \|_{L^2(\Omega;\RR^{3 \times 3}) } &\leq& C  \| e_h(\vect{\varphi}) \|_{L^2(\Omega_s^{h};\RR^{3 \times 3}) }, \quad \forall \vect \varphi \in V(\Omega^{h}_{s}).    
 \end{eqnarray*} 	
\end{lema} 
\begin{proof} 
See \cite[Chapter 4]{Oleinik}. 
\end{proof} 	
	
\begin{lema}\label{lemma 1}
Let $\vect \varphi\in V(\Omega^{h}_{f})$. Then we have
\begin{equation}\label{eq:56}
\|\vect \varphi\|_{L^2(\Omega^{h}_{f};\RR^3) }\leq C\ee \|\nabla_{h}\vect  \varphi\|_{L^2(\Omega^{h}_{f};\RR^{3 \times 3})} \leq  C\ee \|e_{h}(\vect \varphi)\|_{L^2(\Omega^{h}_{f};\RR^{3 \times 3})}.
\end{equation}
\end{lema}
\begin{proof} 
The estimate can be established on each small cube contained in $\Omega$ by rescaling it on the physical domain.  	
\end{proof} 	
 The following proposition yields important Korn-type estimate with respect to the fluid and solid component of the symmetric gradient.
\begin{propozicija}\label{prop 1}
Let $ \vect \xi\in H^1\left(0,T; H^1(\Omega;\RR^3)\right)$ such that $\vect  \xi (0)=0$ and $\int_{\Omega} \kappa^h \vect{\xi}(t) \,dx=0$, for all $t \in [0,T]$. 
Then the following estimate holds for all $t\in [0,T]$, with a constant $C$ independent of $h$ (\CCC recall \ref{pimapping} \BBB),
\begin{equation}\label{eq:57}
\|\pi_{1/h} \vect  \xi(t)\|_{L^2(\Omega;\RR^3)}\leq\displaystyle C\left\{\frac{1}{h}\|e_{h}(\vect  \xi(t))\|_{L^2(\Omega^{h}_{s};\RR^{3 \times 3})}+\frac{\ee}{h}\int_{0}^{t}\|e_h(\partial_{\tau}\vect  \xi(\tau))\|_{L^2(\Omega^{h}_{f};\RR^{3 \times 3})}\,d\tau\right\}.
\end{equation}
\end{propozicija}

\begin{proof}
For every $t\in[0,T]$, let $ \widehat{\vect \xi}(t)$ be the extension of $\left.\vect \xi(t)\right|_{\Omega^{h}_{s}}$, given by Lemma \ref{lext}. We define $ \vect  z(t)=\widehat{\vect \xi}(t)-\vect \xi(t)$ on $\Omega^{h}_{f}$ and zero elsewhere. Then for every $t\in(0,T)$ we have $\vect  z(t)\in V(\Omega^{h}_{f})$. As a consequence of Lemma \ref{lemma 1} we have
\begin{equation}\label{eq:58}
\|\vect z(t)\|_{L^2(\Omega^{h}_{f};\RR^3) }\leq C\ee \|\nabla_h \vect z(t)\|_{L^2(\Omega^{h}_{f};\RR^{3  \times 3}) }\leq C\ee \|e_h(\vect z(t))\|_{L^2(\Omega^{h}_{f};\RR^{3  \times 3}) },
\end{equation}
for all $t\in(0, T)$. We conclude
\begin{equation} \label{eq:59} 
\|\widehat{\vect \xi}(t)-\vect \xi(t)\|_{L^2(\Omega^{h}_{f};\RR^3) } \leq C\ee\left\{\|e_{h}(\widehat{ \vect \xi}(t))\|_{L^2(\Omega^{h}_{f};\RR^{3  \times 3}) }+\|e_{h}(\vect \xi(t))\|_{L^2(\Omega^{h}_{f};\RR^{3  \times 3}) }\right\}.
\end{equation}
 Next, from (\ref{eq:59}) and using the fact that $\ee, h\leq 1$ and Proposition \ref{propapp1}  we obtain
{\allowdisplaybreaks
\begin{align}
\nonumber &\|\pi_{1/h} \vect \xi(t)\|_{L^2(\Omega;\RR^3)}\leq \|\pi_{1/h}\widehat{\vect \xi}(t)\|_{L^2(\Omega;\RR^3)}+\|\pi_{1/h}(\widehat{\vect \xi}(t)-\vect\xi(t))\|_{L^2(\Omega;\RR^3)}\\*
\nonumber &\leq \|\pi_{1/h}\widehat{ \vect \xi}(t)\|_{L^2(\Omega;\RR^3)}+\frac{1}{h}\|\widehat{\vect \xi}(t)-\vect \xi(t)\|_{L^2(\Omega;\RR^3)}\\
\nonumber &\leq \|\pi_{1/h}\widehat{\vect \xi}(t)\|_{L^2(\Omega;\RR^3)}+\frac{C\ee}{h}\left\{\|e_{h}(\widehat{\vect \xi}(t))\|_{L^2(\Omega^{h}_{f};\RR^{3  \times 3}) }+\|e_{h}(\vect \xi(t))\|_{L^2(\Omega^{h}_{f};\RR^{3  \times 3}) }\right\}\\
\nonumber &\leq
C \Bigg\{\frac{1}{h} \|e_h(\widehat{\vect  \xi}(t))\|_{L^2(\Omega;\RR^{3 \times 3})}+\left|\int_{\Omega} \pi_{1/h}(\kappa^h \widehat{\vect \xi}(t))\right|+ \frac{\ee}{h}\|e_{h}(\widehat{\vect \xi}(t))\|_{L^2(\Omega^{h}_{f};\RR^{3  \times 3}) } \\* 
\label{eqnak3} &\hspace{+5ex}+\frac{\ee}{h}\|e_{h}(\vect \xi(t))\|_{L^2(\Omega^{h}_{f};\RR^{3  \times 3}) }  \Bigg\}
\end{align}
}
Note that since $\int_{\Omega}\kappa^h \vect \xi(t)\,dx=0$,  we have that 
\begin{equation} \label{eqnak1}  
\int_{\Omega} \pi_{1/h} (\kappa^h \widehat{\vect \xi} (t)) )\,dx=\int_{\Omega} \pi_{1/h} (\kappa^h\widehat{\vect \xi} (t)-\kappa^h\vect \xi(t)) )\,dx.   
\end{equation} 
From \eqref{eq:59} and \eqref{eqnak1} we conclude 
\begin{equation} \label{eqnak2} 
\left|\int_{\Omega} \pi_{1/h} (\kappa^h\widehat{\vect \xi}(t)))\,dx\right| \leq \frac{C\ee}{h} \left\{\|e_{h}(\widehat{ \vect \xi}(t))\|_{L^2(\Omega^{h}_{f};\RR^{3  \times 3}) }+\|e_{h}(\vect \xi(t))\|_{L^2(\Omega^{h}_{f};\RR^{3  \times 3}) }\right\}. 
\end{equation} 
Next, we remark that using $\vect{\xi} (0)=0$, we have
\begin{equation} \label{eq:61pedro}
\|e_h(\vect \xi(t))\|_{L^2(\Omega^{h}_{f};\RR^{3 \times 3}) } =\left\|\int_{0}^{t}e_{h}\left(\partial_{\tau}\vect \xi(t)\right)d\tau\right\|_{L^2(\Omega^{h}_{f};\RR^{3 \times 3})}\leq \displaystyle\int_{0}^{t}\left\|e_{h}\left(\partial_{\tau}\vect \xi(t)\right)\right\|_{L^2(\Omega^{h}_{f};\RR^{3 \times 3})}d\tau.
\end{equation}

\eqref{eq:57} now follows from \eqref{eqnak3}, \eqref{eqnak2} and \eqref{eq:61pedro}  by using Lemma \ref{lext}, and $\ee \leq 1$.
\end{proof}
The following proposition gives us necessary \emph{a priori} estimates.  
\begin{propozicija}\label{prop 2may10}
Let us suppose that
\footnote{this can always be achieved by translation of coordinate system for every $t$. \CCC This is a standard argument: If \eqref{kirill3} is not satisfied and if we denote by $C(t)=\int_{\Omega} \vect{F}^h(t)\, dx$, then it is easy to see  that the solution $\breve{\vect{u}}^h(t)= \vect{u}^h(t)-C_1(t)$, $\breve{p}^h(t)=p^h(t)$ corresponds to the loads $\vect{F}^h(t)-C(t)$ which obviously satisfy \eqref{kirill3}. Here $C_1:[0,T] \to \R^3$ solves
$$C_1''(t)\eta \int_{\Omega} \kappa^h\, dx=C(t), \ C_1(0)=0, \ C_1'(0)=0,  $$ which are the equations Newton's second law. One can then make the appropriate conclusions for $(\vect{u}^h,p^h)$ by analyzing $(\breve{\vect{u}}^h,\breve{p}^h)$.  \BBB } 
\begin{equation} \label{kirill3} 
	\int_{\Omega} \vect{F}^h(t)\,dx=0, \quad \textrm{ for every } t \in [0,T]. 
\end{equation} 
\begin{enumerate} 
\item In the case when $\eta=\eta(h)$ is bounded from below by  a positive constant we assume
\begin{equation} \label{forcesassumptions0}
\|F^h_3\|_{L^2(0,T;L^2(\Omega))}+h\sum_{\alpha=1,2}\| {F}^h_{\alpha}\|_{H^1(0,T; L^2(\Omega))}
\leq C,
\end{equation}
where $C$ doesn't depend on $h$.
\item In the case when  $\eta(h) \to 0$, we assume that 
\begin{equation} \label{forcesassumptions1}
\|\pi_h \vect{F}^h\|_{H^{1}(0,T; L^2(\Omega;\RR^3))}\leq C,
\end{equation} 
where $C$ doesn't depend on $h$. 
\end{enumerate} 
 Then we have:
\begin{equation}
\label{eq:58*may10}
    \eta^{\frac{1}{2}} \left\|\partial_t\vect{u}^h\right\|_{L^{\infty}(0,T; L^2(\Omega;\RR^3))}+\displaystyle\frac{1}{h}\left\|e_h(\vect{u}^h)\right\|_{L^{\infty}(0,T; L^2(\Omega^{h}_{s};\RR^{3 \times 3}))}+\frac{\ee}{h^2} \left\|e_h\left(\partial_t\vect{u}^h\right)\right\|_{L^{2}(0,T; L^2(\Omega^{h}_{f};\RR^{3  \times 3}))}\leq C,
\end{equation}
where $C$ doesn't depend on $h$. Here $(\vect{u}^h,p^h)$ is a solution of \eqref{varformporous2gamma0may10} with the initial condition \eqref{eq:51rescaledver}. 
\end{propozicija}

\begin{proof}
As a consequence of \eqref{kirill3}, we have $\int_{\Omega} \kappa^h \vect{u}^h(t)\, dx=0$, for all $t \in [0,T]$. 
From \eqref{forcesassumptions0} we conclude
\begin{equation} \label{nak1010} 
h\sum_{\alpha=1,2}\|F^h_{\alpha}\|_{L^{\infty}(0,T,L^2(\Omega;\RR^3))} \leq C,
\end{equation}
and from \eqref{forcesassumptions1} we can conclude 
\begin{equation} \label{nak10} 
\|\pi_h\vect{F}^h\|_{L^{\infty}(0,T,L^2(\Omega;\RR^3))} \leq C,
\end{equation} 
with $C$ independent of $h$. Firstly we assume that $\vect{F}^h \in H^1(0,T; L^2(\Omega;\RR^3))$. In this case we can use Remark \ref{obsjosip1} and
we take $\varphi=\partial_{t}\vect{u}^{h}$ as test function in \eqref{varformporous2gamma0may1000}. This yields for almost every $t \in (0,T)$
\begin{align}
&\displaystyle\frac{1}{2}\frac{d}{dt}\left(\displaystyle\int_{\Omega}\eta \kappa^{h}\nonumber|\partial_{t}\vect{u}^h(t)|^2\,dx+\frac{1}{h^{2}}\int_{\Omega_{s}^{h}}\mathbb{A}^h e_h(\vect{u}^h(t)):e_h(\vect{u}^h(t))\,dx\right)\\
&\label{eq:61may10}+\displaystyle\frac{\ee^{2}}{h^{4}}\int_{\Omega^h_{f}}2 |e_h(\partial_{t}\vect{u}^h(t))|^2\,dx=\int_{\Omega}\vect{F}^h(t)\partial_{t}\vect{u}^h(t)\,dx.
\end{align}
By integration over $[0,t]$ for $t \in [0,T]$, we obtain
\begin{equation}\label{eq:62may10}
    \begin{split}
        &\displaystyle\frac{1}{2}\displaystyle\int_{\Omega}\eta \kappa^{h}|\partial_{t}\vect{u}^h(t)|^2\,dx+\frac{1}{h^{2}}\int_{\Omega_{s}^{h}}\mathbb{A}^h e_h(\vect{u}^h(t)):e_h(\vect{u}^h(t))\,dx \\
&+\displaystyle\frac{\ee^{2}}{h^{4}}\int_{0}^{t}\int_{\Omega^h_{f}}2 |e_h(\partial_{t}\vect{u}^h(\tau))|^2\,dx\,d\tau=\int_{0}^{t}\int_{\Omega}\vect{F}^h(\tau)\partial_{t}\vect{u}^h(\tau)\,dx\,d\tau.
    \end{split}
\end{equation}
On the other hand, we have
{\allowdisplaybreaks
\begin{equation} 
    \label{josipispravak}
    \begin{split}       &\left|\displaystyle\int_{0}^{t}\int_{\Omega}\vect{F}^h(\tau)\partial_{t}\vect{u}^h(\tau)\,dx\,d\tau\right|\\*
	&\leq \sum_{\alpha=1,2} \left( \left|\displaystyle\int_{\Omega}F_{\alpha}^h(t){u}^h_{\alpha}(t)\,dx\right|+\left|\displaystyle\int_{0}^{t}\int_{\Omega}\partial_{t }F^h_{\alpha}(\tau){u}^h_{\alpha}(\tau)\,dx\,d\tau\right|\right)+C\|F^h_3\|_{L^2(\Omega)}\|\partial_t u_3^h(t)\|_{L^2(\Omega)}\\*
	&\leq C\Bigg(\sum_{\alpha=1,2}\left(\|hF^h_{\alpha}(t)\|_{L^2(\Omega)}\|\frac{u_{\alpha}^h(t)}{h}\|_{L^2(\Omega)}+\|h\partial_tF^h_{\alpha}\|_{L^{2}(0, t; L^2(\Omega))}\|\frac{u_{\alpha}^h}{h}\|_{L^2(0,t; L^2(\Omega;\RR^3))}\right)\\ 
    & \hspace{+5ex} +\|F^h_3(t)\|_{L^2(\Omega)}\|\partial_t u_3^h(t)\|_{L^2(\Omega)}\Bigg).
    \end{split}  
\end{equation}
}
\CCC In case when $\eta=\eta(h)$ is bounded from below the estimate 
\eqref{eq:58*may10} follows from  \eqref{nak1010}, \eqref{nak10}, \eqref{eq:62may10}, \eqref{josipispravak} and \eqref{eq:63may10}  by using Young's inequality, Gronwall's lemma and Proposition \ref{prop 1}.
Note that, after doing Young's inequality in the last term of \eqref{josipispravak} in the way: 
\begin{equation} 
\|F^h_3\|_{L^2(\Omega)}\|\partial_t u_3^h(t)\|_{L^2(\Omega)} \leq \frac{1}{\alpha} \|F_3^h(t)\|^2_{L^2(\Omega)}+\alpha\|\partial_t u_3^h(t)\|^2, 
\end{equation} 	
for $\alpha>0$ small enough, the last term can be absorbed by the left hand side of \eqref{eq:62may10}. This is not possible in the case when $\eta=\eta(h)\to 0$, 
when we do the following estimate: \BBB 
\begin{equation}
\label{eq:63may10}
    \begin{split}   &\left|\displaystyle\int_{0}^{t}\int_{\Omega}\vect{F}^h(\tau)\partial_{t}\vect{u}^h(\tau)\,dx\,d\tau\right|\leq \left|\displaystyle\int_{\Omega}\vect{F}^h(t)\vect{u}^h(t)\,dx\right|+\left|\displaystyle\int_{0}^{t}\int_{\Omega}\partial_{t }\vect{F}^h\vect{u}^h(\tau)\,dx\,d\tau\right|\\
&\leq  C\left(\|\pi_{h}\vect{F}^h(t)\|_{L^2(\Omega;\RR^3)}\|\pi_{1/h}\vect{u}^h(t)\|_{L^2(\Omega;\RR^3)}+\|\pi_{h}\partial_t\vect{F}^h\|_{L^{2}(0, t; L^2(\Omega;\RR^3))}\|\pi_{1/h}\vect{u}^h\|_{L^2(0,t; L^2(\Omega;\RR^3))}\right). 
    \end{split}
\end{equation}
Again by using \eqref{nak1010}, \eqref{nak10}, \eqref{eq:62may10}, \eqref{josipispravak} and \eqref{eq:63may10}, the estimate
 \eqref{eq:58*may10} also follows in this case by using Young's inequality, Gronwall's lemma and Proposition \ref{prop 1}.
 
\eqref{eq:58*may10} can also be proved in the case when $\eta(h)$ doesn't converge to zero and $F_3^h \notin H^1(0,T; L^2(\Omega))$. This can be done by approximation of the loads using stability estimate of Proposition \ref{newpropmay10} 1. \CCC Although the constant in the stability estimate depends on $h$, the constant on the right hand side of \eqref{eq:58*may10} depends on $L^2$ norm of $F_3^h$ \footnote{\CCC Actually from the proof it follows that this constant is of the form $C_1  C$ where $C_1$ doesn't depend on $h$ or $\vect{F}^h$ and $C$ is given by \eqref{forcesassumptions0}.\BBB}  and the approximation is done for every fixed $h>0$. More precisely, we take $(F_3^{h,\eps})_{\eps>0} \subset H^1(\Omega)$ such that $F_3^{h,\eps} \to F_3^h$ in $L^2$  and then let $\eps \to 0$ in \eqref{eq:58*may10},  for fixed $h$.\BBB  
\end{proof} 
\CCC \begin{remarkica} 
It is standard in the derivation of the plate theory in the context of linearized elasticity that one needs to scale differently in-plane and vertical components of the loads to obtain the limit equations. This is connected with the different scaling of in-plane and vertical components  of the displacement (in the {\it a priori} estimate) and the physically observed fact that for the plate it is much easier to bend than to stretch (see also \cite{ciarlet2,Marinvelciczubrinic2022}).   
\end{remarkica} \BBB	
\subsubsection{A priori estimate for pressure}

With  Proposition \ref{prop 2may10} at hand it is now possible for us to establish the\emph{ a priori }estimate for a pressure. We will need the following proposition \CCC and corollary. \BBB
\begin{propozicija}\label{rescaldivpress}
For every  $g\in L^2(\Omega)$ and every $h$, there exists $\vect{v}^h\in H^1_{\#}(\Omega;\mathbb{R}^3)$, such that (\CCC recall \eqref{defscaledg}) \BBB
\begin{equation}\label{rev1} 
\diver_h\vect{v}^h=g \quad\textrm{ and }\quad \|\vect{v}^h\|_{L^2(\Omega;\RR^3)}+\|\nabla_h\vect{v}^h\|_{L^2(\Omega;\RR^3)}\leq C\|g\|_{L^2(\Omega)},
\end{equation}
where $C$ is independent of $h$. 
\end{propozicija}

\begin{proof}
\CCC We follow the idea of \cite{braides} in a different context. 
In order to get the control of scaled gradients, \BBB we divide $\Omega$ in small plates $\{\Omega^h_i\}_{i=1}^\infty\subset \Omega$ of size $h$, such that $\Omega^h_i\cap \Omega^h_j=\emptyset$ for all $i\neq j$. Here 
\begin{equation}
\Omega^h_k:=\omega\times \left(-h/2+kh, h/2+kh\right) \cap \Omega, \quad k \in \mathbb{Z}. 
\end{equation}
Notice that for every $h$ only finitely number of $\Omega_i^h$ are non-empty, i.e. for every $h>0$ there exists $n(h)$ such that $\Omega_k^h \neq \emptyset$ for $|k| \leq n(h)$ and  $\Omega_k^h = \emptyset$, for $|k|>n(h)$. \CCC Note that for $|k|<n(h)$ we have that 
$$ \omega \times (-h/2+kh, h/2+kh) \subset \Omega.   $$ \BBB
We define $g^h$ on $\Omega$, by defining it in the following way
$$
g^h(x_1^h,x_2^h,x_3^h)=g(x_1^h,x_2^h,x_3^h/h), \quad \textrm{if } (x_1^h,x_2^h,x_3^h) \in \Omega_0^h, $$
and we extend it by periodicity in $x_3$ to whole $\Omega$. We find $\vect{k}^h  \in H^1_{\#}(\Omega;\RR^3) $,
 which satisfies
 \begin{equation} \label{nak20} 
 \diver \vect{k}^h =g^h,\quad \|\vect{k}^h \|_{L^2(\Omega;\RR^3)}+\|\nabla\vect {k}^h\|_{L^2(\Omega;\RR^{3\times 3})}\leq C \|g^h\|_{L^2(\Omega)}
 \end{equation} 
 This can be obtained by solving $\Delta  \phi^h=g^h$ on $\Omega$ with zero boundary condition on $\omega \times \{-\frac{1}{2},\frac{1}{2} \}$ (and periodic on $\partial \omega \times I$) and putting $\vect{k}^h =\nabla \phi^h$. 
From \eqref{nak20} it follows that there exists $\Omega_k^h$, for $|k|<n(h)$, such that for $\widetilde{\vect{v}}^h=\vect{k}^h |_{\Omega_k^h}$, it is satisfied 
$${\rm div }\,     \widetilde{\vect{v}}^h= g^h|_{\Omega_k^h}, \quad  \|\widetilde{\vect{v}}^h\|_{L^2(\Omega_k^h;\RR^3)}+\|\nabla\widetilde{\vect{v}}^h\|_{L^2(\Omega_k^h;\RR^{3\times 3})}\leq C \|g^h\|_{L^2(\Omega_k^h)},$$ 
where $C$ is independent of $h$ (\CCC with constant $C$ posssibly larger than in \eqref{nak20}). \BBB $\vect v^h$ can then be defined by translating $\widetilde{\vect v}^h$ on $\Omega_0^h$  and rescaling it to whole $\Omega$, \CCC $\vect{v}^h(x):=\widetilde{\vect{v}}^h(\hat{x},hx_3+kh)$. Then we have that \eqref{rev1} is satisfied, since $\nabla_h \vect{v}^h (x)=\nabla \widetilde{\vect{v}}^h(\hat{x},hx_3+kh)$.  \BBB 
\end{proof}
\CCC
The following corollary can be proved in the same way as Proposition \ref{rescaldivpress}. 
Before stating it, we introduce  the space  \begin{equation} 
	H_T^1L^2:=\left\{\varphi \in H^1(0,T;L^2(\Omega)) : \varphi(T)=0\right\}.
\end{equation}	\BBB
\CCC
\begin{corollary}\label{correscaldivpress}
	For every  $g\in H_T^1L^2$ and every $h$, there exists $\vect{v}^h\in H^1(0,T;H^1_{\#}(\Omega;\mathbb{R}^3))$, such that $\vect{v}(T)=0$ and 
	\begin{equation}\label{rev1} 
		\diver_h\vect{v}^h(t)=g(t), \ \forall t \in [0,T], \quad\textrm{ and }\quad \|\vect{v}^h\|_{H^1(0,T;L^2(\Omega;\RR^3))}+\|\nabla_h\vect{v}^h\|_{H^1(0,T;L^2(\Omega;\RR^3))}\leq C\|g\|_{H_T^1L^2},
	\end{equation}
	where $C$ is independent of $h$. 
\end{corollary} 
\begin{proof} 
	We define $\Omega_k^h$, $g^h$, $\vect{k}^h $ and $\phi^h$ in the same way as in Proposition \ref{rescaldivpress}, for every $t \in [0,T]$. Note that $\Delta \partial_t \phi^h= \partial_t g^h$ and $\partial_t \vect{k}^h=\partial_t \phi^h$.  
	\eqref{nak20} becomes
	\begin{equation} \label{naknak20} 
		\diver \vect{k}^h(t) =g^h(t), \ \forall t \in [0,T], \quad \|\vect{k}^h \|_{H^1(0,T;L^2(\Omega;\RR^3))}+\|\nabla\vect {k}^h\|_{H^(0,T;L^2(\Omega;\RR^{3\times 3})}\leq C \|g^h\|_{H_T^1L^2}.
	\end{equation}
	Again it follows that there exists $\Omega_k^h$, for $|k|<n(h)$, such that for $\widetilde{\vect{v}}^h=\vect{k}^h |_{\Omega_k^h}$, it is satisfied 
	$${\rm div }\,    \widetilde{\vect{v}}^h(t)= g^h(t)|_{\Omega_k^h},  \ \forall t\in [0,T], \  \|\widetilde{\vect{v}}^h\|_{H^1(0,T;L^2(\Omega_k^h;\RR^3))}+\|\nabla\widetilde{\vect{v}}^h\|_{H^1(0,T;L^2(\Omega_k^h;\RR^{3\times 3}))}\leq C \|g^h\|_{H^1(0,T;L^2(\Omega_k^h))}.$$
The rest of the proof follows the argument of Proposition \ref{rescaldivpress}.
\end{proof} 	
\BBB 
Next proposition establishes \emph{a priori} estimate for the pressure. 

\begin{propozicija}\label{pressureregularity}
Let conditions of Proposition \ref{prop 2may10} be satisfied.
Let  $(\vect{u}^h,p^h)$ be a solution of \eqref{varformporous2gamma0may10} with the initial condition \eqref{eq:51rescaledver}. 
The following is satisfied: 
\begin{equation} 
\|\widehat{p}^h\|_{(H_T^1L^2)'} \leq Ch,
\end{equation}  
where $C$ is independent of $h$ and $\widehat{p}^h$ is the extension by zero of $p^h$ onto the whole $\Omega$.
\end{propozicija} 
\begin{proof} 

Using \CCC Corollary \BBB \ref{correscaldivpress} we take for \CCC $g \in H_T^1L^2$ \BBB and $h>0$,  \CCC $\vect v^h \in H^1(0,T;H^1_{\#}(\Omega;\mathbb{R}^3))$ that satisfies \CCC $\vect{v}^h(T)=0$ and
$$ {\rm div}_h\, \vect{v}^h(t)= g(t), \forall t \in [0,T],\quad \|\vect{v}^h\|_{H^1(0,T;L^2(\Omega;\RR^3))}+\|\nabla_h\vect{v}^h\|_{H^1(0,T;L^2(\Omega;\RR^3))}\leq C\|g\|_{H_T^1L^2}. $$
\BBB
By \CCC taking $\vect{v}^h$ as a test function in \BBB \eqref{varformporous2gamma0may10} and using Proposition \ref{prop 2may10} we conclude that 
$$\CCC _{(H_T^1L^2)'}\big\langle \frac{\widehat{p}^h}{h}, g \big\rangle_{H_T^1L^2}\BBB \leq C \CCC \|g\|_{H_T^1 L^2}, \BBB $$
where $C>0$ is independent of $h$, which finishes the proof. 
\end{proof}

\section{Quasi-static case} 
\label{secetanula} 
In this section,  we  analyze the case when $\lim_{h \to 0} \eta (h)=0$.  In Section \ref{compactness}  we provide the compactness result for the sequence of solutions of \eqref{varformporous2gamma0may10}. 
\CCC In order to deal with the limit problem, we have to use the rescaled two-scale convergence, which is defined in Appendix \ref{secaptwoscale}. The reason is that, after rescaling on the canonical domain, the characteristic cell has size $\eps$ in in-plane direction and $\eps/h$ in the vertical direction. Thus one needs to use the test functions that oscillate with different period in in-plane and vertical direction.  Auxiliary compactness claims about rescaled two-scale convergence are also given in Appendix \ref{secaptwoscale}.      \BBB
In Section \ref{effective} we define the effective tensors appearing in the limit problem and prove some of their properties. In Section \ref{limit} we  obtain the limit model, while in Section \ref{analysis} we prove the existence and uniqueness result and energy-type equality for the limit problem.  
\CCC For the definition of the weak solution, existence and uniqueness result we use the results from \cite{gurvich}. Here, however, one needs to put an additional effort to define the appropriate operators, since the limit equations do not decouple, see the proof of Theorem \ref{existeo} below.  \BBB
In Section \ref{secstrongcon} we prove the strong convergence of the solutions of \eqref{varformporous2gamma0may10} to the solution of the limit problem  with appropriate correctors, while in Section \ref{secgen} we discuss the possible generalization of Assumptions \ref{assumption on regions}, as well as the possibility of having surface loads, non-zero initial conditions, non-periodic boundary conditions at the transverse boundary 
and the situation when fluid part touches the upper and lower boundary of $\Omega^h$. \CCC We recall that $\bar{\kappa}= \int_I \kappa \, dx_3$ (see \ref{defbar}) and $\kappa$ is the weak star limit of $\kappa^h$ as $h \to 0$, consequently a function on $\omega$, bounded from below and above by a positive constant. \BBB

We additionally introduce the following notation.
We denote by 
$$L_{\kappa,0}:=\left\{(\vect{\psi}_1,\psi_2) \in H^1_{\#}(\omega;\mathbb{R}^2) \times H^2_{\#}(\omega) : \int_{\Omega} \kappa(\vect{\psi}_1-x_3 \nabla \psi_2) \, dx=0, \, \int_{\omega} \bar{\kappa}\psi_2 \, d\hat{x}=0\right\},$$
and by $M=L^2(\omega; H^1(J_K)\oplus L^2(J_p\backslash J_K))$.
\CCC Note that 
$$L_{1,0}=\dot{H}^1_{\#}(\omega;\R^2) \times \dot{H}^2_{\#}(\omega).    $$
\BBB

 For $\vect{w} \in H^1(\mathcal{Y};\RR^3)$, $\vect{g} \in \RR^3$, we define \CCC $\mathfrak{C}_{\infty}(\vect{w}, \vect{g}) \in L^2(Y;\R^{3 \times 3}_{\sym}) $ in the following way: \BBB
\begin{equation} \label{revdef1} 
\mathfrak{C}_{\infty}(\vect{w}, \vect{g}):=e_{y}(\vect{w})+\text{sym}\left(0|0|\vect{g}\right). \end{equation} 
The definition naturally extends to $H^1(\mathcal{Y};\RR^3) \times \RR^3$ valued functions, \CCC for e.g. if $(\vect{w},\vect{g}) \in L^{\infty}(0,T;L^2(\Omega;$ $H^1(\mathcal{Y};\RR^3)\times \RR^3 ))$, then taking \eqref{revdef1} pointwise in $(x,t)$, we obtain   
$\mathfrak{C}_{\infty}(\vect{w}, \vect{g}) \in L^{\infty}(0,T;L^2(\Omega;L^2(Y;\R^{3 \times 3}_{\sym})))$. \BBB 
\label{secquasi-static}
\subsection{Compactness result} \label{compactness} 
\CCC When dealing with evolution problems, if we want to use the results from elliptic (static) problems, we often do averaging in time. \BBB 
For a given $\nu \in L^1(0,T)$ and a Hilbert space $X$  we define the operator $\cdot^{\nu}: L^{\infty}(0,T;X)\to L^{\infty} (X)$ in the following way
$$\vect{u}^{\nu}: =\int_0^T \vect{u}(t) \nu(t)\,dt. $$ 
Furthermore, we define:
\begin{equation} \label{defH1H1gggg} 
	H_T^1L^2(\Omega \times \mathcal{Y}):=\{\varphi \in H^1(0,T;L^2(\Omega \times \mathcal{Y})): \varphi(T)=0\}.
\end{equation}	
The following theorem gives us the compactness result. 
Appropriate definitions of two-scale convergences are given in Section \ref{secaptwoscale}. Recall also the definition of the spaces $\Omega \times \mathcal{Y}_{f/s}^{x_3}$ given in \eqref{yxspace}.

 \begin{teorem}\label{gammazeromay11}
 Let assumptions \eqref{kirill3} and  \eqref{forcesassumptions1} and 
 Assumption \ref{assumption on regions} be satisfied. 
Assume that \footnote{For the definitions of $\xrightharpoonup{t,2-r,p}$, $\xrightharpoonup{2-r}$ see Appendix \ref{secaptwoscale}.} $\pi_h\vect{F}^{h}\xrightharpoonup{t,2-r,2} \vect{F},$
where $\vect{F} \in H^1(0,T; L^2(\Omega \times \mathcal{Y};\mathbb{R}^3))$. 
The following statements hold: Let $(\vect{u}^h,p^h)$ be the solution of (\ref{varformporous2gamma0may10}) with initial condition \eqref{eq:51rescaledver}. Then there exist limits (on a subsequence)
\footnote{\CCC  The identities $$ p=0 \textrm{ on } \Omega \times \mathcal{Y}_{s}^{x_3}, \,  p=p(x,t)\chi_{\Omega \times \mathcal{Y}_{f}^{x_3}}$$ should be understood in the weak sense. Namely, they imply: 
$$_{(H_T^1L^2(\Omega \times \mathcal{Y}))'}\langle p, \varphi\rangle_{H_T^1L^2(\Omega \times \mathcal{Y})}=
_{(H_T^1L^2(\Omega \times \mathcal{Y}))'}\langle p, \int_{\mathcal{Y}_f(x_3)}\varphi(t,x,\cdot) \, dy\rangle_{H_T^1L^2(\Omega \times \mathcal{Y})}
, \quad \forall \varphi \in H_T^1L^2(\Omega \times \mathcal{Y}). $$	  \BBB} : 
\begin{align}
\label{nakk101}(\mathfrak{a},\mathfrak{b}) &\in L^{\infty}(0,T; L_{\kappa,0}),\\ \label{rev1111} 
 \vect{w}&\in L^{\infty} (0, T; L^2(\Omega, \dot{H}^1(\mathcal{Y};\RR^3))) ,\\ \label{kirill10}
 \vect{g}&\in L^{\infty}(0,T; L^2(\Omega;\RR^3)),\\  \label{kreso1} \CCC 
 \vect{u}^0_f & \CCC \in H^1(0,T; L^2(\Omega; H^1(\mathcal{Y};\RR^3))),\ \diver_{y} \vect{u}_{f}^{0}=0,\  \vect{u}_f^0 =0 \textrm{ on } [0,T] \times \Omega \times \mathcal{Y}_s^{x_3},\,\vect{u}_{f}^{0}(0)=0, \BBB\\
  \label{nakk106}  p&\in(H_T^1L^2(\Omega \times \mathcal{Y}))',\ p=0 \textrm{ on } \Omega \times \mathcal{Y}_{s}^{x_3}, \,  p=p(x,t)\chi_{\Omega \times \mathcal{Y}_{f}^{x_3}}
        \end{align} 
 such that for the sequence of the solutions $(\vect u^h, p^h)$ to \eqref{varformporous2gamma0may10},  we have (on a subsequence) \footnote{\CCC $p^{\nu} \in L^2(\Omega \times \mathcal{Y})$  is defined in the following way: 
 	$$\langle p^{\nu},\varphi\rangle_{L^2(\Omega \times \mathcal{Y})}=_{(H_T^1L^2(\Omega \times \mathcal{Y}))'}\langle p, \nu\varphi \rangle_{H_T^1L^2(\Omega \times \mathcal{Y})}, \quad \forall \varphi \in L^2(\Omega \times \mathcal{Y}) .$$
 	\CCC Similarly  $(\widehat{p}^h{)^\nu} \in L^2(\Omega)$ is defined in the following way: 
 	$$ \langle (\widehat{p}^h)^{\nu}, \varphi\rangle_{L^2(\Omega)} := _{(H_T^1L^2)'}\langle p, \nu\varphi \rangle_{H_T^1L^2}, \quad \forall \varphi \in L^2(\Omega).  $$ 
 	\BBB
 }
\begin{align}
&\label{crgamma0}(h^{-1}u^{h}_{\alpha})^{\nu}\overset{L^2}{\rightarrow} (\mathfrak{a}_{\alpha}-x_3\partial_{\alpha}\mathfrak{b})^{\nu},\quad (u^{h}_{3})^{\nu}\overset{L^2}{\rightarrow} (\mathfrak{b})^{\nu},\quad \forall \nu \in L^1(0,T)\\
&\label{cr1gamma000} h^{-2}\vect{u}^h_f\xrightharpoonup{t,2-r,2} \vect{u}^0_f,\\
&\label{cr2gamma0}h^{-1}e_h (\widehat{\vect{u}}^{h})\xrightharpoonup{t,2-r,2} \iota (e_{\widehat{x}}(\mathfrak{a})-x_3\nabla^2_{\widehat{x}}\mathfrak{b})+\mathfrak{C}_{\infty}(\vect{w}, \vect{g}),\\
&\label{cr3gamma0}\ee h^{-2} e_h ( \vect{u}^h_f)\xrightharpoonup{t,2-r,2} e_y(\vect{u}^0_f),\\
&\label{cr5gamma00}\ee h^{-2} e_h ( \partial_t \vect{u}^h )\CCC \chi_{\Omega^h_f} \BBB\xrightharpoonup{t,2-r,2} e_y(\partial_t\vect{u}^0_f) \\
&\label{cr5gamma0}h^{-1}(\widehat{p}^{h})^{\nu}\xrightharpoonup{2-r} p^{\nu},\quad \forall \nu \in H^1(0,T), \textrm{ such that } \nu(T)=0. 
\end{align}
Here $\widehat{p}^h$ is the extension of $p^h$ by zero on $\Omega$ and 
$ \vect{u}^h=\widehat{\vect{u}}^h+\vect{u}_f^h;$ 
$\widehat{\vect{u}}^h$
is the extension defined in Lemma \ref{lext} \CCC and $\vect{u}_f^h$ is simply the difference $\vect{u}_f^h:= \vect{u}^h-\widehat{\vect{u}}^h$.\BBB \footnote{\CCC Note that $\vect{u}^h_f$ is zero on the solid part. \BBB}
Moreover,  the following identity is valid for every $\varphi \in L^2(\omega;H^1(J_p))$ for a.e. $t \in [0,T]$:
\begin{equation}
\label{cr4gamma0}
    \begin{split}
        &\displaystyle -\int_{\Omega_p} |\mathcal{Y}_{f}(x_3)|\diver_{\widehat{x}}\,\mathfrak{a}(\widehat{x},t)\varphi(x)\,dx+\int_{\Omega_p}|\mathcal{Y}_{f}(x_3)|x_3\diver_{\widehat{x}} \,\nabla_{\widehat{x}}\mathfrak{b}(\widehat{x},t) \varphi(x)\,dx\\
&\quad+ \displaystyle\displaystyle \int_{\Omega_p}\int_{\mathcal{Y}_{f}(x_3)} u^{0}_{f,3}(x,y,t)dy \partial_3 \varphi(x)\,dx-\displaystyle\int_{\Omega_p} \int_{\mathcal{Y}_{f}(x_3)}\left[\trace\,\mathfrak{C}_{\infty}\left(\vect{w},\vect{g}\right)(x,y,t)\right]\,dy\varphi(x) \,dx=0.
    \end{split}
\end{equation}
 
\end{teorem}

 \begin{proof}
 	\CCC The first part of the proof uses Griso's decomposition and its consequences from Appendix \ref{griso}. Griso's decomposition enables us to obtain the two-scale limits of sequences with bounded symmetrized scaled gradients. These results are fundamental for dimension reduction in linearized elasticity and they have been used also in \cite{Vel14a,mbukalvelcic2017,Marinvelciczubrinic2022}to obtain compactness result. Here one needs to adapt them to deal with fluid part and for the case of oscillations of the material across the thickness. \BBB

As in the proof of Proposition \ref{prop 2may10} we have that $\int_{\Omega} \kappa^h \vect{u}^h(t)=0$, for every $t \in [0,T]$. 
\CCC We have from \BBB Proposition \ref{prop 2may10} 
\begin{equation}\label{eq:77march21}
\displaystyle\frac{1}{h}\left\|e_h\left(\vect{u}^h\right)\right\|_{L^{\infty}(0,T;L^2(\Omega^{h}_{s};\RR^{3  \times 3})) }\leq C.
\end{equation}
Thus, using the extension operator properties from Lemma \ref{lext}, we have:
\begin{equation}\label{eq:78march21}
\displaystyle\left\|e_h\left(\frac{1}{h}\widehat{\vect{u}}^{h}\right)\right\|_{L^{\infty}(0,T;L^2(\Omega;\RR^{3  \times 3}) }\leq C\displaystyle\frac{1}{h}\left\|e_h\left(\vect{u}^h\right)\right\|_{L^{\infty}(0,T;L^2(\Omega;\RR^{3  \times 3}) }\leq C.
\end{equation}
We obtain that on a subsequence 
\begin{equation} \label{kirill1} 
e_h\left(\frac{1}{h}\widehat{\vect{u}}^{h}\right) \xrightharpoonup{t,2-r,2} \mathbb{L}, 
\end{equation} 
where $\mathbb{L} \in L^\infty(0,T; L^2(\Omega \times \mathcal{Y};\RR^{3 \times 3}_{\rm sym}))$. Next we take an arbitrary $\nu \in L^1(0,T)$. Lemma \ref{lemmaA8} and Remark \ref{remnak1} yield  the following decomposition of the sequence $(\widehat{\vect{u}}^{h})^{\nu}$:
\begin{equation}
\label{convergenceextensionmarch21}
    \begin{split}
        \frac{1}{h}(\widehat{\vect{u}}^{h})^{\nu}&= (
      -x_3\partial_1\mathfrak{b}(\nu) ,
       -x_3\partial_2\mathfrak{b}(\nu) ,
       h^{-1}\mathfrak{b}(\nu)
   )^T+
      (\mathfrak{a}_1 (\nu) ,
       \mathfrak{a}_2(\nu) ,
       0
   )^T+\vect{\psi}^h(\nu)+C^h(\nu),\\  \displaystyle\frac{1}{h}e_{h}\left((\widehat{\vect{u}}^{h})^{\nu}\right)&=\iota \left(-x_3\nabla^2_{\widehat{x}}\mathfrak{b}(\nu)+e_{\widehat{x}}(\mathfrak{a}(\nu))\right)+e_h\left(\vect{\psi}^h(\nu)\right),
    \end{split}
\end{equation}
where $\mathfrak{b}(\nu)\in H^2_{\#}(\omega), \mathfrak{a}(\nu)\in  H^1_{\#}(\omega;\RR^2), (\vect{\psi}^h(\nu))_{h>0}\subset  H^1_{\#}(\Omega;\RR^3)$,
$
 h\pi_{1/h}\vect{\psi}^h(\nu)\overset{L^2}{\rightarrow}0$
  and  $C^h(\nu) \in \mathbb{R}^3$ is chosen such that
\begin{equation} \label{kirill6}
\int_{\Omega} \kappa^h \left(\frac{1}{h}(\widehat{\vect{u}}^h)^{\nu}-C^h(\nu)\right)\, dx=0.
\end{equation}
\CCC In order to see this we firstly choose $C^h(\nu)$ that satisfies the expression \eqref{kirill6} and then apply Lemma \ref{lemmaA8} and Remark \ref{remnak1} to the sequence $(1/h(\widehat{\vect{u}}^h)^{\nu}-C^h(\nu))$. \BBB
From \eqref{kirill6} it follows that
\begin{equation} \label{kirill7} 
\int_{\Omega} \kappa (\mathfrak{a} (\nu)-x_3 \nabla \mathfrak{b}(\nu))\, dx=0, \quad \int_{\omega} \bar{\kappa} \mathfrak{b} (\nu)\, d\widehat{x}=0. 
\end{equation}  
Furthermore, by applying Lemma \ref{lemmaA10} and Remark \ref{remnak2}  to $\vect{\psi}^h(\nu)$, there exists another subsequences $(\varphi^h(\nu))_{h>0}\subset  H^2_{\#}(\omega)), (\widetilde{\vect{\psi}}^{h}(\nu))_{h>0}\subset H^1_{\#}(\Omega;\RR^3)), (o^h(\nu))_{h>0}\subset$ $  L^2(\Omega;\RR^{3\times 3})$ such that
\begin{equation} \label{refk8} 
e_h\left(\vect{\psi}^h(\nu)\right)=\iota\left(-x_3\nabla^2_{\widehat{x}}\varphi^h(\nu)\right)+e_h\left(\widetilde{\vect{\psi}}^{h}(\nu)\right)+o^h(\nu),
\end{equation}
where
\begin{equation}
\varphi^h(\nu)\overset{H^1}{\rightarrow}0, \quad  \left\|\nabla^2_{\widehat{x}}\varphi^h(\nu)\right\|_{L^2}\leq C,
\quad \widetilde{\vect{\psi}}^{h}(\nu)\overset{L^2}{\rightarrow}0, \quad \left\|\nabla_h \widetilde{\vect{\psi}}^{h}(\nu)\right\|_{L^2}\leq C,\quad o^h(\nu)\overset{L^2}{\rightarrow}0, \label{eq:80}
\end{equation}
where $C>0$ is independent of $h$.  
By using Lemma \ref{lemmaA17} and Lemma \ref{newlemmamarch212}, we have that there exist $\varphi(\nu) \in L^2(\omega; \dot{H}^2(\CCC \hat{\mathcal{Y}} \BBB))$ (recall Section \ref{notation}), $\vect{g}(\nu)\in L^2(\Omega,\RR^3)$  and $\vect{w}^1(\nu)\in L^2(\Omega,\dot{H}^1(\mathcal{Y},\RR^3))$ such that
\begin{equation}\label{refk9} 
\nabla^2_{\widehat{x}}\varphi^h(\nu)(\widehat{x},t)\xrightharpoonup{2-r} \nabla^2_{\CCC \hat{y} \BBB} \varphi(\nu) (\widehat{x},\CCC \hat{y} \BBB,t),\quad 
e_h \left(\widetilde{\vect{\psi}}^{h}(\nu)\right)\xrightharpoonup{2-r} \mathfrak{C}_{\infty}\left(\vect{w}^1(\nu), \vect{g}(\nu)\right).
\end{equation}
By introducing the function
\begin{equation} \label{refk10} 
\vect{w}(\nu)(x,y,t):=(
      -x_3\partial_{y_1}\varphi (\nu) (\widehat{x},\CCC \hat{y} \BBB,t) ,
       -x_3\partial_{y_2}\varphi (\nu) (\widehat{x},\CCC \hat{y} \BBB,t) ,
       \varphi (\nu) (\widehat{x},\CCC \hat{y} \BBB,t)
   )^T+\vect{w}^1(\nu)(x,y,t),
\end{equation}
we have
\begin{equation}\label{refk11} 
\mathfrak{C}_{\infty}(\vect{w}(\nu)(x,y,t),\vect{g}(\nu)(x,t))=-x_3 \iota\left(\nabla^2_{\CCC \hat{y} \BBB} \varphi (\nu) (\widehat{x}, \CCC \hat{y} \BBB,t)\right)+\mathfrak{C}_{\infty}\left(\vect{w}^1(\nu)(x,y),\vect{g}(\nu)(x)\right).
\end{equation}
From \eqref{kirill1} we conclude \CCC by using \eqref{convergenceextensionmarch21}, \eqref{refk8}, \eqref{refk9}, \eqref{refk10}, \eqref{refk11}: 
$$ \int_{[0,T]} \mathbb{L}\nu=\iota\left(e_{\hat{x}} (\mathfrak{a}(\nu))-x_3 \nabla_{\hat{x}}^2 \mathfrak{b}(\nu)\right)+\mathfrak{C}_{\infty}(\vect{w}(\nu)(x,y,t),\vect{g}(\nu)(x,t)).   $$
From this we have by integration: \BBB
\begin{align}
    \label{kirill2}  
e_{\widehat{x}}(\mathfrak{a}(\nu))_{\alpha\beta}&=\int_{[0,T]\times I \times \cal Y} \mathbb{L}_{\alpha\beta} \nu  \, dy \,dx_3\,dt, \\ \label{kirill4}  (\nabla^2_{\widehat{x}}\mathfrak{b}(\nu))_{\alpha\beta}&=-12\int_{[0,T]\times I \times \cal Y} x_3\mathbb{L}_{\alpha\beta} \nu  \, dy \,dx_3\,dt, \, \alpha, \beta \in \{1,2\}, \\ \label{kirill5} 
\CCC \vect{g}_{\alpha3}(\nu) &=\CCC 2\int_{[0,T] \times \cal Y} \mathbb{L}_{\alpha3} \nu  \, dy \,dt, \ \alpha \in\{1,2\}, \BBB \\ 
 \vect{g}_{i3}(\nu) &= \int_{[0,T] \times \cal Y} \mathbb{L}_{\alpha3} \nu  \, dy\,dt, \\ \label{kirill11} 
 e_y\left(\vect{w}(\nu)\right) &= \int_{[0,T]}\mathbb{L}\nu\, dt-\iota \left(-x_3\nabla^2_{\widehat{x}}\mathfrak{b}(\nu)+e_{\widehat{x}}(\mathfrak{a}(\nu))\right)- \text{sym}\left(0|0|\vect{g}(\nu)\right).  
\end{align}
\CCC Note that when the subsequence in $h$ is chosen such that  \eqref{kirill1} is satisfied, we directly have that $\mathfrak{a}(\nu)$, $\mathfrak{b}(\nu)$, $\vect{w}(\nu)$, $\vect{g}(\nu)$ do not depend on further subsequence, since the solutions of \eqref{kirill2}-\eqref{kirill11} that satisfy \eqref{kirill7}  are unique for given right hand side (to see it for \eqref{kirill2} and \eqref{kirill11} we use Remark \ref{korrem}).    \BBB

Since $\nu \in L^1(0,T)$ is arbitrary, the existence of $\mathfrak{a}$, $\mathfrak{b}$, $\vect{w}$, $\vect{g}$ in appropriate spaces (see \eqref{nakk101}-\eqref{kirill10}) that satisfy \eqref{cr2gamma0} follows from \eqref{kirill7} and  
\eqref{kirill2}-\eqref{kirill11}.
\CCC Indeed,  to  obtain $\mathfrak{a}$ we firstly 
conclude from Lemma \ref{kornper} that the space 
$$ \mathcal{S}:= \{e_{\hat{x}} (\mathfrak{s})\ : \ \mathfrak{s} \in H_{\#}^1(\omega;\R^2)  \} $$
is closed in $L^2(\omega;\R^{2 \times 2})$ with respect to strong (and thus by convexity  with respect to weak) topology. 
From this we have by Lebesgue theorem  (see e.g. \cite[Theorem 3.20]{Folland})
and by taking $\nu=\chi_{[t,t+\delta]}$ (recall Section \ref{notation}) in \eqref{kirill2}
 that for a.e. $t \in [0,T]$ there exists  $\mathfrak{a} \in L^{\infty}(0,T;H^1_{\#}(\omega;\R^2))$  such that 
$$ (e_{\hat{x}}(\mathfrak{a(t)}))_{\alpha\beta}=\lim_{\delta \to   0}\frac{1}{\delta} \int_{[t,t+\delta]\times I \times \mathcal{Y}} \mathbb{L}_{\alpha\beta}(\tau)\,d\tau \,dy \,dx_3=\int_{I \times \mathcal{Y}} \mathbb{L}_{\alpha\beta}(t)\, dy\,dx_3,  $$
where the limit is taken with respect to weak topology in $L^2(\omega;\R^{2 \times 2})$. According to Remark \ref{korrem} $\mathfrak{a}(t)$ is unique up to a constant.  
In the similar way we obtain the existence of $\mathfrak{b}$, $\vect{w}$, $\vect{g}$. Uniqueness of $(\mathfrak{a},\mathfrak{b})$ follows from the fact that $(\mathfrak{a},\mathfrak{b})$ have to belong to the space $L_{\kappa,0}$ (see \ref{kirill6}), uniqueness for $\vect{w}$ follows from Remark \ref{korrem}, while for $\vect{g}$ is direct.   

 \BBB

\BBB

In order to obtain \eqref{cr1gamma000} and (\ref{cr3gamma0}), we have from Proposition \ref{prop 2may10} 
\begin{equation} \label{nakk2000} 
\displaystyle\frac{\ee}{h^2}\left\|e_h\left(\vect{u}^h\right)\right\|_{H^1(0,T;L^2(\Omega^{h}_{f};\RR^{3  \times 3})) }\leq C,
\end{equation}
therefore from \eqref{eq:78march21}, using $\eps \ll h$ 
\begin{equation}\label{eq:81march21}
\displaystyle\frac{\ee}{h^2}\left\|e_h\left(\vect{u}^h_f\right)\right\|_{L^{\infty}(0,T;L^2(\Omega;\RR^{3  \times 3})) }\leq C.
\end{equation}
From Lemma \ref{lemma 1}, we have
\begin{equation}\label{eq:82march21}
\displaystyle\frac{1}{h^2}\|\vect{u}^h_f\|_{L^{\infty}(0,T;L^2(\Omega;\RR^3))}\leq C\frac{\ee}{h^2} \|\nabla_h \vect{u}^h_f\|_{L^{\infty}(0,T;L^2(\Omega;\RR^{3  \times 3})) }\leq C\frac{\ee}{h^2} \|e_h(\vect{u}^h_f)\|_{L^{\infty}(0,T;L^2(\Omega;\RR^{3  \times 3})) }.
\end{equation}
 By applying  Corollary \ref{cornewlemmamarch21} to $\displaystyle\frac{1}{h^2}\vect{u}^h_f$, we obtain that there exists a function $\vect{u}^{0}_f\in L^{\infty}(0,T;L^2(\Omega; H^1(\mathcal{Y};\RR^3)))$ and a subsequence (not relabeled) such that \eqref{cr1gamma000} and (\ref{cr3gamma0}) hold.

The fact that $\vect{u}_f^0$ is supported in $\Omega \times \mathcal{Y}_{f}^{x_3}$ follows from the fact that $\vect{u}_f^h$ is zero outside $\Omega_f^h$. \eqref{cr5gamma00} \CCC follows \BBB from \eqref{cr3gamma0}, \eqref{eq:78march21} and \eqref{nakk2000}. \CCC To see this, note that \eqref{nakk2000} gives the necessary compactness and that the derivative in time can be moved to the test function by integration by parts in time variable and then we can use \eqref{cr3gamma0} and \eqref{eq:78march21}.   
This also gives that  $e_y(\vect{u}_f^0)\in H^1(0,T;L^2(\Omega \times Y;\RR^{3\times 3}))$. To conclude that  $\vect{u}^0_f \in H^1(0,T; L^2(\Omega; H^1(\mathcal{Y}_f;\RR^3))) $ we use  an approximation with convolution (in time) and Korn's inequality. Consequently, we have \eqref{cr5gamma00}. 
\BBB
The property $\diver_y \vect{u}_f^0={\rm tr } \,e_y(\vect{u}_f^0)=0$ follows from convergence \eqref{cr3gamma0}, \eqref{eq:78march21} and the property $\frac{\varepsilon}{h^2} \diver \vect{u}^h=0$ on $\Omega^h_f$. The fact that  $\vect{u}_f^0(0)=0$ follows from \eqref{eq:51rescaledver}, \eqref{cr3gamma0}, \eqref{cr5gamma00} and \eqref{eq:78march21}
\CCC after passing to the limit in 
$$\frac{\eps}{h^2} e_h(\vect{u}^h(t)) \chi_{\Omega_f^h}=\frac{\eps}{h^2}\int_0^t e_h(\partial_t \vect{u}^h)(\tau)) \chi_{\Omega_f^h}\, d\tau,    $$
and using Korn's inequality. 
\BBB
 This establishes \eqref{kreso1}. 
Since for every $\nu \in L^1(0,T)$ we have  from  \eqref{kirill6} and \eqref{eq:82march21} 
$$ 0=\pi_{1/h}\int_{\Omega} \kappa^h(\vect{u}^h)^{\nu}= \pi_{1/h}\int_{\Omega} \kappa^h(\widehat{\vect{u}}^h+\vect{u}_f^h)^{\nu}\to  \lim_{h \to 0} \left(h\pi_{1/h} C^h(\nu)\int_{\Omega}\kappa^h\right),$$
we conclude that  $ \lim_{h \to 0} \left(h\pi_{1/h} C^h(\nu)\right)=0$.  From \eqref{convergenceextensionmarch21} we conclude \eqref{crgamma0}.

 
 \CCC Next we show (\ref{cr4gamma0}). 
 The approach of \cite{clopeaumikelic2001} needs to be adapted to this more complex framework. \BBB
  We take $\varphi \in C^1_{\#}(\bar{\Omega})$, $\nu \in L^2(0,T)$ and compute \CCC
{\allowdisplaybreaks
\begin{align}
        0&=-\displaystyle\int_0^T\int_{\Omega^{h}_f}\left(\frac{\partial_1 u^{h}_1}{h}+\frac{\partial_2 u^{h}_2}{h}+\frac{\partial_3 u^{h}_3}{h^2}\right)\varphi(x)\nu(t)\,dx\,dt=-\displaystyle\int_0^T\int_{\Omega^{h}_f}\text{div}_h\frac{\vect{u}^h}{h}\varphi(x)\nu(t)\,dx\,dt\\*
    &=-\displaystyle\int_0^T\int_{\Omega^{h}_f}\left(\frac{\partial_1 \widehat{u}^{h}_1}{h}+\frac{\partial_2 \widehat{u}^{h}_2}{h}+\frac{\partial_3 \widehat{u}^{h}_3}{h^2}\right)\varphi(x)\nu(t)\,dx\,dt-\displaystyle\int_0^T\int_{\Omega^{h}_f}\left(\frac{\partial_1 u^{h}_{f,1}}{h}+\frac{\partial_2 u^{h}_{f,2}}{h}+\frac{\partial_3 u^{h}_{f,3}}{h^2}\right)\varphi(x)\nu(t)\,dx\,dt\\
    &=-\displaystyle\int_0^T\int_{\Omega^{h}_f}\left(\frac{\partial_1 \widehat{u}^{h}_1}{h}+\frac{\partial_2 \widehat{u}^{h}_2}{h}+\frac{\partial_3 \widehat{u}^{h}_3}{h^2}\right)\varphi(x)\nu(t)\,dx\,dt\\ & \quad +\displaystyle\int_0^T\int_{\CCC \Omega\BBB}\left(\frac{ u^{h}_{f,1}}{h}\partial_1\varphi(x)+\frac{ u^{h}_{f,2}}{h}\partial_2\varphi(x)+\frac{u^{h}_{f,3}}{h^2}\partial_3 \varphi(x)\right)\nu(t)\,dx\,dt\\
    &\rightarrow-\displaystyle\int_0^T\int_{\Omega}|\mathcal{Y}_f(x_3)|\text{div}_{\widehat{x}}\mathfrak{a}(\widehat{x},t)\,\varphi (x)\nu(t)\,dx\,dt+\displaystyle\int_0^T\int_{\Omega}|\mathcal{Y}_f(x_3)|x_3\text{div}_{\widehat{x}}\nabla_{\widehat{x}}\mathfrak{b}(\widehat{x},t)\,\varphi (x)\nu(t)\,dx\,dt\\*
    & \quad-\displaystyle\int_0^T\int_{\Omega}\int_{\mathcal{Y}_f(x_3)}\text{div}_{y}\vect{w}(x,y,t)\varphi (x)\nu(t)\,dy\,dx\,dt\\*
    & \label{conv2may11} \quad -\displaystyle\int_0^T\int_{\Omega}|\mathcal{Y}_f(x_3)|g_3(x,t)\varphi (x)\nu(t)\,dx\,dt+\displaystyle\int_0^T\int_{\Omega}\int_{\mathcal{Y}_{f}(x_3)} u^{0}_{f,3} (x,y,t)\partial_3\varphi (x)\nu(t)\,dy\,dx\,dt.
\end{align}
}
\BBB
Note that in \CCC concluding \eqref{conv2may11} \BBB   we can always take test functions $\varphi$ that belong to $C^1_{\#}(\Omega_p) \cap H^1(\Omega_p)$. 
(\ref{cr4gamma0}) then follows directly from (\ref{conv2may11}) by the density argument and the arbitrariness of $\nu$.

The convergence \eqref{cr5gamma0} (on a subsequence) to some $p\in (H_T^1L^2(\Omega \times \mathcal{Y}))'$ is a direct consequence of Proposition \ref{pressureregularity}.\footnote{To obtain the compactness statement in the context of two-scale convergence one can use Riesz representation theorem for the elements of $H_T^1L^2$.} 	The fact that $p=0$ on $\Omega \times \mathcal{Y}_{s}^{x_3}$ is then the direct consequence of the fact that $\widehat{p}^h=0$ on $\Omega^h_s$. 
To prove its independence of $y$ i.e.
\begin{equation}\label{pressureindependencemay11}
p=p(x,t)\,\,\text{in}\,\, \Omega\times\mathcal{Y}_f^{x_3}\times(0,T),
\end{equation}
we can choose as test function $\vect{v}^{h}(x,t)=\ee h \sum_{i=1}^m \vect{\tau}^i\left(\displaystyle\frac{\widehat{x}}{\ee},\frac{x_3}{\frac{\ee}{h}}\right)\xi^i(x)\nu(t)$,  such that $\nu \in C^1([0,T])$, $\nu(T)=0$, $\xi^i \in C_c^1 (\Omega)$, $\textrm{supp } \xi^i\subset U_i$,  $\vect{\tau}^i \in \CCC C_c^\infty \BBB(Y_f^i;\RR^3)$,
for $i=1,\dots,m$
and plug it in \eqref{varformporous2gamma0may10}. 
We obtain as a consequence of  Proposition \ref{prop 2may10} \CCC
\begin{equation}
0\leftarrow-_{(H_T^1L^2_f)'}\big\langle \frac{p^h}{h}\cdot\frac{\text{div}_{h}\vect{v}^{h}}{h}\big\rangle_{H_T^1L^2_f}=- _{(H_T^1L^2)'}\big\langle \frac{p^h}{h}\chi_{\Omega^{h}_f}\cdot\frac{\text{div}_{h}\vect{v}^{h}}{h}\big\rangle_{H_T^1L^2}.
\end{equation}
Consequently
\begin{align*}
& 0=-\lim_{h \to 0} {} _{(H_T^1L^2)'}\big\langle \frac{p^h}{h}\chi_{\Omega^{h}_f}\cdot\frac{\text{div}_{h}\vect{v}^{h}}{h}\big\rangle_{H_T^1L^2}=\sum_{i=1}^m \int_{\Omega} \int_{Y}  \chi_{Y_f(x_3)}(y)\cdot p^{\nu}(x,y)\cdot\text{div}_{y}\vect{\tau}^i(y)\xi^i(x) \, dx \,dy\\
&=\displaystyle\sum_{i=1}^m\int_{\Omega}\int_{{Y}_f(x_3)}p^\nu(x, y)\cdot\text{div}_{y}\vect{\tau}^i(y)\xi^i(x) \,dy\,dx=-\displaystyle\sum_{i=1}^m \int_{\Omega} {}_{\mathcal{D}'({Y}_f(x_3);\R^3)}\langle \nabla_{y}p^{\nu}(x,y),\vect{\tau}^i(y)\xi^i(x)\rangle_{\mathcal{D}({Y}_f(x_3);\R^3)}\,dx.
\end{align*}
Here $\nabla_yp^{\nu}$ is the distributional gradient of $p^{\nu}$ \CCC and $\mathcal{D}'(Y_f(x);\R^3)$ is the space of Schwartz distributions on ${Y}_f(x_3)$, while $\mathcal{D}({Y}_f(x_3);\R^3) \equiv C_c^\infty({Y}_f(x_3);\RR^3)$ with appropriate structure. Since the distributional gradient (for a.e. $x\in \Omega$) of  of the restriction of $p^\nu(x,\cdot)$ on ${Y}_f(x_3)$ is zero, this proves the independence of pressure on $y$, see for details e.g., \cite{Schwartz} (notice that, however, its support is $y$-dependent). \BBB
\end{proof}  
\begin{remarkica} \label{remboundary}
In \eqref{conv2may11} we also used the fact that $u_{f,3}^h$ is zero on $\omega \times \{-1/2,1/2\}$ which is a consequence of our geometrical assumption \eqref{nodownup}. We will discuss in  Section \ref{secgen}  the cases when \eqref{nodownup} is not satisfied. 
\end{remarkica}

\subsection{Effective tensors} \label{effective}

\BBB
Before stating the result on limit equations we need to define the effective tensors. 
\begin{enumerate} 
\item	
We define the effective fourth order tensor: For  $\vect{A}, \vect{B},\vect{C},\vect{D}\in\RR^{2\times 2}_{\text{sym}}$, we  define (recall \eqref{defAAA})
\begin{align}\nonumber 
	&\mathbb{A}^{\rm hom}(\vect{A}, \vect{B}):(\vect{C}, \vect{D}):=\\ &\displaystyle\quad \label{gassmantensorapril01}\int_{I}\int_{\mathcal{Y}_s(x_3)}\mathbb{A}(x_3,y)\left[\iota(\vect{A}-x_3\vect{B})+\mathfrak{C}_{\infty}(\vect{w}_{1,\vect{A},\vect{B}}, \vect{g}_{1,\vect{A},\vect{B}})\right]:\left[\iota(\vect{C}-x_3\vect{D})\right]dydx_3,
\end{align}
where  $(\vect{w}_{1,\vect{A},\vect{B}},\vect{g}_{1,\vect{A},\vect{B}}) \in L^2(I;\dot{H}^1(\mathcal{Y}_s(x_3);\mathbb{R}^3)) \times L^2(I;\mathbb{R}^3)$ is the  unique solution of
\begin{align}
	&\label{limitequationhom1}\displaystyle\int_{I}\int_{\mathcal{Y}_s(x_3)}\mathbb{A}(x_3,y) \left[\iota(\vect{A}-x_3\vect{B})+\mathfrak{C}_{\infty}(\vect{w}_{1,\vect{A},\vect{B}}, \vect{g}_{1,\vect{A},\vect{B}})\right]: \mathfrak{C}_{\infty}(\vect{\zeta}, \vect{r})\,dy\,dx_3=0;\\
& \nonumber
\forall \vect{\zeta} \in L^2(I;H^1_{\#}(\mathcal{Y}_s(x_3);\mathbb{R}^3)), \ \vect{r} \in L^2(I;\mathbb{R}^3). 
\end{align} 
\item For \CCC $x_3 \in J_p$,\BBB we  define the second order tensor
\begin{equation}\label{btensorapril1}
	\mathbb{B}^{H}(x_3):=\displaystyle\int_{\mathcal{Y}_s(x_3)}\mathbb{A}(x_3,y)\,\mathfrak{C}_{\infty}\left(\widetilde{\vect{w}}_{2,x_3},\widetilde{\vect{g}}_{2,x_3}\right)\,dy,
\end{equation}
where $(\widetilde{\vect{w}}_{2,x_3}, \widetilde{\vect{g}}_{2,x_3})\in \dot{H}^1(\mathcal{Y}_s(x_3);\mathbb{R}^3) \times \mathbb{R}^3$ is the  unique solution of 
\begin{equation}
\label{eq unique uationapril16}
    \begin{split}
        \displaystyle\int_{\mathcal{Y}_s(x_3)}\mathbb{A}(x_3,y)\, \mathfrak{C}_{\infty}(\widetilde{\vect{w}}_{2,x_3}, \widetilde{\vect{g}}_{2,x_3}): \mathfrak{C}_{\infty}(\vect{\zeta}, \vect{r})\,dy&=\displaystyle\int_{\mathcal{Y}_f(x_3)}\, \left(\text{div}_{y}\, \CCC \widehat{\vect{\zeta}} \BBB+r_3\right)\,dy,\\  & \forall \vect{\zeta} \in {H}^1(\mathcal{Y}_s(x_3);\mathbb{R}^3),\ \vect{r} \in \mathbb{R}^3,
    \end{split}
\end{equation}
\CCC where $\widehat{\vect{\zeta}}$ is any extension of $\zeta$ to $H^1(\mathcal{Y};\R^3)$. \footnote{\CCC The right hand side of \eqref{eq unique uationapril16}is actually independent of this extension, which can be seen by using Gauss formula that converts the integral over $\mathcal{Y}_f(x_3)$ in the integral over $\partial \mathcal{Y}_f(x_3)=\partial \mathcal{Y}_s(x_3)$. However we write it in this way because of the way how Theorem \ref{effectiveequationsmay11} is proved.  \BBB} \BBB
\item For $x_3 \in J_p$, we  define the second order tensor 
\begin{equation}\label{ctensorapril1}
	\mathbb{C}^{H}_{ij}(x_3):=-\displaystyle\int_{\mathcal{Y}_f(x_3)}\text{tr}\,\mathfrak{C}_{\infty}(\CCC \widehat{\vect{w}}^{ij}_{x_3}\BBB, \vect{g}^{ij}_{x_3} \BBB)\,dy,
\end{equation}
where for $i,j=1,2,3$
$(\vect{w}^{ij}_{x_3}, \vect{g}^{ij}_{x_3}) \in \dot{H}^1(\mathcal{Y}_s(x_3);\mathbb{R}^3) \times \RR^3$ is the  unique solution of the following cell problem:
\begin{equation}
\displaystyle\int_{\mathcal{Y}_s(x_3)}\mathbb{A}(x_3,y) \left[\vect{e}^{i}\odot \vect{e}^{j}+\mathfrak{C}_{\infty}(\vect{w}_{x_3}^{ij}, \vect{g}_{x_3}^{ij})\right]: \mathfrak{C}_{\infty}(\vect{\zeta}, \vect{r})\,dy=0, \quad \forall \vect{\zeta} \in {H}^1(\mathcal{Y}_s(x_3);\mathbb{R}^3),\ \vect{r} \in \mathbb{R}^3.\label{cellproblem1april1}  \end{equation}
\CCC and $\widehat{w}^{ij}_{x_3}$ is any extension of 
${w}^{ij}_{x_3}$ to $H^1(\mathcal{Y};\R^3)$ \footnote{\CCC Again, it can be shown, by using Gauss formula, that the expression \eqref{ctensorapril1} is independent of this extension.   \BBB}. 
\BBB
\item For $x_3 \in J_p$, we  define the second order tensor
\begin{equation}\label{ktensorapril}
	\mathbb{K}_{ij}(x_3):=\displaystyle\int_{\mathcal{Y}_f(x_3)}q^{i}_{j,x_3}\,dy.\quad \text{(the permeability tensor)}, \quad i,j=1,2,3,
\end{equation}
where for $i=1,2,3$, $\mathbf{q}^i_{x_3} \in H_0^1(\mathcal{Y}_f(x_3);\mathbb{R}^3)$, $\pi^i_{x_3} \in L^2(\mathcal{Y}_f(x_3))$ ($\equiv L^2(Y_f(x_3))$) are the weak solutions on torus of\footnote{\CCC Recall that a function belonging to  $H_0^1(\mathcal{Y}_f(x_3);\mathbb{R}^3)$, when extended by zero, belongs to $H^1(\mathcal{Y};\R^3)$.  \BBB} .  
\begin{equation}\label{diferformfluid2march301}
	\begin{array}{lcl}
		\begin{cases}
		-\Delta_{y}\vect{q}^{i}_{x_3}+\nabla_{y}\pi^i_{x_3}=\vect{e}^i, \quad \text{in}\, \mathcal{Y}_f(x_3) &\\
		\text{div}_{y}\vect{q}^i_{x_3}=0, \quad \text{in }\, \mathcal{Y}_f(x_3).
	\end{cases}
		& \textrm{i.e.} & 
		\begin{cases}
			\int_{Y_f(x_3)} \nabla_{y}\vect{q}^i_{x_3}\nabla_y \vect{\psi}\, dy-  \int_{Y_f(x_3)}\pi_{x_3}^i\textrm{div}_y \vect{\psi}\, dy=\int_{Y_f(x_3)} \psi_3\,dy,& \\ \quad \forall \vect{\psi} \in H_0^1(\mathcal{Y}_f(x_3);\R^3); &\\
			\diver_{y}\vect{q}_{x_3}^i=0, \quad \text{in }\, \mathcal{Y}_f(x_3).
		\end{cases}
	\end{array} 
\end{equation}

\item For $x_3 \in J_p$, we define scalar
\begin{equation} \label{defM00}  M_0(x_3):=\displaystyle\int_{\mathcal{Y}_f(x_3)}\text{tr}\,\mathfrak{C}_{\infty}(\CCC \widehat{\widetilde{\vect{w}}}_{2,x_3}\BBB, \widetilde{\vect{g}}_{2,x_3})\,dy, \end{equation} 
where again $(\widetilde{\vect{w}}_{2,x_3}, \widetilde{\vect{g}}_{2,x_3})\in \dot{H}^1(\mathcal{Y}_s(x_3);\mathbb{R}^3) \times \mathbb{R}^3$ are unique solutions of 
\eqref{eq unique uationapril16} \CCC and $\widehat{\widetilde{\vect{w}}}_{2,x_3}$ is any extension of $\widetilde{\vect{w}}_{2,x_3}$ to $H^1(\mathcal{Y};\R^3)$. \footnote{\CCC Again, by using Gauss formula, the expression \eqref{defM00} is independent of the extension. \BBB}   \BBB
\end{enumerate} 
The following proposition gives us the important properties of the tensors defined above
\begin{propozicija} 
\label{tensorprop} 
The solutions of \eqref{limitequationhom1}, \eqref{eq unique uationapril16}, \eqref{cellproblem1april1}, \eqref{diferformfluid2march301}  exist and are unique.     The tensors  $\mathbb{B}^H$, $\mathbb{C}^H$, $\mathbb{K}$, $M_0$ are piecewise constant and following properties are satisfied: 
\begin{enumerate} 
	\item The tensor $\mathbb{A}^{\rm hom}$ is symmetric and there exists  $\nu_{\mathbb{A}^{\rm hom}}$ such that 
	$$ \nu_{\mathbb{A}^{\rm hom}}^{-1} (|\vect{A}|^2+|\vect{B}|^2) \leq \mathbb{A}^{\rm hom}(\vect{A}, \vect{B}):(\vect{A}, \vect{B}) \leq \nu_{\mathbb{A}^{\rm hom}} (|\vect{A}|^2+|\vect{B}|^2),\quad   \forall \vect{A},\vect{B} \in \mathbb{R}_{\rm sym} ^{2 \times 2}. $$
\item For all $x_3 \in J_p$, the tensors $\mathbb{B}^H(x_3)$, $\mathbb{C}^H(x_3)$ are symmetric and $\mathbb{B}^H(x_3)=\mathbb{C}^H(x_3)$. 
\item For all $x_3 \in J_p$, the tensor $\mathbb{K}(x_3)$ is symmetric and positive semidefinite and there exists $\nu_{\mathbb{K}}>0$ such that
 \begin{equation} \label{tonci1}
 \begin{split}  
  &\nu_{\mathbb{K}} \leq \mathbb{K}_{33}(x_3) \leq  \nu_{\mathbb{K}}^{-1},\quad i=1,2,3, \, \forall x_3 \in J_K; \\ &
  \mathbb{K}_{33} (x_3)=0, \quad \forall x_3 \in J_p \backslash J_K. 
  \end{split} 
  \end{equation}
 \item There exists $\nu_{M_0}>0$ such that $ M_0(x_3) \geq  \nu_{M_0}>0, \, \forall x \in \Omega_p.$ 
\end{enumerate} 	
\end{propozicija} 	
	\begin{proof} 
	The existence and uniqueness of the solutions of \eqref{limitequationhom1}, \eqref{eq unique uationapril16}, \eqref{cellproblem1april1}, \eqref{diferformfluid2march301} go by standard arguments, using Lax-Milgram and \eqref{tensorA}.  
	To prove the claim 1., notice that  
	\begin{equation}\label{eq:abril1}
		\left(\vect{w}_{1,\vect{A},\vect{B}}(x_3,\cdot), \vect{g}_{1,\vect{A},\vect{B}}(x_3,\cdot)\right)=\displaystyle\sum_{\alpha,\beta=1,2}[\vect{A}_{\alpha\beta}\left(\vect{w}^{\alpha\beta}_{x_3},\vect{g}^{\alpha\beta}_{x_3}\right)-x_3\vect{B}_{\alpha\beta}\left(\vect{w}^{\alpha\beta}_{x_3},\vect{g}^{\alpha\beta}_{x_3}\right)],
	\end{equation} 
where $\left(\vect{w}^{\alpha\beta}_{x_3},\vect{g}^{\alpha\beta}_{x_3}\right)$ are defined in \eqref{cellproblem1april1}.
\CCC This can be seen by taking $\vect{\zeta}=\xi\widetilde{\vect{\zeta}}$, $\vect{r}=\xi \widetilde{\vect{r}}$, where  $\xi \in L^2(I)$, $\widetilde{\vect{\zeta}} \in H^1(\mathcal{Y}_s(x_3);\R^3) $, $\widetilde{\vect{r}} \in \R$ in \eqref{limitequationhom1}, and use linearity of \eqref{cellproblem1april1}.   \BBB 
By using the extension on $\mathcal{Y}_f(x_3)$ (see \cite[Chapter 4]{Oleinik}), for each $\vect{w}^{\alpha\beta}_{x_3}$, for $\alpha, \beta=1,2$, the claim 1. follows the proof of \cite[Proposition 3.4]{Marinvelciczubrinic2022}. 
To prove the claim 2., notice that the symmetricity of $\mathbb{C}^H$ is obvious, since $\vect{w}^{ij}_{x_3}=\vect{w}^{ji}_{x_3}$, $\vect{g}^{ij}_{x_3}=\vect{g}^{ji}_{x_3}$, for $i,j=1,2,3$. We prove $\mathbb{B}^h=\mathbb{C}^H$. 

Using $\left(\widetilde{\vect{w}}_{2,x_3},\widetilde{\vect{g}}_{2,x_3}\right)$ as test functions in (\ref{cellproblem1april1}) and by virtue of (\ref{eq unique uationapril16}), we obtain that
\begin{align}
\mathbb{B}^{H}_{ij}(x_3)&=\displaystyle\int_{\mathcal{Y}_s(x_3)}\left(\mathbb{A}(x_3,y)\mathfrak{C}_{\infty}\left(\widetilde{\vect{w}}_{2,x_3},\widetilde{\vect{g}}_{2,x_3}\right)\right)_{ij}\,dy=\displaystyle\int_{\mathcal{Y}_s(x_3)}\mathbb{A}(x_3,y)\mathfrak{C}_{\infty}\left(\widetilde{\vect{w}}_{2,x_3},\widetilde{\vect{g}}_{2,x_3}\right):\left(\vect{e}^{i}\odot \vect{e}^{j}\right)\,dy\\
&=-\displaystyle\int_{\mathcal{Y}_s(x_3)}\mathbb{A}(x_3,y)\mathfrak{C}_{\infty}\left(\widetilde{\vect{w}}_{2,x_3},\widetilde{\vect{g}}_{2,x_3}\right):\mathfrak{C}_{\infty}\left(\vect{w}^{ij}_{x_3}, \vect{g}^{ij}_{x_3}\right)\,dy=-\displaystyle\int_{\mathcal{Y}_f(x_3)}\text{tr}\,\mathfrak{C}_{\infty}(\widehat{\vect{w}}^{ij}_{x_3}, \vect{g}^{ij}_{x_3})\,dy=\mathbb{C}^{H}_{ij}(x_3).
\end{align}
Next, we prove the claim 3.:  
By testing \eqref{diferformfluid2march301} with $\mathbf{q}^j_{x_3}$, for $j=1,2,3$, we obtain 
$$\mathbb{K}_{ij}(x_3)= \langle \nabla \mathbf{q}^i_{x_3}, \nabla \mathbf{q}^j_x \rangle_{L^2(\mathcal{Y}_f(x_3);\mathbb{R}^3)}, \quad i,j=1,2,3.$$
\eqref{tonci1} follows from the fact that $\mathbf{q}^3_{x_3}$ is constant 
only if it is zero. On the other hand  $\mathbf{q}^3_{x_3}=0$ is a solution of \eqref{diferformfluid2march301} for $i=3$ if and only if $\mathcal{Y}_f(x_3)$ doesn't intersect $\{y_3=0\}$ (in this case $\pi_{x_3}^i(y)=y_3$). 
This can be seen directly from the weak form of \eqref{diferformfluid2march301}, by doing integration by part in the second term on the left hand side. \BBB
 
Uniformity with respect to $x_3$ is the consequence of the geometric assumptions, where only finite number of different domains $\mathcal{Y}_f(x_3)$ appear. 
It remains to prove 4. By virtue of (\ref{eq unique uationapril16}) we have
\begin{equation} \label{defM0}
M_0(x_3)=\displaystyle\int_{\mathcal{Y}_s(x_3)}\mathbb{A}(x_3,y)\, \mathfrak{C}_{\infty}(\widetilde{\vect{w}}_{2,x_3}, \widetilde{\vect{g}}_{2,x_3}): \mathfrak{C}_{\infty}(\widetilde{\vect{w}}_{2,x_3}, \widetilde{\vect{g}}_{2,x_3})\,dy>0.
\end{equation}
Again, the uniformity with respect to $x_3$, is the consequence of the geometric assumptions. 
\end{proof} 

\subsection{Limit equations} 
\label{limit}
Before proving the main theorem which gives the limit problem  we prove Lemma \ref{lmtestfunc} which is fundamental for constructing the appropriate test functions from which we will be able to conclude the condition that the limit pressure satisfies at the interface. 
For proving Lemma \ref{lmtestfunc}  we need a certain construction of Lipschitz sets that is used in Lemma \ref{lmtestfunc} and that is provided in Lemma \ref{konstrakta1}. 

For $p\in \R^2$ and $r>0$ we denote by $G(p,r)\subset \R^3$ the infinite cylinder whose base is $B(p,r)$
\CCC
\begin{lema} \label{konstrakta1} 
	Let $C_i \subset Y$, $i=1,2$, be connected, Lipschitz sets that are associated sets of open, Lipschitz sets  $\mathcal{C}_i \subset \mathcal{Y}$, $i=1,2$, respectively. Assume that $C_1 \cap C_2 \cap\{y_3=0\} \neq \emptyset$. Then there exist sets $\widetilde{C}_i\subset Y$, $i=1,2$, that are associates sets of
	$\widetilde{\mathcal{C}}_i \subset \mathcal{Y}$, $i=1,2$, respectively such that:
	satisfy: 
	\begin{itemize} 
		\item[(i)] For $i=1,2$, $\widetilde{C}_i$ are connected, with smooth boundary and $\overline{\widetilde{C}}_i \backslash \{y_3=1\} \subset C_i$;
		\item[(ii)] For $i =1,2$, $\overline{\widetilde{C}}_i$ doesn't intersect the boundary of $\partial Y$ except at $\{y_3 \in \{0,1\}\}$, i.e. 
		we have
		$$ \overline{\widetilde{C}_i}  \cap \partial Y =\overline{\widetilde{C}_i} \cap \{y_3 \in \{0,1\}\}, \quad i=1,2. $$ 
		\item[(iii)] 
		There exists $p \in (0,1)^2$, $\eps>0$ such that 
		$$\widetilde{C}_1 \cap \{y_3\in[0,\eps] \cup [1-\eps,1)\}=\widetilde{C}_2\cap \{y_3\in[0,\eps] \cup [1-\eps,1)\}=G(p,\eps) \cap \{y_3\in[0,\eps] \cup [1-\eps,1)\}. $$
		
	\end{itemize} 
	\begin{proof} 
		We will construct the sets $\widetilde{C}_i$ as cylindrical sets around certain simple (non-intersecting) smooth curve that is contained in $C_i$ and that doesn't intersect $\partial Y$ or $\partial C_i$ (except at endpoints). To this end we take $p \in  (0,1)^2$ that $(p,0)  \in C_1 \cap C_2$. As a consequence of the fact that $\mathcal{C}_1,\mathcal{C}_2$ are open, there exists $\widetilde{r}>0$ such that $(B((p,0),\widetilde{r}) \cup B((p,1),\widetilde{r}))\cap Y \subset C_1 \cap C_2$. Since $C_i$ is connected, for $i=1,2$, there exists a polygonal line $(p,0)s^i_1\dots s^i_n(p,1)$   contained in $C_i$ that doesn't intersect the boundary of $C_i$ except at points $(p,0),(p,1)$. Without loss of generality we can assume that 
		\begin{equation} \label{uvdist1} 
			\dist(s^i_j,(p,0))> \widetilde{r}, \ \dist(s^i_j,(p,1))> \widetilde{r} \quad  j=1,\dots,n.  
		\end{equation} 
		Since both $C_i$ are open, by slightly perturbing the points $s^i_j$ for $j=1,\dots,n$,  we can assume that  there are no four coplanar points in the set  $\{(p,0),s^i_1,\dots,s^i_{n},(p,1)\}$ and that the condition \eqref{uvdist1} is satisfied. Since there are no four coplanar points the polygonal line $(p,0)s_1^i\dots s_n^i(p,1)$ is simple. We add four more points in the line in the following way: we define $t^i$ as intersection of the line $(p,0)s^i_1$ with $B((p,0),\widetilde{r})$ and $q^i$ as intersection of the line $s^i_{n-1}(p,1)$ with $B((p,1),\widetilde{r})$. We define the polygonal  line $(p,0)(p,\widetilde{r}/2)t^is^i_1\dots s^i_n q^i (p,1-\widetilde{r}/2)(p,1)$. The new polygonal line is contained in $C_i$, doesn't intersect the boundary of $C_i$ except at the endpoints and is simple. Then we can smoothen it at its corners in a way that it keeps these properties, i.e. we can, for $i=1,2$, find a smooth curve $c^i:[a,b]\to \R^3$ ($a,b \in \R$, $a<b$)  that satisfies 
		\begin{itemize}
			\item $c^i(a)=(p,0)$, $c^i(b)=(p,1)$; 
			\item $c^i$ is injective, its image (except at the last point) is contained in $C_i$; 
			\item $\{c^i(t):t \in (a,b)\} \cap \partial Y=\{c^i(t):t \in (a,b)\} \cap \partial C_i=\emptyset$; 
			\item there exists $\widetilde{\eps}>0$ such that 
			$$ \{c^i(t):t \in [a,b]\} \cap \{y_3\in[0,\widetilde{\eps}] \cup [1-\widetilde{\eps},1)\} =\{c^i(t):t \in [a,a+\widetilde{\eps}]\cup[b-\widetilde{\eps},b]\}=(p,0)(p,\widetilde{\eps}) \cup (p,1-\widetilde{\eps})(p,1).  $$
		\end{itemize} 	 
		The domain $\widetilde{C}_i$ we can define as the (small enough) cylindrical neighbourhood of the curve $c^i$, using its Frenet frme, This is diffeomorphic to cylinder, i.e. there exists $\widehat{\eps}>0$ such that 
		$\overline{\widetilde{C}}_i$ is diffeomorphic to
		$\overline{B(0,\widehat{\eps})}\times [a,b]$ 
		(see \cite{tambaca} for a short elementary proof or \cite[Theorem 3.1-1]{ciarlet3} for the proof in the context of shells).   
	\end{proof} 
	\BBB		 
\end{lema} 
\CCC
\begin{lema} \label{lmtestfunc}
	\begin{enumerate} 
		\item If $\mathcal{C} \subset \mathcal{Y}$  is an open Lipschitz set that has empty intersection with $\{y_3=0\}$,
		then  every $\vect{\tau} \in H_0^1(\mathcal{C};\RR^3)$\footnote{\CCC Again, recall that a function $\vect{\tau}$, when extended by zero, belongs to $H^1(\mathcal{Y};\R^3)$. \BBB} that satisfies $\diver_y\vect{\tau}=0$ satisfies also $\int_Y \tau_3 \, dy=0$. 	
		\item Assume that $\mathcal{C}_i \subset \mathcal{Y}$, for $i=1,2$ are non-empty, open Lipschitz sets that satisfy and $|\mathcal{Y} \backslash \mathcal{C}_i|>0$  and 
		whose associated set $C_i\subset Y$ is connected. Moreover, we assume that  	
		that $\mathcal{C}_1 \cap \mathcal{C}_2 \cap\{y_3=0\} \neq \emptyset$. Then for $i=1,2$ there exists $\vect{\tau}^i \in C_c^{\infty}(\mathcal{C}_i;\R^3)$ such that  for the extensions by zero of $\vect{\tau}^i$ on whole $\mathcal{Y}$ (still denoted by $\vect{\tau}^i$), we have $\diver_y \vect{\tau}^1=\diver_y \vect{\tau}^2=0$, $\vect{\tau}^1=\vect{\tau}^2$ on the set $\{y_3=0\}$ and $\int_Y\tau^1_3 \, dy=\int_Y\tau^2_3 \, dy\neq 0$.    
	\end{enumerate} 		
\end{lema} 	
\begin{proof} 
	To prove the first claim  we take $\vect{\tau} \in H_0^1(\mathcal{C}; \RR^3)$ such that $\diver \vect{\tau}=0$ and by doing integration by parts we conclude
	$$0= \int_C \diver_y \vect{\tau} \cdot y_3 \,dy=- \langle \vect{\tau}, \nabla_y y_3\rangle_{L^2(C;\RR^3)}=-\int_Y \tau_3 \,dy,$$
	which concludes the proof. 
	
	To prove the second claim we take for $i=1,2$ the sets $C_i$ which are associated sets of $\mathcal{C}_i$ and make the construction from Lemma \ref{konstrakta1}, obtaining the sets $\widetilde{C}_i$ that satisfy (i)-(iii) of Lemma \ref{konstrakta1}.  
	We take an arbitrary non-zero smooth non-negative function $g:\R^2 \to \R$  with support on $B(p,\eps)$ (recall (iii) of Lemma \ref{konstrakta1}) and extend it to $\R^3$ by taking it independent of $x_3$ variable (this extension we still denote by $g$). 
	By using \cite[Lemma 2.2]{giraultraviart1986} we can for $i=1,2$ find $\widetilde{\vect{\tau}}^i \in H^1(\widetilde{C}^{\eps/2}_i;\R^3)$  such that 
	$$ \diver_y \widetilde{\vect{\tau}}^i=0,\quad \widetilde{\vect{\tau}}^i|_{\partial \widetilde{C}^{\eps/2}_i \cap \{y_3=\eps/2 \textrm{ or } y_3=1-\eps/2\}}=g \vect{e}^3,\, \widetilde{\vect{\tau}}^i|_{\partial \widetilde{C}^{\eps/2}_i \backslash \{y_3=\eps/2 \textrm{ or } y_3= 1-\eps/2\}}=0.      $$
	Here  $$\widetilde{C}_i^{\eps/2}:=\widetilde{C}_i \backslash \left(G(p,\eps) \cap \{y_3\in[0,\eps/2] \cup [1-\eps/2,1)\}\right),  $$
	which is obviously Lipschitz and connected. 
	Next we extend $\widetilde{\vect{\tau}}^i$ on whole $\widetilde{C}_i$, denoting this extension also by $\widetilde{\vect{\tau}}^i$, by taking 
	$$\widetilde{\vect{\tau}}^i:=g \vect{e}^3 \textrm{ on } G(p,\eps) \cap \{y_3\in[0,\eps/2] \cup [1-\eps/2,1)\}. $$  
	Due to the periodic boundary condition we easily see that these extensions satisfy $\widetilde{\vect{\tau}}^i \in H_0^1(\widetilde{\mathcal{C}}_i;\R^3)$, for $i=1,2$,  
	and are also 
	divergence free. To obtain $\vect{\tau}^i \in C_c^{\infty} (\mathcal{C}_i;\R^3)$, we can use convolution with mollifier supported in a small enough neighbourhood of zero. It is important to note that the condition $\widetilde{\vect{\tau}}^1=\widetilde{\vect{\tau}}^2$ on $\{y_3=0\}$ is preserved with convolution (for mollifier supported in small enough neighbourhood of zero) since $\widetilde{\vect{\tau}}^1=\widetilde{\vect{\tau}}^2$ in a neighbourhood of $\{y_3=0\}$. To check the last property note that by integration by parts for any $\vect{\tau} \in H_0^1(\mathcal{C}_i;\R^3)$ we have 
	$$0= \int_{C_i} \diver_y \vect{\tau} \cdot y_3 \,dy=- \langle \vect{\tau}, \nabla_y y_3\rangle_{L^2(C_i;\RR^3)}+\int_{\partial C_i \cap\{y_3=1\}} \tau_3\, d\hat{y} =-\int_{C_i} \tau_3 \,dy+\int_{\partial C_i \cap\{y_3=1\}} \tau_3\, d\hat{y},$$
	from which it follows, using periodicity, that 
	$$ \int_{\mathcal{C}_1} \tau_3^1\, dy= \int_{\mathcal{C}_1 \cap\{y_3=0\}} \tau^1_3\, d\hat{y}=\int_{\mathcal{C}_2 \cap\{y_3=0\}} \tau^2_3\, d\hat{y}= \int_{\mathcal{C}_2} \tau_3^2\, dy. $$
\end{proof} 	
\begin{remarkica} \label{pojednost} 
	If we apply Lemma \ref{lmtestfunc} 2 in the case $\mathcal{C}_1=\mathcal{C}_2=\mathcal{C}$, where $\mathcal{C}\subset \mathcal{Y}$ is an open set with Lipschitz boundary that has non-empty intersection with $\{y_3=0\}$, we have the existence of $\vect{\tau} \in C_c^{\infty} (\mathcal{C};\R^3)$ that satisfies $\int_{\mathcal{C}} \tau_3 \,dy \neq 0$.  
\end{remarkica} 	
 
Next we give the result on the limit problem. 
 \begin{teorem}\label{effectiveequationsmay11}
Let assumptions \eqref{kirill3} and \eqref{forcesassumptions1}  be satisfied. We also suppose that Assumptions \ref{assumption on regions} is satisfied and $	\pi_h\vect{F}^{h}\xrightharpoonup{t,2-r,2} \vect{F} \footnote{Note that as a consequence of \eqref{kirill3} we have $\int_{\omega} \overline{\langle \vect{F}\rangle_{Y}}\,d\widehat{x}=0.$},$
where $\vect{F} \in H^1(0,T; L^2(\Omega \times \mathcal{Y};\mathbb{R}^3))$.  
Then for the limit $p$ obtained in  \eqref{nakk106} of Theorem \ref{gammazeromay11}, we have the additional regularity $p \in L^{\infty}(0,T; L^2(\Omega))\cap L^2(0,T;M)$  and the  limits \eqref{nakk101}-\eqref{nakk106}  satisfy the following equations\footnote{In the second equation, in the last term we can take $\varphi \in  L^2(0,T;M)\cap H^1(0,T;L^2(\Omega_p))$, since $\mathbb{K}_{33}=0$ on $\Omega_p\backslash \Omega_K$ and thus the integral over $\Omega_p$  can be interpreted as the integral over $\Omega_K$}(\CCC recall \eqref{defbar}, \eqref{deflangle}\BBB): 
\begin{equation}\label{diferformpen1*may11}\begin{split}
    &\displaystyle\int_{\omega}\mathbb{A}^{\rm hom}(e_{\widehat{x}}(\mathfrak{a}), \nabla^2_{\widehat{x}}\mathfrak{b}):(e_{\widehat{x}}(\vect{\theta}_*), \nabla^2_{\widehat{x}}\theta_3)d\widehat{x}-\displaystyle\int_{\omega}\int_{J_p}\left(|\mathcal{Y}_f(x_3)| \mathbb{I}-\mathbb{B}^{H}(x_3)\right)p \, dx_3:\left[\iota (e_{\widehat{x}}(\vect{\theta}_{*}))\right]\,d\widehat{x}\\
    &+\displaystyle\int_{\omega} \int_{J_p}\left(|\mathcal{Y}_f(x_3)| \mathbb{I}-\mathbb{B}^{H}(x_3)\right)x_3 p \, dx_3:\left[\iota (\nabla^2_{\widehat{x}}\theta_3)\right]\,d\widehat{x} =\displaystyle\int_{\omega} \overline{\langle\vect{F}\rangle_{Y}}\cdot(\vect{\theta}_*,\theta_3)\, d\widehat{x}-\displaystyle\int_{\omega} \overline{\langle x_3\vect{F}_{*}\rangle_{Y}}\cdot\nabla_{\widehat{x}}\theta_3\, d\widehat{x},\\
    &\quad\forall \left(\vect{\theta}_{*},{\theta}_3\right) \in L_{\kappa,0}, \quad \textrm{for a.e. } t \in (0,T),
\end{split}\end{equation}
\begin{equation}\label{diferformpen0*may11}\begin{split}
    & \displaystyle-\int_0^T\int_{\Omega_p}M_0(x_3) p\, \partial_t \varphi dx\,dt-\displaystyle\int_0^T\int_{\omega}\int_{J_p}(|\mathcal{Y}_{f}(x_3)| \mathbb{I}-\mathbb{B}^{H}(x_3))\partial_t\varphi \, dx_3:\iota\left(e_{\widehat{x}}\left(\mathfrak{a}\right)\right) d\widehat{x}\,dt \\&+\displaystyle\int_0^T\int_{\omega}\int_{J_p}\left(|\mathcal{Y}_f(x_3)|\mathbb{I}-\mathbb{B}^H(x_3)\right)x_3\partial_t\varphi \, dx_3:\iota\left(\nabla^2_{\widehat{x}}\mathfrak{b}\right)\,d\widehat{x}\,dt +\displaystyle\int_0^T\int_{\Omega_p}\mathbb{K}_{33}(x_3)\partial_{3}p\,\partial_{3}  \varphi dx\,dt=0,\\&  \forall  \varphi \in L^2(0,T;M)\cap H^1(0,T;L^2(\Omega_p)) \textrm{ such that } \varphi(T)=0. 
\end{split}\end{equation}
\end{teorem}

 \begin{proof}
 
 To obtain 
 the limit equations, we choose the test function $\vect{v}^h\varphi$ in \eqref{varformporous2gamma0may10}, where $\varphi \in  C^1([0,T])$ such that $\varphi(T)=0$ and $\vect{v}^h$ is defined with \footnote{\CCC As usual the test function is chosen to accommodate compactness result: $h \vect{\theta}^h$ satisfies classical Kirchoff-Love ansatz,  $h\ee\vect{\zeta}\left(x,\frac{\widehat{x}}{\ee},\frac{x_3}{\frac{\ee}{h}}\right)+h^2\displaystyle\int_{0}^{x_3}\vect{r}(x)\,dx_3 $ is the corrector on the elastic matrix and $h^2\vect{\xi}\left(x,\frac{\widehat{x}}{\ee},\frac{x_3}{\frac{\ee}{h}}\right)$ is the corrector on the fluid part. \BBB} \CCC 

\begin{equation}
\vect{v}^{h}(x)=\displaystyle h\vect{\theta}^{h}(x)+h\ee\vect{\zeta}\left(x,\frac{\widehat{x}}{\ee},\frac{x_3}{\frac{\ee}{h}}\right)+h^2\vect{\xi}\left(x,\frac{\widehat{x}}{\ee},\frac{x_3}{\frac{\ee}{h}}\right)+h^2\displaystyle\int_{0}^{x_3}\vect{r}(x)\,dx_3,
\end{equation} 
Here
\begin{equation}
\vect{\theta}^{h}(x)= \left(\theta_1(\widehat{x})-x_3\partial_1\theta_3(\widehat{x}), \theta_2(\widehat{x})-x_3\partial_2\theta_3(\widehat{x}),  h^{-1}\theta_3(\widehat{x}) \right)^T,
\end{equation}
$\vect{\theta}_{*}\in \dot{C}^{1}_{\#}(\omega;\RR^2)$, ${\theta}_3\in \dot{C}^2_{\#}(\omega)$, $\vect{\zeta}\in C^{1}_{c}(\Omega;C^{1}(\mathcal{Y};\RR^3))$, $\vect{\xi}\in C^{1}_{c}(\Omega;C^{1}_{\#}(\mathcal{Y};\RR^3))$, $\vect{\xi}(x,y):=\sum_{i=1}^m\vect{\tau}^i\left(y\right)\xi^i(x)$, is such that $\xi^i \in C_c^1 (U_i)$, $\textrm{supp } \xi^i\subset U_i$,  $\vect{\tau}^i \in C^1(\mathcal{Y};\RR^3)$, $\text{supp } \vect{\tau}^i \subset \mathcal{Y}_f^i$, $\textrm{div}_y \vect{\tau}^i=0$,\,
for $i=1,\dots,m$ and $\vect{r}\in C^1_{c}(\Omega)^3$. 
Note that 
\begin{eqnarray*} 
 e_h(v^h)(x)&=&h\left[\iota\left(e_{\hat{x}}(\vect{\theta}_*)(\hat{x})-x_3 \nabla^2_{\hat{x}}\theta_3 (\hat{x})\right)+\mathfrak{C}_{\infty}(\vect{\zeta},\vect{r})\left(x,\frac{\widehat{x}}{\ee},\frac{x_3}{\frac{\ee}{h}}\right)\right]+\frac{h^2}{\eps} e_y(\vect{\xi})\left(x,\frac{\widehat{x}}{\ee},\frac{x_3}{\frac{\ee}{h}}\right)\\ & &+h\sum_{i=1}^3\partial_3 \xi_i \left(x,\frac{\widehat{x}}{\ee},\frac{x_3}{\frac{\ee}{h}}\right)\vect{e}^i\odot \vect{e}^3+O(\max\{\ee,h^2\}),\end{eqnarray*}
 where $\|O(\max\{\ee,h^2\})\|_{L^{\infty{}}}\leq C\max\{\ee,h^2\}. $
Then as a consequence of \eqref{depra}, Proposition \ref{prop 2may10} and Theorem \ref{gammazeromay11} we have
\CCC 
\begin{equation} \label{limeq1*may11}      -\eta \int_0^T\int_{\Omega}\kappa^{h}\partial_t \vect{u}^h(t)\vect{v}^{h}(x)\varphi'(t)\,dx\,dt  \rightarrow 0,
\end{equation}
\BBB
\begin{equation}\label{limeq2*may11}
    2\displaystyle\frac{\ee}{h^2}\int_0^T\int_{\Omega_{f}^h} e_h\left(\partial_t \vect{u}^h\right): \frac{\ee}{h^2}e_h\left(\vect{v}^{h}\right)\varphi\,dx\,dt
\rightarrow \displaystyle 2\int_0^T\int_{\Omega}\int_{\mathcal{Y}_f}e_{y}\left(\partial_t \vect{u}^{0}_f\right): e_{y}(\vect{\xi})\varphi dydxdt,
\end{equation} 
\begin{equation}\label{limeq3*may11}
    \begin{split}
        &\displaystyle\int_0^T \int_{\Omega_s^h}\mathbb{A}^h(x) \frac{1}{h}e_h(\widehat{\vect{u}}^{h}):\CCC \frac{1}{h}e_h (\vect{v}^h)\varphi\, dx\, dt\\ \BBB
& \rightarrow \displaystyle\int_0^T\int_{\Omega}\int_{\mathcal{Y}_s(x_3)}\mathbb{A}(x_3,y) \left[\iota (e_{\widehat{x}}(\mathfrak{a})-x_3\nabla^2_{\widehat{x}}\mathfrak{b})+\mathfrak{C}_{\infty}(\vect{w}, \vect{g})\right]:  \left[\iota (e_{\widehat{x}}(\vect{\theta}_*)-x_3\nabla^2_{\widehat{x}}\theta_3)+\mathfrak{C}_{\infty}(\vect{\zeta},\vect{r})\right]\varphi\,dy\,dx\,dt,
    \end{split}
\end{equation} 
\CCC
\begin{equation}\label{limeq4*may11}
\begin{split}
    &-_{(H_T^1L^2_f)'} \Big \langle \frac{p^h}{h}, \frac{1}{h}\text{div}_h\,\vect{v}^h\varphi\Big \rangle_{H_T^1L^2_f} \\
&\rightarrow\displaystyle -_{(H_T^1L^2(\Omega \times \mathcal{Y} ))'} \Big \langle \chi_{\Omega \times \mathcal{Y}^{x_3}_f}p,\left(\text{div}_{\widehat{x}}\,\left(\vect{\theta}_*-x_3\nabla_{\widehat{x}}\theta_3\right)+\text{div}_{y}\, \vect{\zeta}+r_3+\partial_3\xi_3\right)\varphi\Big \rangle_{H_T^1L^2(\Omega \times \mathcal{Y} )} \\
\end{split}
\end{equation}
\BBB
\begin{equation}\label{limeq5*may11} \CCC
    \begin{split}
        &\displaystyle\int_0^T\int_{\Omega}(\pi_{1/h} \vect{F}^h)\cdot \vect{v}^h \varphi\,dx\,dt\\
&\rightarrow \displaystyle\int_0^T\int_{\omega} \overline{\langle \vect{F}\rangle_{Y}}\cdot(\vect{\theta}_*,\theta_3)\varphi\, d\widehat{x}\,dt-\displaystyle\int_0^T\int_{\omega} \overline{\langle x_3\vect{F}_{*}\rangle_{\mathcal{Y}}}\cdot\nabla_{\widehat{x}}\theta_3\varphi\, d\widehat{x}\, dt.
    \end{split}
\end{equation}
\CCC 
By taking $\vect{\theta}_*=0$, $\theta_3=0$, $\vect{\xi}=0$ we conclude from \eqref{varformporous2gamma0may10} that:  
\begin{equation}\label{homoperatormarch30add}
	\begin{split} &\displaystyle\int_0^T\int_{\Omega}\int_{\mathcal{Y}_s(x_3)}\mathbb{A}(x_3,y) \left[\iota (e_{\widehat{x}}(\mathfrak{a})-x_3\nabla^2_{\widehat{x}}\mathfrak{b})+\mathfrak{C}_{\infty}(\vect{w}, \vect{g})\right]: \mathfrak{C}_{\infty}(\vect{\zeta}, \vect{r})\varphi\,dy\,dx\, dt\\ & \hspace{+2ex} -\displaystyle_{(H_T^1L^2(\Omega \times \mathcal{Y}))'}\langle \chi_{\Omega \times \mathcal{Y}_f^{x_3}}p, \left(\text{div}_{y}\, \vect{\zeta}+r_3\right) \varphi\rangle_{H_T^1L^2(\Omega \times \mathcal{Y})}=0,\\ &
		\quad \forall (\vect{\zeta},\vect{r})\in L^2(\Omega, \dot{H}^1(\mathcal{Y};\RR^3))\times L^2(\Omega;\RR^3), \varphi \in H^1(0,T) \textrm{ such that } \varphi(T)=0.
	\end{split}
\end{equation}
\CCC
From \eqref{homoperatormarch30add} we conclude that
 (see \eqref{limitequationhom1}, \eqref{eq unique uationapril16}, \eqref{cellproblem1april1},\eqref{eq:abril1})\footnote{E.g. one can take $\vect{\zeta}=\zeta_1(x) \vect{\zeta}_2(y)$, $\vect{r}=\vect{\zeta}_1(x)r_1 $ where $\zeta_1 \in L^2(\Omega)$, $\vect{\zeta}_2 \in \dot{H}^1(\mathcal{Y};\R^3)$, $r_1 \in \R$ in \eqref{homoperatormarch30add} and use linearity property to decompose the solution.}
\begin{align}
	\left(\vect{w},\vect{g}\right)(x,y,t)&\nonumber=p (x,t)\left(\widetilde{\vect{w}}_{2,x_3}(y),\widetilde{\vect{g}}_{2,x_3}\right) \\
	&\label{eq:abril123}\quad+\displaystyle\sum_{\alpha,\beta=1,2}\left[\left(e_{\widehat{x}}(\vect{a}(t,\widehat{x}))\right)_{\alpha\beta}\left(\vect{w}^{\alpha\beta}_{x_3}(y),\vect{g}^{\alpha\beta}_{x_3}\right)-x_3\partial_{\alpha\beta}\mathfrak{b}(t,\widehat{x})\left(\vect{w}^{\alpha\beta}_{x_3}(y),\vect{g}^{\alpha\beta}_{x_3}\right)\right],
\end{align}
from which it follows that  $p \in L^{\infty}(0,T;L^2(\Omega))$.
\BBB
From this, \CCC using \eqref{varformporous2gamma0may10} and again  convergences in the expression \eqref{limeq5*may11} and taking now firstly $\vect{\zeta}=\vect{r}=\vect{\xi}=0$ and secondly $\vect{\theta}_*=\theta_3=\vect{\zeta}=\vect{r}=0 $ \BBB we conclude for a.e. $t \in[0,T]$:
\begin{equation}\label{wandgmarch300}
    \begin{split}       &\displaystyle\int_{\Omega}\int_{\mathcal{Y}_s(x_3)}\mathbb{A}(x_3,y) \left[\iota (e_{\widehat{x}}(\mathfrak{a})-x_3\nabla^2_{\widehat{x}}\mathfrak{b})+\mathfrak{C}_{\infty}(\vect{w}, \vect{g})\right]:\left[\iota (e_{\widehat{x}}(\vect{\theta}_{*})-x_3\nabla^2_{\widehat{x}}\theta_3)\right]\,dy\,dx\\
&\quad +\int_{\Omega}|\mathcal{Y}_f(x_3)|x_3 p\,(\text{div}_{\widehat{x}}\nabla_{\widehat{x}}\theta_3)\,d\widehat{x}-\displaystyle\int_{\Omega}|\mathcal{Y}_f(x_3)|p\,\text{div}_{\widehat{x}}\left(\vect{\theta}_*\right)\,dx\\
&=\displaystyle\int_{\omega} \overline{\langle \vect{F}\rangle_{Y}}\cdot(\vect{\theta}_*,\theta_3)\, d\widehat{x}-\displaystyle\int_{\omega} \overline{\langle x_3\vect{F}_{*}\rangle_{Y}}\cdot\nabla_{\widehat{x}}\theta_3\, d\widehat{x}, \quad\forall \left(\vect{\theta}_{*},{\theta}_3\right) \in L_{\bar{\kappa},0},
    \end{split}
\end{equation}
\begin{align}
    &\displaystyle 2\int_{\Omega}\int_{\mathcal{Y}_f(x_3)}e_{y}\left(\partial_t \vect{u}^{0}_f\right): e_{y}(\vect{\xi})dydx\displaystyle -\displaystyle\int_{\Omega}\int_{\mathcal{Y}_f(x_3)}p\,\partial_3\xi_3\,dx\,dy=0,\quad \forall \vect{\xi}\in L^2(\omega; H_{0} ^1(I_i;H^1(\mathcal{Y};\RR^3))),  \\
&i=1,\dots,m, \text{ such that }\,\, \vect{\xi}=0\quad\text{on}\,\,\Omega\times\mathcal{Y}^{x_3}_s\quad\text{and}\quad \text{div}_y\,\vect{\xi}=0\quad\text{on}\quad \Omega\times\mathcal{Y}^{x_3}_f. \label{fluidpressureeqmarch30}
\end{align}

By substituting (\ref{eq:abril123}) into \eqref{wandgmarch300} \CCC and using \eqref{gassmantensorapril01} and \eqref{btensorapril1} \BBB we obtain \eqref{diferformpen1*may11}. 
\CCC

In order to conclude that $p \in L^2(0,T;M)$ we will use 
\eqref{fluidpressureeqmarch30} and an additinal observation. 
Firstly we conclude that $p \in L^2(0,T;L^2(\omega;H^1(I_i)))$ for every $I_i \subset J_K$ (recall \eqref{intervals}).  This can be done directly from \eqref{fluidpressureeqmarch30}. Namely,  
by using  
Remark \ref{pojednost} we take $\vect{\tau} \in  C_c^{\infty}(\mathcal{Y}_f^i;\R^3) $,  such that $\textrm{div}_y  \vect{\tau}=0$ and $\int_Y \tau_3\neq 0$ and define $\vect{\xi}(x,y)=\xi^i(x)\vect{\tau} (y)  $ where $\xi^i \in C_c^1(U_i)$. The claim follows by taking $\vect{\xi}$ as a test function in \eqref{fluidpressureeqmarch30}.

To obtain that $p \in L^2(0,T;M)$ it is enough to prove $p \in L^2(0,T;L^2(\omega;H^1(J_K))$, i.e. that we have the continuity on the interface  between $U_i$ and $U_{i+1}$, for $i=1,\dots,m-1$, where there is a flow between $U_i$ and $U_{i+1}$. 
This cannot be concluded from \eqref{fluidpressureeqmarch30}, since the test functions are zero in the neighborhood of enpoint of every $I_i$, $i=1,\dots,m$.  

Again, by using Lemma \ref{lmtestfunc} 2
we take  functions $\vect{\tau}^i,\vect{\tau}^{i+1}$ 
such that $\vect{\tau}^i \in C_c^{\infty} (\mathcal{Y}_f^i;\R^3)$, $\vect{\tau}^{i+1} \in C_c^{\infty} (\mathcal{Y}_f^{i+1};\R^3)$ and 
that  for extensions  by zero on whole $\mathcal{Y}$ 
of these functions
(we denote these extensions also by $\vect{\tau}^i$ and $\vect{\tau}^{i+1}$) we have 
\begin{equation} 
\label{osveta1} \diver_y\vect{\tau}^{i}=\diver_y\vect{\tau}^{i+1}=0, \quad \int_{\mathcal{Y}_f^i} \tau_3^i\,dy=\int_{\mathcal{Y}_f^{i+1}} \tau_3^{i+1}\,dy\neq0; 
\end{equation} 
\begin{equation} 
\label{pecenje1}  
	\vect{\tau}^i=\vect{\tau}^{i+1} \quad \textrm{ on } \{y_3=0\}. 
\end{equation}	
We also take $\xi \in C_c^1(\omega\times( d_{i+1}-\delta,d_{i+1}+\delta))$ for $\delta>0$ small enough (recall \eqref{intervals}).  	  
We define (recall \eqref{intervals2} and \eqref{intervals3}):
\begin{equation}\label{deftestsub} 
\vect{\xi}^h(x)=\chi_{(d_i^h,d_{i+1}^h]}(x_3)\xi(x)\vect{\tau}^i\left(\frac{\widehat{x}}{\ee},\frac{x_3}{\frac{\ee}{h}}\right)+\chi_{(d_{i+1}^h,d_{i+2}^h]}(x_3)\xi(x)\vect{\tau}^{i+1}\left(\frac{\widehat{x}}{\ee},\frac{x_3}{\frac{\ee}{h}}\right).
\end{equation} 	
As a consequence of \eqref{pecenje1} we have that 
$\vect{\xi}^h \in W^{1,\infty} (\Omega;\R^3)$.
Moreover we easily conclude that 
\begin{align} \label{konnn1}  
	\eps \nabla_h \vect{\xi}^h &\xrightarrow{2-r}\chi_{(d_i,d_{i+1}]}(x_3)\xi(x)\nabla_y \vect{\tau}^i(y)+\chi_{(d_{i+1},d_{i+2}]}(x_3)\xi(x)\nabla_y \vect{\tau}^{i+1}(y), \\
	h\diver_h \vect{\xi}^h &\xrightarrow{2-r}  \chi_{(d_i,d_{i+1}]}(x_3)\partial_3\xi(x) \tau^i_3(y)+\chi_{(d_{i+1},d_{i+2}]}(x_3)\partial_3 \xi(x) \tau^{i+1}_3(y)
\end{align} 	
Next we plug in $(\ref{varformporous2gamma0may10})$
the function $h^2 \vect{\xi}^h(x)\varphi(t) $,
where $\varphi \in  C^1([0,T])$ satisfies $\varphi(T)=0$. As in the case for \eqref{fluidpressureeqmarch30}, using \eqref{osveta1}, we conclude that 
\begin{equation}\label{fluidpressureeqmarch3000}
	\begin{split}
		& 2\int_{\omega\times (d_i,d_{i+1})}\xi  \int_{\mathcal{Y}_f^i}e_{y}\left(\partial_t \vect{u}^{0}_f\right): e_{y}(\vect{\tau}^i)\,dy\,dx\\
		&+ 2\int_{\omega\times (d_{i+1},d_{i+2})}\xi   \int_{\mathcal{Y}_f^{i+1}} e_{y}\left(\partial_t \vect{u}^{0}_f\right): e_{y}(\vect{\tau}^{i+1})\,dy\,dx
		\\&
		 -\int_{\omega\times (d_i,d_{i+2})}p \partial_3 \xi \int_{\mathcal{Y}_f^1} \tau^i_3(y)\,dy\,dx=0.  
	\end{split}
\end{equation}
From this we easily conclude that $\partial_3 p \in L^2(0,T;L^2(\omega \times (d_{i+1}-\delta,d_{i+1}+\delta)))$. Together with previous conclusion,  this implies that $p \in L^2(0,T;L^2(\omega;H^1(J_K))$.

\BBB
\CCC
Note that for given $p \in L^2(0,T;M)$, the solution of \eqref{fluidpressureeqmarch30} that satisfies \eqref{kreso1} is unique, which can be easily checked by subtraction. 
If $x_3 \in I_i \subset J_p \backslash J_k$, for some $i=1,\dots,m$, we easily conclude from 
\eqref{fluidpressureeqmarch30} and Lemma \ref{lmtestfunc} 1 that  $ \vect{u}_f^0(x,y,t)=0$. 
This can be seen by taking in \eqref{fluidpressureeqmarch30} the test functions of the form 
$\vect{\xi}(\hat{x},x_3,y)=\xi_1(\hat{x}) \xi_2(x_3)\vect{\tau}(y)$, where $\xi_1 \in L^2(\omega)$, $\xi_2 \in H_0^1(I_i)$ and $\vect{\tau} \in H_0^1(\mathcal{Y}_f^i;\R^3)$ such that $\textrm{div}_y \vect{\tau}=0 $. 

For $x_3 \in J_K$ the differential form of (\ref{fluidpressureeqmarch30}) is (for fixed $x \in \Omega$)
\begin{equation}\label{diferformfluid1march30}
-\Delta_{y}\partial_t \vect{u}^0_f+(0,0,\partial_3 p)+\nabla_{y}\pi=(0,0,0),
\end{equation}
which we rewrite as \footnote{\CCC Recall that $\partial_t \vect{u}_f^0$ belongs to $H^1(\mathcal{Y};\R^3)$, when extended by zero outside $\mathcal{Y}_f(x_3)$.  \BBB} 
\begin{equation}\label{diferformfluid1march301}
\begin{cases}
-\Delta_{y}\displaystyle\partial_t \vect{u}^0_f+\nabla_{y}\pi=-(0,0,\partial_3 p), &\\
\text{div}_{y}\displaystyle\partial_t \vect{u}^0_f=0, \quad
\left.\vect{u}^0_f\right|_{\partial\mathcal{Y}_f(x_3)}=0.  &
\end{cases}
\end{equation}
We conclude that, for $x_3 \in J_K$, $\partial_{t}\vect{u}^{0}_f$ and $\pi$ have a representation:
\begin{equation}\label{diferformfluid2march302}
\begin{cases}
\displaystyle\partial_t \vect{u}^0_f(x,y,t)=-\vect{q}_{x_3}^{3}(y)\partial_3 p(x,t), &\\
\displaystyle\pi=-\pi^{3}_{\gamma}(y)\partial_3 p(x,t), &
\end{cases}
\end{equation}
where $\vect{q}^3_{x_3}$ is defined in \eqref{diferformfluid2march301} and extended by zero outside $\mathcal{Y}_f(x_3)$. 
Note that the first equation in \eqref{diferformfluid2march302} is also valid for $x_3 \in J_p\backslash J_K$, since $\vect{q}_{x_3}^3=0$ there (cf. proof of Proposition \ref{tensorprop} 3).

To prove \eqref{diferformpen0*may11} we use \eqref{cr4gamma0}. We take  $\varphi\in H^1(0,T;L^2(\omega;H^1(J_p))$, $\varphi(T)=0$, and use $\partial_t \varphi(t)$ as a test function in \eqref{cr4gamma0} for every $t \in (0,T)$ and integrate over the interval $(0,T)$. We then integrate by parts the previous to the last term on the left hand side of \eqref{cr4gamma0} and use \eqref{eq:abril123}, \eqref{diferformfluid2march302} and  Proposition \ref{tensorprop} 2 \& 3 as well as the density argument.

 \end{proof}
 \begin{remarkica} \label{remdarcy}
 	The first equation in \eqref{diferformfluid2march302} tells us that in the effective Darcy's law only the derivative of the pressure in the vertical direction is present (unlike in the bulk model where the full gradient is present). By integrating the first equation in \eqref{diferformfluid2march302} over $\mathcal{Y}$ we obtain that the mean fluid velocity is given by 
  $$\int_{\mathcal{Y}_f}\partial_t \vect{u}_f^0(\cdot,y)\, dy=-\partial_3 p\int_{\mathcal{Y}_f} \vect{q}_x^3(y) \, dy.  $$
 \end{remarkica} 		

\begin{remarkica} \label{remgotovo} \CCC
In order to conclude that pressure belongs to $ L^2(0,T; L^2(\omega; H^1(J_K))$ we needed to argument the continuity of traces over the interface where there is a flow (see Section \ref{geometry}). This was done by using test functions that are supported in the neighbourhood of the interface. Such test functions need to be  be weakly differentiable, to have support on the fluid part  \CCC and to satisfy the condition that the third component has non-zero mean value over the torus.  Note that, as a consequence of Lemma \ref{lmtestfunc} 1, it is not possible that the support  of such fastly oscillating test function doesn't intersect the part of $\mathcal{Y}$ given by  $\{y_3=0\}$. Thus we needed to be be careful to constructing such test functions at the interface, since the changing of the domain of the fluid required that we change them from one region to the other, while keeping the continuity property. 

The conclusion that the limit pressure  $p$ belongs to the space $L^2(0,T;M)$ has the following consequences (cf. Remark \ref{marindod1}): 
\begin{enumerate}
\item on the interface between $U_i$ and $U_{i+1}$ that satisfies the property that there is no flow, the trace of the pressure may have jump. However from the space of test functions we conclude \CCC that $\partial_3 p=0$ on the interface (Neumann boundary condition); 
\item on the interface between $U_i$ and $U_{i+1}$ that satisfies  the property that there is flow, the trace of the pressure is continuous. 
\end{enumerate} 
\BBB
\end{remarkica}
\CCC
\begin{remarkica} \label{usporedba} 
Here we compare our model with the ones obtained in \cite{taber,Mikelic2015}. In \cite{taber} the author derives the flexural (bending) plate equations (i.e., the equations for $\mathfrak{b}$), starting from 3D  Biot's equations (using isotropic elasticity and assuming the force term being zero, but with some boundary conditions), under the following assumptions: 
\begin{enumerate} 
\item Normals to the middle surface of the solid skeleton ($x_3=0$) remain straight and normal during deformation (this is equivalent to assuming Kirchoff-Love ansatz, cf. \eqref{crgamma0} and see \cite{ciarlet2}). 
\item The plate is in a state of approximate plane stress.
\item In plane fluid-velocity gradients relative to the solid are small compared to the transverse fluid-velocity gradient (this is eqzivalent to assuming that the term $\partial_{3} p$ is dominant in $\nabla p$ ). 
\end{enumerate} 	 
The equations obtained in this way are then justified in \cite{Mikelic2015} from 3D Biot's (quasi-static) theory (again assuming isotropic elasticity) by using the approach from \cite{ciarlet2}.  Moreover, in \cite{Mikelic2015} the  membrane equations (for the part of in-plane components $\mathfrak{a}$) are obtained and they are decoupled from the bending equations (the equations for $\mathfrak{b}$). 

Note that in the case when $M_0$, $|\mathcal{Y}_f|$ and $\mathbb{B}^H$ are $x_3$ independent. one can also separate $\mathfrak{a}$ and $\mathfrak{b}$ appearing in \eqref{diferformpen0*may11} by using test functions that are independent of $x_3$ and the ones that are $x_3$ dependent and satisfy $\int_I\varphi \,dx_3=0$ for every $\hat{x} \in \omega$ (this is how \eqref{diferformpen0*may11} is written  in \cite{Mikelic2015}). 

It is not unusual that non-decoupling of membrane and bending plate equations doesn't happen in our case, when there are heterogeneities across thickness. This is also obtained in \cite{mbukalvelcic2017} and in \cite{Marinvelciczubrinic2022} in some regimes in a different context.   
\end{remarkica} 
 \BBB
 \subsection{Analysis of the limit equations} \label{analysis} 
We will slightly modify the system \eqref{diferformpen1*may11}-\eqref{diferformpen0*may11}.
 For  $\vect{F} \in H^1(0,T;L_{1,0}')$, $\vect{G} \in L^2(0,T;M')$ and $t_0 \in L^2(\Omega)$ we will analyze the following problem:  find  $(\mathfrak{a},\mathfrak{b}, p) \in L^2(0,T;L_{1,0})  \times L^2(0,T;M)$ that satisfy
\begin{equation}\label{diferformpen1*may111}\begin{split}
    &\displaystyle\int_{\omega}\mathbb{A}^{\rm hom}(e_{\widehat{x}}(\mathfrak{a}), \nabla^2_{\widehat{x}}\mathfrak{b}):(e_{\widehat{x}}(\vect{\theta}_*), \nabla^2_{\widehat{x}}\theta_3)d\widehat{x}-\displaystyle\int_{\omega}\int_{J_p}\left(|\mathcal{Y}_f(x_3)| \mathbb{I}-\mathbb{B}^{H}(x_3)\right)p \, dx_3:\left[\iota (e_{\widehat{x}}(\vect{\theta}_{*}))\right]\,d\widehat{x}
    \\&+\displaystyle\int_{\omega} \int_{J_p}\left(|\mathcal{Y}_f(x_3)| \mathbb{I}-\mathbb{B}^{H}(x_3)\right)x_3 p \, dx_3:\left[\iota (\nabla^2_{\widehat{x}}\theta_3)\right]\,d\widehat{x} =\displaystyle _{L_{1,0}'}\langle\vect{F}(t), (\vect{\theta}_{*},{\theta}_3) \rangle_{L_{1,0}}, 
    \\& \quad\forall \left(\vect{\theta}_{*},{\theta}_3\right) \in L_{1,0},  \textrm{for a.e. } t \in (0,T), 
\end{split}\end{equation}
\begin{equation}\label{diferformpen0*may111}\begin{split}
    &\displaystyle-\int_0^T\int_{\Omega}M_0(x_3) p\, \partial_t \varphi dx\,dt-\displaystyle\int_0^T\int_{\omega}\int_{J_p}(|\mathcal{Y}_{f}(x_3)| \mathbb{I}-\mathbb{B}^{H}(x_3))\partial_t\varphi \, dx_3:\iota\left(e_{\widehat{x}}\left(\mathfrak{a}\right)\right) d\widehat{x}\,dt
    \\&+\displaystyle\int_0^T\int_{\omega}\int_{J_p}\left(|\mathcal{Y}_f(x_3)|\mathbb{I}-\mathbb{B}^H(x_3)\right)x_3\partial_t\varphi \, dx_3:\iota\left(\nabla^2_{\widehat{x}}\mathfrak{b}\right)\,d\widehat{x}\,dt +\displaystyle\int_0^T\int_{\Omega_p}\mathbb{K}_{33}(x_3)\partial_{3}p\,\partial_{3}  \varphi dx\,dt  
    \\& =\int_0^T {_{M'}}\langle \vect{G}, \varphi \rangle_{M}\,dt+ \int_{\Omega_p} t_0 \cdot \varphi(0) \,dx, \quad  \forall  \varphi \in L^2(0,T;M)\cap H^1(0,T;L^2(\Omega_p)) \textrm{ such that } \varphi(T)=0. 
\end{split}\end{equation}
Thus the loads $\overline{\langle \vect{F}\rangle_{Y}}$, $\overline{\langle x_3\vect{F}_*\rangle_{\mathcal{Y}}}$ for $\vect{F} \in H^1(0,T;L^2(\Omega \times \mathcal{Y};\mathbb{R}^3))$ of the system \eqref{diferformpen1*may11}-\eqref{diferformpen0*may11} are naturally replaced with $\vect{F} \in H^1(0,T;L_{1,0}')$, and we have additional terms on the right hand side of \eqref{diferformpen0*may111} when compared with \eqref{diferformpen0*may11}. 
\CCC We have also replaced $L_{{\kappa},0}$ by $L_{1,0}$. This can be done without loss of generality, since the solution in space $L_{\kappa,0}$ can be easily obtained from the solution in space $L_{1,0}$ by a translation for every $t \in [0,T]$. \BBB
 We will give the abstract framework for this problem.
  To this end, we adapt the approach from \cite{gurvich}.
  As we will see, adding the additional term on the right hand side of \eqref{diferformpen0*may11} will change the initial condition for the pressure.

 Let $V$ be a separable Hilbert space with dual $V'$ (which is not identified with $V$ here). Assume that $V$ is densely and continuously embedded into another Hilbert space $H$, which is identified with its dual: 
 $$V\hookrightarrow H \equiv H' \hookrightarrow V'.   $$
 Consequently $H$ is densely and continuously embedded into $V'$. 
We denote the inner product on $H$ by $\langle\cdot,\cdot\rangle_H$ and the norm induced by that scalar product by $\|\cdot\|_H$. We also denote by $_{V'}\langle \cdot,\cdot\rangle_V$ the duality pairing between $V$ and $V'$ and by $\| \cdot\|_V$ the norm on $V$. Suppose that $\mathcal{B}_0 \in L(H)$ and $\mathcal{A}_0 \in L(V,V')$.  Let $u_0 \in H$ and $S \in L^2(0,T;V')$ be given.
We consider the following Cauchy problem:  Find $u \in L^2(0,T;V)$ such that
\begin{equation}  \label{abspbm}
\begin{cases} 
\frac{d}{dt} [\mathcal{B}_0u]	+\mathcal{A}_0 u=S \in L^2(0,T;V') \\ 
[\mathcal{B}_0u](0)=\CCC \mathcal{B}_0\BBB u_0 \in H. 
\end{cases} 	   
\end{equation} 
The time derivative in \eqref{abspbm} is taken in the sense of distributions. 
We give the following definition of the weak solution of \eqref{abspbm} on the time interval $(0,T)$ (\CCC see \cite{gurvich}). \BBB 
\begin{definition}\label{absdef} 
	The function $u \in L^2(0,T;V)$ such that 
	\begin{equation} \label{krnj10} 
	-\int_0^T \langle \mathcal{B}_0 u(t), v'(t) \rangle_H\, dt+\int_0^T {_{V'}} \langle \mathcal{A}_0 u(t),v(t) \rangle_V\, dt= \int_0^T {_{V'}} \langle S(t), v(t) \rangle_V \, dt+\langle \mathcal{B}_0 u_0, v(0) \rangle_H, 
	\end{equation} 
	holds for all $ v \in \{w\in L^2(0,T;V) \cap H^1(0,T;H): w(T)=0\} $
	is called a weak solution to \eqref{abspbm}. 	
\end{definition}
\begin{remarkica} \label{refkrnj1} 
Notice that the solution of \eqref{abspbm} satisfies $\mathcal{B}_0u \in H^1(0,T;V')$ and  thus  we can give a meaning to the initial value in \eqref{abspbm}. 	
	\end{remarkica} 
	The following assumption will be needed for the existence result.
\begin{assumption} \label{absassum} We assume that 
	\begin{enumerate} 
		\item 
		$\mathcal{A}_0 \textrm{ is monotone on } V$, i.e.
		$  _{V'}\langle \mathcal{A}_0 v,v\rangle_V \geq 0, \quad \forall v \in V; $ 
		\item $\mathcal{B}_0$ is self-adjoint positive semidefinite on $H$;
		\item There exist constants $\lambda,c>0$ such that
		\begin{equation} \label{krnj20} _{V'} \langle \mathcal{A}_0 v,v\rangle_V+\lambda \langle \mathcal{B}_0v,v\rangle_H \geq c\|v\|_V^2, \quad \forall v \in V.    \end{equation}   
	\end{enumerate} 	 
\end{assumption} 
We have the following theorem as a consequence of \cite[Theorem 2.1, Theorem 2.2]{gurvich}. \CCC We will use this theorem to prove the existence and uniqueness result. \BBB
\begin{teorem} \label{absteo} 
Under the Assumption \ref{absassum}, there exists a unique solution of \eqref{abspbm} in the sense of Definition \ref{absdef}. The solution satisfies the following stability estimate: 
\begin{equation} \label{stab1} 
 \|u\|^2_{L^2(0,T;V)} \leq C(c,\lambda) \left[\|S\|^2_{L^2(0,T;V')}+\langle\mathcal{B}_0u_0,u_0\rangle_H \right], 
 \end{equation} 
where $C(c,\lambda)>0$ depends only on $c$ and $\lambda$. 
\end{teorem} 
\begin{remarkica} \label{remkrnj1} 
If, in addition, we have that $\mathcal{B}_0$ is coercive on $H$, then it follows from \cite[Remark 1, Chapter XVIII.5]{dautraylions}	that the solution of \eqref{krnj10} is in $C([0,T];H)$ and there exists $C(c,\lambda)>0$ such that  
\begin{equation}  \label{stab11} 
	\|u\|^2_{C(0,T;H)} \leq C(c,\lambda) \left[\|S\|^2_{L^2(0,T;V')}+\langle\mathcal{B}_0u_0,u_0\rangle_H \right].
\end{equation} 	
\end{remarkica} 	
Next we state and prove the main result of this section. 
\begin{teorem} \label{existeo} 
For $\vect{F} \in H^1(0,T;L_{1,0}')$,  $\vect{G} \in L^2(0,T;M')$  there exist a unique solution $(\mathfrak{a},\mathfrak{b}, p) \in L^2(0,T;L_{1,0}) \times L^2(0,T;M)$ of \eqref{diferformpen1*may111}- \eqref{diferformpen0*may111}.
	\end{teorem} 
\begin{proof} 
\CCC We will put our problem in the framework of Theorem \ref{absteo}. \BBB
For $\widehat{p} \in L^2(\Omega)$ and $\widehat{\vect{F}} \in L_{1,0}'$ we introduce $(\mathfrak{a}^{\widehat{p}}_{0,1},\mathfrak{b}^{\widehat{p}}_{0,1}) \in L_{1,0}$ and $(\mathfrak{a}_{0,2}^{\widehat{\vect{F}}},\mathfrak{b}_{0,2}^{\widehat{\vect{F}}}) \in L_{1,0}$ respectively that satisfy: 
\begin{align}
      \nonumber &  \displaystyle\int_{\omega}\mathbb{A}^{\rm hom}(e_{\widehat{x}}(\mathfrak{a}_{0,1}^{\widehat{p}}), \nabla^2_{\widehat{x}}\mathfrak{b}_{0,1}^{\widehat{p}}):(e_{\widehat{x}}(\vect{\theta}_*), \nabla^2_{\widehat{x}}\theta_3)d\widehat{x}=\displaystyle\int_{\omega}\int_{J_p}\left(|\mathcal{Y}_f(x_3)| \mathbb{I}-\mathbb{B}^{H}(x_3)\right)\widehat{p}\, dx_3:\left[\iota (e_{\widehat{x}}(\vect{\theta}_{*}))\right]\,d\widehat{x}\\ \label{defap}
&-\displaystyle\int_{\omega} \int_{J_p}\left(|\mathcal{Y}_f(x_3)| \mathbb{I}-\mathbb{B}^{H}(x_3)\right)x_3 \widehat{p} \, dx_3:\left[\iota (\nabla^2_{\widehat{x}}\theta_3)\right]\,d\widehat{x},  \quad\forall \left(\vect{\theta}_{*},{\theta}_3\right) \in L_{1,0};
\end{align}
\begin{equation} \label{defaf}
    \begin{split}
        &\displaystyle\int_{\omega}\mathbb{A}^{\rm hom}(e_{\widehat{x}}(\mathfrak{a}_{0,2}^{\widehat{\vect{F}}}), \nabla^2_{\widehat{x}}\mathfrak{b}_{0,2}^{\widehat{\vect{F}}}):(e_{\widehat{x}}(\vect{\theta}_*), \nabla^2_{\widehat{x}}\theta_3)d\widehat{x}=_{L_{1,0}'} \langle \widehat{\vect{F}},(\vect{\theta}_{*},{\theta}_3) \rangle_{L_{1,0}}, \quad\forall \left(\vect{\theta}_{*},{\theta}_3\right) \in L_{1,0}.
    \end{split}
\end{equation}
Assuming $\widehat{p}$ and $\widehat{\vect{F}}$ are known, the solutions of \eqref{defap}-\eqref{defaf} are unique as a consequence of Proposition \ref{tensorprop} and Lax-Milgram. 
\CCC By assuming $p$ known we decompose the solution of \eqref{diferformpen1*may111} as $(\mathfrak{a}^p_{0,1},\mathfrak{b}^p_{0,1})+(\mathfrak{a}^{\vect{F}}_{0,2},\mathfrak{b}^{\vect{F}}_{0,2})$.    
By plugging this in\eqref{diferformpen0*may111}, \BBB  the solution of \eqref{diferformpen0*may111} we write in the following way: 
{\allowdisplaybreaks
\begin{align}
     &-\displaystyle\int_0^T\int_{\Omega_p}M_0(x_3) p\,\partial_t \varphi dx\,dt -\displaystyle\int_0^T\int_{\omega}\int_{J_p}(|\mathcal{Y}_{f}(x_3)| \mathbb{I}-\mathbb{B}^{H}(x_3))\partial_t\varphi \, dx_3:\iota\left(e_{\widehat{x}}\left(\mathfrak{a}_{0,1}^p\right)\right) d\widehat{x}\,dt\\*
    &+\displaystyle\int_0^T\int_{\omega}\int_{J_p}\left(|\mathcal{Y}_f(x_3)|\mathbb{I}-\mathbb{B}^H(x_3)\right)x_3\partial_t\varphi \, dx_3:\iota\left(\nabla^2_{\widehat{x}}\mathfrak{b}_{0,1}^p\right)\,d\widehat{x}\,dt+\displaystyle\int_0^T\int_{\Omega_p}\mathbb{K}_{33}(x_3)\partial_{3}p\,\partial_{3} \varphi dx\,dt\\ & =
    \displaystyle\int_0^T\int_{\omega}\int_{J_p}(|\mathcal{Y}_{f}(x_3)| \mathbb{I}-\mathbb{B}^{H}(x_3))\partial_t \varphi \, dx_3:\iota\left(e_{\widehat{x}}\left(\mathfrak{a}_{0,2}^{\vect{F}}\right)\right) d\widehat{x}\,dt \\* & -\displaystyle\int_0^T\int_{\omega}\int_{J_p}\left(|\mathcal{Y}_f(x_3)|\mathbb{I}-\mathbb{B}^H(x_3)\right)x_3\partial_t \varphi \, dx_3:\iota\left(\nabla^2_{\widehat{x}}\mathfrak{b}_{0,2}^{\vect{F}}\right)\,d\widehat{x}\,dt+\int_0^T {_{M'}}\langle G,\varphi\rangle_{M}\,dt\\* \label{defp}  & + \int_{\Omega_p} t_0\cdot \varphi(0)\, dx, \quad    \forall  \varphi \in L^2(0,T;M) \cap H^1(0,T;L^2(\Omega_p)) \textrm{ such that } \varphi(T)=0.  
\end{align}
}
We put \eqref{defp} in the framework of weak solution of \eqref{abspbm}, after integration by parts on the right hand side. We now put $V:=M$, $H:=L^2(\Omega_p)$. The operator $\mathcal{A}_0$ is defined through bilinear form 
\begin{equation} \label{krnj3} 
 _{V'} \langle \mathcal{A}_0 v_1,v_2\rangle_V:=\int_{\Omega_p} \mathbb{K}_{33}(x_3) \partial_3 v_1(x)\cdot \partial_3 v_2(x) \,dx, \quad \forall v_1,v_2 \in V. 
 \end{equation} 
Obviously 1 of Assumption \ref{absassum} is satisfied. 
The operator $\mathcal{B}_0$ is also defined  through bilinear form
\begin{align}
         \langle \mathcal{B}_0 h_1, h_2 \rangle_H&:= \int_{\Omega_p} M_0(x_3) h_1(x) \cdot h_2(x) \, dx +\displaystyle\int_{\omega}\int_{J_p}(|\mathcal{Y}_{f}(x_3)| \mathbb{I}-\mathbb{B}^{H}(x_3))h_2 \, dx_3:\iota\left(e_{\widehat{x}}\left(\mathfrak{a}_{0,1}^{h_1}\right)\right) d\widehat{x}\\ & -\displaystyle\int_{\omega}\int_{J_p}\left(|\mathcal{Y}_f(x_3)|\mathbb{I}-\mathbb{B}^H(x_3)\right)x_3h_2 \, dx_3:\iota\left(\nabla^2_{\widehat{x}}\mathfrak{b}_{0,1}^{h_1}\right)\,d\widehat{x},\label{krnj2} 
\end{align}
where $(\mathfrak{a}_{0,1}^{h_1},\mathfrak{b}_{0,1}^{h_1}) \in L_{1,0}$ are defined as the solutions of \eqref{defap} once $\widehat{p}$  is replaced by $h_1$. 
From \eqref{defap} it follows: 
\begin{equation*}
    \begin{split}
        &\displaystyle\int_{\omega}\int_{J_p}(|\mathcal{Y}_{f}(x_3)| \mathbb{I}-\mathbb{B}^{H}(x_3))h_2 dx_3:\iota\left(e_{\widehat{x}}\left(\mathfrak{a}_{0,1}^{h_1}\right)\right) d\widehat{x} 
 -\displaystyle\int_{\omega}\int_{J_p}\left(|\mathcal{Y}_f(x_3)|\mathbb{I}-\mathbb{B}^H(x_3)\right)x_3h_2 dx_3:\iota\left(\nabla^2_{\widehat{x}}\mathfrak{b}_{0,1}^{h_1}\right)\,d\widehat{x} \\	& = \displaystyle\int_{\omega}\mathbb{A}^{\rm hom}(e_{\widehat{x}}(\mathfrak{a}_{0,1}^{h_1}), \nabla^2_{\widehat{x}}\mathfrak{b}_{0,1}^{h_1}):(e_{\widehat{x}}(\mathfrak{a}_{0,1}^{h_2}), \nabla^2_{\widehat{x}}\mathfrak{b}_{0,1}^{h_2})d\widehat{x}\\ &=\displaystyle\int_{\omega}\int_{J_p}(|\mathcal{Y}_{f}(x_3)| \mathbb{I}-\mathbb{B}^{H}(x_3))h_1  dx_3:\iota\left(e_{\widehat{x}}\left(\mathfrak{a}_{0,1}^{h_2}\right)\right) d\widehat{x}  
 -\displaystyle\int_{\omega}\int_{J_p}\left(|\mathcal{Y}_f(x_3)|\mathbb{I}-\mathbb{B}^H(x_3)\right)x_3h_1  dx_3:\iota\left(\nabla^2_{\widehat{x}}\mathfrak{b}_{0,1}^{h_2}\right)\,d\widehat{x}.
    \end{split}
\end{equation*}
Using this and 4 of Proposition \ref{tensorprop} we have the positive definitness of $\mathcal{B}_0$ and thus Assumption \ref{absassum} 2 is satisfied. 
Also, it is easy to see that Assumption \ref{absassum} 3 is satisfied. 
The operator $S(t) \in V'$ is defined for $t\in [0,T]$ in the following way
\begin{eqnarray*} 
 _{V'}\langle S(t),v \rangle_V &:=&
-\displaystyle\int_{\omega}\int_{J_p}(|\mathcal{Y}_{f}(x_3)| \mathbb{I}-\mathbb{B}^{H}(x_3)) v\, dx_3:\iota\left(e_{\widehat{x}}\left(\mathfrak{a}_{0,2}^{\partial_t\vect{F}(t)}\right)\right) d\widehat{x} \\ & & +\displaystyle\int_{\omega}\int_{J_p}\left(|\mathcal{Y}_f(x_3)|\mathbb{I}-\mathbb{B}^H(x_3)\right)x_3v \, dx_3:\iota\left(\nabla^2_{\widehat{x}}\mathfrak{b}_{0,2}^{\partial_t \vect{F}(t)}\right)\,d\widehat{x} +_{V'}\langle \vect{G}(t),v\rangle_V, \ \forall v \in V,
\end{eqnarray*} 
where for $t \in [0,T]$, $(\mathfrak{a}_{0,2}^{\partial_t\vect{F}(t)}, \mathfrak{b}_{0,2}^{\partial_t\vect{F}(t)}) \in L_{1,0}$ are defined  as solutions of \eqref{defaf} once $\widehat{\vect{F}}$ is replaced by $\partial_t \vect{F}(t)$. 
 The pressure
$p_0$ is defined through the following relation: 
\begin{equation}\label{krnj1} 
    \begin{split}
        & \langle \mathcal{B}_0 p_0,h \rangle_H = \displaystyle\int_{\omega}\int_{J_p}(|\mathcal{Y}_{f}(x_3)| \mathbb{I}-\mathbb{B}^{H}(x_3))h \, dx_3:\iota\left(e_{\widehat{x}}\left(\mathfrak{a}_{0,2}^{\vect{F}(0)}\right)\right) d\widehat{x} \\  & - \displaystyle\int_{\omega}\int_{J_p}\left(|\mathcal{Y}_f(x_3)|\mathbb{I}-\mathbb{B}^H(x_3)\right)x_3 h \, dx_3:\iota\left(\nabla^2_{\widehat{x}}\mathfrak{b}_{0,2}^{\vect{F}(0)}\right)\,d\widehat{x}+\int_{\Omega_p} t_0\cdot h\,dx\\ & =  \displaystyle\int_{\omega}\mathbb{A}^{\rm hom}(e_{\widehat{x}}(\mathfrak{a}_{0,1}^{h}), \nabla^2_{\widehat{x}}\mathfrak{b}_{0,1}^{h}):(e_{\widehat{x}}(\mathfrak{a}_{0,2}^{\vect{F}(0)}), \nabla^2_{\widehat{x}}\mathfrak{b}_{0,2}^{\vect{F}(0)})d\widehat{x}+\int_{\Omega_p} t_0\cdot h\,dx, \quad \forall h \in H,
    \end{split}
\end{equation}
where $(\mathfrak{a}^{\vect{F}(0)}, \mathfrak{b}^{\vect{F}(0)}) \in L_{1,0}$ is defined by \eqref{defaf}. Since $\mathcal{B}_0$ is a positive definite bounded operator this defines unique $p_0 \in H$ by Riesz representation theorem.  We apply Theorem \ref{absteo} to obtain the unique $p \in L^2(0,T;V) $. The uniqueness of $(\mathfrak{a},\mathfrak{b}) \in L^2(0,T;L_{1,0})$ follows from \eqref{diferformpen1*may11} and Lax-Milgram. 
	\end{proof} 
From \eqref{krnj1} we conclude that
the initial condition $p_0$ of the system \eqref{diferformpen1*may11}-\eqref{diferformpen0*may11} is related with the initial condition $\widetilde{p}_0$ of the system \eqref{diferformpen1*may111}-\eqref{diferformpen0*may111}, after simple substitution of the loads, in the following way
\begin{equation} \label{deftp0} 
\widetilde{p}_0=p_0+\mathcal{B}^{-1}_0t_0, 
\end{equation} 
 where $\mathcal{B}_0$ is defined by $\eqref{krnj2}$.   
The following proposition gives us additional regularity of the solution. 
\begin{propozicija} 
If $\vect{F} \in H^1(0,T;L_{1,0}')$, $\vect{G} \in L^2(0,T;V')$, $t_0 \in L^2(\Omega)$, then the solution of \eqref{diferformpen1*may111}- \eqref{diferformpen0*may111}  satisfies $(\mathfrak{a},\mathfrak{b}, p) \in C([0,T];L_{1,0})  \times C([0,T];L^2(\Omega))$.
\end{propozicija} 	
\begin{proof} 
The fact that $p \in C([0,T];L^2(\Omega))$ follows from \eqref{defp} and Remark \ref{remkrnj1}. The regularity of $\mathfrak{a}$, $\mathfrak{b}$ follows then from \eqref{diferformpen1*may111}. 
\end{proof}
\begin{remarkica} 
The previous proposition fills the gap between the result of Theorem  \ref{effectiveequationsmay11} and Theorem \ref{existeo}. Namely as a result of Theorem  \ref{effectiveequationsmay11} we have that the obtained limit is in $L^{\infty}$ w.r.t. time, while Theorem \ref{existeo} guarantees only the existence of solution, which is in $L^2$ w.r.t. time. 
\end{remarkica} 	 	
\CCC The rest of the section is devoted to proving Proposition \ref{propen}, which is crucial for proving the convergence results in Section \ref{secstrongcon}. \BBB
\CCC Proposition \ref{propreg}  will help us to prove Proposition \ref{propen}.
\CCC It gives us additional regularity of solution under the additional regularity of the loads. Proposition \ref{propen} is then proved by using this regularity and approximation argument. \BBB 
Before stating proving Proposition \ref{propreg}, we define the operator $\widetilde{\mathcal{A}}_0$ as a self-adjoint operator on $L^2(\Omega_p)$ defined through the bilinear form: 
\begin{equation} \label{krnj4} 
\widetilde{a}_0(v_1,v_2):=\int_{\Omega_p} \mathbb{K}_{33}(x_3) \partial_3 v_1(x)\cdot \partial_3 v_2(x) \,dx, \quad \forall v_1,v_2 \in V. 
\end{equation} 
Note that $\mathcal{D}(\widetilde{\mathcal{A}}_0)\subset  V$ and the operator $\mathcal{A}_0$ defined by \eqref{krnj3} coincides with the operator $\widetilde{\mathcal{A}}_0$ on $\mathcal{D}(\widetilde{\mathcal{A}}_0)$.
\begin{propozicija}\label{propreg}
If $\vect{F} \in H^2(0,T;L_{1,0}')$, $\vect{G} \in H^1(0,T;V')$  and $\widetilde{p}_0 \in \mathcal{D}(\widetilde{\mathcal{A}}_0)$, where $\widetilde{p}_0$ is defined by \eqref{deftp0} and $\widetilde{\mathcal{A}}_0$ is defined through bilinear form \eqref{krnj4}, then the solution of \eqref{diferformpen1*may111}- \eqref{diferformpen0*may111} satisfies $(\mathfrak{a},\mathfrak{b}, p) \in H^1(0,T;L_{1,0})  \times H^1(0,T;M)$. 
\end{propozicija}
\begin{proof} 
It is easy to see, by taking the derivative in time in \eqref{abspbm}, that if in the equation \eqref{abspbm} one takes $S$ such that $\frac{d}{dt} S \in L^2(0,T;V')$,  $u_0 \in V$ such that $\mathcal{A}_0 u_0 \in H$, then as a consequence we get that $\frac{d}{dt} u \in L^2(0,T;V)$. 
If the conditions of the proposition are satisfied, we are exactly in this situation. This can be seen by writing the system \eqref{diferformpen1*may111} and \eqref{diferformpen0*may111}  in the form of \eqref{abspbm} (cf. \eqref{defp}). 
 From this we conclude the regularity w.r.t. time for $p$. To conclude the regularity w.r.t. time for $\mathfrak{a}$ and $\mathfrak{b}$ we use \eqref{diferformpen1*may111}. 
\end{proof}  

\CCC Finally we prove the last result of this section. It is the energy-type equality valid for the solution of \eqref{diferformpen1*may111}-\eqref{diferformpen0*may111}.As usual,  such identities are needed to prove the strong convergence result.    \BBB
\begin{propozicija} \label{propen}
If 	$\vect{F} \in H^1(0,T;L_{1,0}')$  $\vect{G} \in L^2(0,T;V')$, $t_0 \in L^2(\Omega)$,  then the solution of \eqref{diferformpen1*may111}-\eqref{diferformpen0*may111} satisfies the following identity for every $t \in (0,T)$ :	
\begin{equation}\label{eneq2}
    \begin{split}
        &\frac{1}{2}\displaystyle\int_{\omega}\mathbb{A}^{\rm hom}(e_{\widehat{x}}(\mathfrak{a}(t)), \nabla^2_{\widehat{x}}\mathfrak{b}(t))):(e_{\widehat{x}}(\mathfrak{a}(t)), \nabla^2_{\widehat{x}}\mathfrak{b}(t))d\widehat{x}+\displaystyle\frac{1}{2}\int_{\Omega_p}M_0(x_3)p^2 (t)dx \\ &+\int_0^t\int_{\Omega_p} \mathbb{K}_{33} (x_3)|\partial_3 p|^2 dx \, dt  = \displaystyle_{L_{1,0}'}\langle  \vect{F}(t),(\mathfrak{a}(t),\mathfrak{b}(t)) \rangle_{L_{1,0}}
- \displaystyle_{L_{1,0}'}\langle  \vect{F}(0),(\mathfrak{a}(0),\mathfrak{b}(0))\rangle_{L_{1,0}} \\ 
 &-\int_0^t {_{L_{1,0}'}}\langle \partial_{t}\vect{F}(\tau),(\mathfrak{a}(\tau),\mathfrak{b}(\tau))\rangle_{L_{1,0}}\,d\tau
+\int_0^t{_{V'}}\langle \vect{G}, p\rangle_V\,dt + \frac{1}{2}\int_{\Omega_p} M_0(x_3) (\widetilde{p}_0)^2 \,dx,
    \end{split}
\end{equation}
where $\widetilde{p}_0$ is defined with \eqref{deftp0} (and $p_0$ with \eqref{krnj1}). 
	\end{propozicija}  
\begin{proof}
\CCC	The proof goes directly, if we have additional regularity of the solution.  
\BBB   
First we consider the system \eqref{diferformpen1*may111}-\eqref{diferformpen0*may111} and assume that 	$\vect{F} \in H^2(0,T;L_{1,0}')$,  $\widetilde{p}_0 \in \mathcal{D}(\widetilde{\CCC \mathcal{A}_9 \BBB})$. Using Proposition \ref{propreg} we conclude that the solution of \eqref{diferformpen1*may111} and \eqref{diferformpen0*may111} satisfies $(\mathfrak{a},\mathfrak{b}, p) \in H^1(0,T;L_{1,0})  \times H^1(0,T;V)$. 
By doing integration by parts and utilizing a density argument we conclude from \eqref{diferformpen0*may111} that $p(0)=\widetilde{p}_0$ and for a.e. $t \in (0,T)$ we have
\begin{equation}\label{diferformpen0pommay11} 
    \begin{split}
         &\displaystyle\int_{\Omega_p}M_0(x_3) \partial_t p\,  \varphi dx+\displaystyle\int_{\Omega_p}\mathbb{K}_{33}(x_3)\partial_{3}p\,\partial_{3}  \varphi dx 
+\displaystyle\int_{\omega}\int_{J_p}(|\mathcal{Y}_{f}(x_3)| \mathbb{I}-\mathbb{B}^{H}(x_3))\varphi \, dx_3:\iota\left(e_{\widehat{x}}\left(\mathfrak{\partial_t a}\right)\right) d\widehat{x}\\
&-\displaystyle\int_{\omega}\int_{J_p}\left(|\mathcal{Y}_f(x_3)|\mathbb{I}-\mathbb{B}^H(x_3)\right)x_3\varphi \, dx_3:\iota\left(\nabla^2_{\widehat{x}}\partial_t\mathfrak{b}\right)\,d\widehat{x}= {_{V'}} \langle \vect{G},p\rangle_V\,dt \quad  \forall  \varphi \in V.
    \end{split}
\end{equation}
We take $t \in [0,T]$.
By testing  \eqref{diferformpen1*may111} with $(\partial_{\tau} \mathfrak{a}, \partial_{\tau} \mathfrak{b})$, integrating over $[0,t]$ and adding \eqref{diferformpen0pommay11}  tested with 
$\varphi=p$ and then integrated over $(0,t)$
we obtain
\begin{equation} \label{eneq1}
    \begin{split}
        &\frac{1}{2}\displaystyle\int_{\omega}\mathbb{A}^{\rm hom}(e_{\widehat{x}}(\mathfrak{a}(t)), \nabla^2_{\widehat{x}}\mathfrak{b}(t)):(e_{\widehat{x}}(\mathfrak{a}(t)), \nabla^2_{\widehat{x}}\mathfrak{b}(t))d\widehat{x}dt+\displaystyle\frac{1}{2}\int_{\Omega_p}M_0(x_3)p^2(t) dx dt\\ &+\int_{0}^{t}\int_{\Omega_p} \mathbb{K}_{33} |\partial_3 p|^2 dx \,d\tau =
 \displaystyle\int_0^t{_{L_{1,0}'}} \langle  \vect{F},(\partial_{t} \mathfrak{a},\partial_{t}\mathfrak{b})\rangle_{L_{1,0}}\,d\tau+\int_0^T {_{V'}} \langle \vect{G},p\rangle_V\,dt+\frac{1}{2}\int_{\Omega_p} M_0(x_3) (\widetilde{p}_0)^2 \,dx.
    \end{split}
\end{equation}
We obtain \eqref{eneq2} for these more regular loads (and initial condition) by doing integration by parts in the first term on the right hand side of \eqref{eneq1}. The statement of the proposition then follows by approximation of the loads and initial condition from the stability results \eqref{stab1} and \eqref{stab11} and the stability estimate for the equation \eqref{diferformpen1*may111}, using the fact that $H^2(0,T;L_{1,0}')$ is dense in $H^1(0,T;L_{1,0}')$ and that $\mathcal{D}(\widetilde{\mathcal{A}}_0)$ is dense in $L^2(\Omega)$.
\end{proof} 	

\begin{remarkica} \label{obsshane} 
By using the operator $\widetilde{\mathcal{A}}_0$	and under the assumption that $\mathcal{B}_0$ is coercive, we can use the semigroup theory to obtain the existence for the problem \eqref{abspbm}, cf. Section \ref{anallim2} below. This, however, gives only the existence in $C([0,T];H)$ and in order to conclude the existence in $L^2(0,T;V)$ as well as the stability estimate \eqref{stab1} we would have to do approximation of the solution by the Galerkin method as is done in Section \ref{anallim2}. Note, however, that putting \eqref{abspbm} in the framework \eqref{abstractevolutionproblem} would be simpler then in the case of non-zero inertial term, since the corresponding operator $\mathcal{A}$ would then be positive definite (as usual, we first need to do the transformation of the solution $\vect{u}(t)\mapsto e^{-\lambda t} u(t)$ which converts the condition \eqref{krnj20} in the positive definiteness of the corresponding operator). By using semigroup approach the claim of Proposition \ref{propreg} becomes clear: if the initial condition is regular and the loads are regular the solution of the problem \eqref{abstractevolutionproblem} is, in fact, strong.  

In Section \ref{anallim2} for the case of non-zero inertial term we used both approaches: semigroup approach and Galerkin approximation.
The semigroup approach doesn't give us {\textit a priori} enough regularity to prove the existence of the solution, but it gives us a good way to approximate the initial condition and loads to obtain regular enough solution. This is used to prove energy-type equality. The existence of the solution is proved via Galerkin approximation. 
\end{remarkica} 	

  \subsection{Strong convergence}\label{secstrongcon}  

 The strong convergence is proved under the condition $\vect{F}(0)=0$.
 This condition was not needed to derive the limit equations  and establish weak convergence (\CCC as it was assumed in e.g. \cite{Mikelic2012}\BBB), but for  establishing the strong convergence
 this seems to be reasonable condition, since under that condition we have that $\mathfrak{a}(0)=0$, $\mathfrak{b}(0)=0$, i.e. the initial conditions are preserved.
 We will firstly prove the convergence of energies.
 \begin{teorem} \label{proplenka1}
  Let assumptions \eqref{kirill3} and \eqref{forcesassumptions1} be satisfied. 
  We also suppose that Assumptions \ref{assumption on regions} is satisfied and $\pi_h\vect{F}^{h}\xrightarrow{t,2-r,2} \vect{F},$ $\pi_h\partial_t \vect{F}^{h}\xrightarrow{t,2-r,2} \partial_t \vect{F}$, 
 	where $\vect{F} \in H^1(0,T; L^2(\Omega \times \mathcal{Y};\mathbb{R}^3))$ and 
 	$\vect{F}(0)=0$.  We have the following convergences:
{\allowdisplaybreaks
\begin{align}
\displaystyle\lim_{h\to 0} \frac{1}{h^{2}}\int_0^T\int_{\Omega_{s}^{h}}\mathbb{A}^h(x) e_h(\vect{u}^h(t)):e_h(\vect{u}^h(t))\,dx&\nonumber=\displaystyle\int_0^T\int_{\omega}\mathbb{A}^{\rm hom}(e_{\widehat{x}}(\mathfrak{a}(t)), \nabla^2_{\widehat{x}}\mathfrak{b}(t)):(e_{\widehat{x}}(\mathfrak{a}(t)), \nabla^2_{\widehat{x}}\mathfrak{b}(t))d\widehat{x}\,dt\\* &\label{strong2may16}+\displaystyle \int_0^T\int_{\Omega}M_0(x_3)p^2(t) dx\,dt,\\
\displaystyle\lim_{h\to 0}\frac{\ee^2}{h^4}\int_{0}^{T}\int_{\Omega^{h}_f}|e_h\left( \partial_{t}\vect{u}^h\right)|^2dx\,dt&\label{strong1may16}=\int_{0}^{T}\int_{\Omega\times\mathcal{Y}_f^{x_3}}\left|e_y\left(\partial_{t} \vect{u}^f_0\right)\right|^2dx\,dy\,dt,\\
 \label{lastvel} 
\lim_{h \to 0} \eta \int_0^T \left\|\partial_t\vect{u}^h(t)\right\|^2_{ L^2(\Omega;\RR^3)}\,dt &=0.
\end{align}
}
Here $(\mathfrak{a},\mathfrak{b}, p) \in L^2(0,T;L_{\kappa,0})  \times L^2 (0,T;M)$ are the solutions of \eqref{diferformpen1*may11}-\eqref{diferformpen0*may11} and the function $\vect{u}^f_0 \in  H^1(0,T;L^2(\Omega; H^1(\mathcal{Y}_f(x_3);\RR^3)))$ is defined with \eqref{diferformfluid2march302}. Furthermore,  $(\vect{u}^h,p^h)$ is the solution of \eqref{varformporous2gamma0may10} with initial condition \eqref{eq:51rescaledver}. 
\end{teorem}
 
 \begin{proof}
We will prove the convergences \eqref{strong2may16} and \eqref{lastvel} and instead of \eqref{strong1may16} we will firstly prove
\begin{equation}\label{strong1may16pom}
 \displaystyle\lim_{h\to 0}\frac{\ee^2}{h^4}\int_0^T\int_{0}^{t}\int_{\Omega^{h}_f}|e_h\left( \partial_{t}\vect{u}^h\right)|^2dx\,d\tau \,dt=\int_0^T\int_{0}^{t}\int_{\Omega\times\mathcal{Y}_f^x}\left|e_y\left(\partial_{t} \vect{u}^f_0\right)\right|^2dx\,dy\,d\tau\,dt. 
\end{equation}
By using Remark \ref{obsjosip1} we  take $\vect{v}=\partial_{t}\vect{u}^h$ as test function in (\ref{varformporous2gamma0may1000}). This yields for every $t \in [0,T]$
 \begin{align}
 	&\displaystyle\frac{1}{2}\frac{d}{dt}\left(\displaystyle\int_{\Omega}\eta \kappa^{h}\nonumber|\partial_{t}\vect{u}^h(t)|^2\,dx+\frac{1}{h^{2}}\int_{\Omega_{s}^{h}}\mathbb{A}^h(x) e_h(\vect{u}^h(t)):e_h(\vect{u}^h(t))\,dx\right)\\
 	&\quad+\displaystyle\frac{\ee^{2}}{h^{4}}\int_{\Omega^{h}_{f}}2 |e_h(\partial_{t}\vect{u}^h(t))|^2\,dx=\int_{\Omega}\vect{F}^h\partial_{t}\vect{u}^h\,dx.
 \end{align}
 Fixing $t \in [0,T]$,   integrating over interval $[0,t]$ and then over interval $[0,T]$ we obtain
 \begin{align}
 	&\displaystyle\frac{1}{2}\displaystyle\int_0^T\int_{\Omega}\eta \kappa^{h}\nonumber|\partial_{t}\vect{u}^h(t)|^2\,dx\, dt+\frac{1}{2 h^{2}}\int_0^T\int_{\Omega_{s}^{h}}\mathbb{A}^h(x) e_h(\vect{u}^h(t)):e_h(\vect{u}^h(t))\,dx\, dt \\
 	&\label{eq:62may16}+\displaystyle\frac{\ee^{2}}{h^{4}}\int_0^T\int_{0}^{t}\int_{\Omega^h_{f}} |e_h(\partial_{t}\vect{u}^h(\tau))|^2\,dx\,d\tau\, dt=\int_0^T\int_{0}^{t}\int_{\Omega}\vect{F}^h(\tau)\partial_{t}\vect{u}^h(\tau)\,dx\,d\tau\,dt.
 \end{align}
 Notice that
 \begin{align}     & \int_0^T\int_{0}^{t}\int_{\Omega}\vect{F}^h(\tau)\partial_{t}\vect{u}^h(\tau)\,dx\,d\tau\,dt= \int_0^T\int_{\Omega} \vect{F}^h(t)\vect{u}^h(t)\,dx\,dt-\int_0^T\int_{0}^{t}\int_{\Omega}\partial_{t}\vect{F}^h(\tau)\vect{u}^h(\tau)\,dx\,d\tau\,dt 
 	\\  &\to\displaystyle\int_0^T\int_{\omega} \overline{\langle \vect{F}\rangle_{Y}}\cdot(\mathfrak{a},\mathfrak{b})\, d\widehat{x}\,dt  +\displaystyle\int_0^T\int_0^t\int_{\omega} \overline{\langle x_3\partial_{t}\vect{F}_{*}(\tau)\rangle_{Y}}\cdot\nabla_{\widehat{x}}\mathfrak{b}(\tau)\, d\widehat{x}\,d\tau\,dt   \\  
 	& -\displaystyle\int_0^T\int_{\omega} \overline{\langle x_3\vect{F}_{*}\rangle_{Y}}\cdot\nabla_{\widehat{x}}\mathfrak{b}\, d\widehat{x}\,dt-\displaystyle\int_0^T\int_0^t\int_{\omega} \overline{\langle\partial_{t}\vect{F}(\tau)\rangle_{Y}}\cdot(\mathfrak{a}(\tau),\mathfrak{b}(\tau))\, d\widehat{x}\,d\tau\,dt. \label{preen}
 \end{align}
 Notice also that 
 \begin{equation}\label{ter2} 
\int_{\Omega\times\mathcal{Y}_f^(x_3)}\left|e_y\left(\partial_{t}\vect{u}^0_f(x,y)\right)\right|^2dx\,dy=\int_{\Omega}K_{33}(x)|\partial_3 p(x,\tau)|^2 dx,
 \end{equation}
 where we used \eqref{diferformfluid2march301} and \eqref{diferformfluid2march302} and the fact that 
 \begin{equation} \label{idk}  \mathbb{K}_{33}(x_3)=\int_{\mathcal{Y}_f(x_3)}q^3_{3,x_3}(y)\,dy=\langle \nabla \vect{q}_{x_3}^3,\nabla \vect{q}_{x_3}^3\rangle_{L^2(\mathcal{Y}_f(x_3);\mathbb{R}^3)} =\int_{\mathcal{Y}_f(x_3)} e_y(\vect{q}_{x_3}^3):e_y(\vect{q}_{x_3}^3)\, dy. 
 \end{equation} 
 \CCC To conclude the last equality we used divergence free condition of $\vect{q}^3_{x_3}$. \BBB
 
 Using \eqref{eneq2} integrated over $[0,T]$ (with $\mathfrak{a}(0)=0$, $\mathfrak{b}(0)=0$, $\widetilde{p}_0=0$), from \eqref{preen} and \eqref{idk} we conclude that 
 \begin{align} 
 	&\displaystyle\frac{1}{2}\displaystyle\int_0^T\int_{\Omega}\eta \kappa^{h}\nonumber|\partial_{t}\vect{u}^h(t)|^2\,dx\, dt+\frac{1}{2 h^{2}}\int_0^T\int_{\Omega_{s}^{h}}\mathbb{A}^h(x) e_h(\vect{u}^h(t)):e_h(\vect{u}^h(t))\,dx\, dt \\
 	&\nonumber+\displaystyle\frac{\ee^{2}}{h^{4}}\int_0^T\int_{0}^{t}\int_{\Omega^h_{f}} |e_h(\partial_{t}\vect{u}^h(\tau))|^2\,dx\,d\tau\, dt \to \frac{1}{2}\displaystyle\int_0^T\int_{\omega}\mathbb{A}^{\rm hom}(e_{\widehat{x}}(\mathfrak{a}), \nabla^2_{\widehat{x}}\mathfrak{b}):(e_{\widehat{x}}(\mathfrak{a}), \nabla^2_{\widehat{x}}\mathfrak{b})d\widehat{x}dt\\  &\label{pred1}+\displaystyle\frac{1}{2}\int_0^T\int_{\Omega}M_0(x_3)p^2 dx\,dt  +\int_{0}^{T}\int_0^t\int_{\Omega} \mathbb{K}_{33} |\partial_3 p|^2 dx \,d\tau\, dt.
 \end{align} \CCC
 Note that as a consequence of \eqref{gassmantensorapril01}, \eqref{cellproblem1april1}, \eqref{eq:abril1},  \eqref{defM0} and \eqref{eq:abril123} 
 we have that 
 \begin{align} \nonumber
 & \int_0^T \int_{\Omega\times \mathcal{Y}_s^{x_3}} \mathbb{A} (x_3,y)\left( 	\iota (e_{\widehat{x}}(\mathfrak{a})-x_3\nabla^2_{\widehat{x}}\mathfrak{b})+\mathfrak{C}_{\infty}(\vect{w}, \vect{g})\right):\left(	\iota (e_{\widehat{x}}(\mathfrak{a})-x_3\nabla^2_{\widehat{x}}\mathfrak{b})+\mathfrak{C}_{\infty}(\vect{w}, \vect{g})\right)\,dx\,dy\, dt\\
 &= \label{ter1}  \displaystyle\int_0^T\int_{\omega}\mathbb{A}^{\rm hom}(e_{\widehat{x}}(\mathfrak{a}), \nabla^2_{\widehat{x}}\mathfrak{b}):(e_{\widehat{x}}(\mathfrak{a}), \nabla^2_{\widehat{x}}\mathfrak{b})d\widehat{x}\,dt+ \int_0^T\int_{\Omega}M_0(x_3)p^2 dx\,dt.
 \end{align} 	
 \BBB
 By using lower semicontinuity of the convex energy we have the following inequalities (using \eqref{depra}, \eqref{cr2gamma0}, \eqref{cr5gamma00}, \eqref{ter2} and \eqref{ter1}): 
  \begin{equation} \label{odozddruga}
  \begin{split}
      & \liminf_{h\to 0} \frac{1}{ h^{2}}\int_0^T\int_{\Omega_{s}^{h}}\mathbb{A}^h(x) e_h(\vect{u}^h(t)):e_h(\vect{u}^h(t))\,dx\, dt \geq  \\ &   \displaystyle\int_0^T\int_{\omega}\mathbb{A}^{\rm hom}(e_{\widehat{x}}(\mathfrak{a}), \nabla^2_{\widehat{x}}\mathfrak{b}):(e_{\widehat{x}}(\mathfrak{a}), \nabla^2_{\widehat{x}}\mathfrak{b})d\widehat{x}dt+ \int_0^T\int_{\Omega}M_0(x_3)p^2 dx\,dt,
  \end{split}
 \end{equation}
  \begin{equation} 
 	\label{odozdtreca} 
 	\liminf_{h\to 0} \displaystyle\frac{\ee^{2}}{h^{4}}\int_0^T\int_{0}^{t}\int_{\Omega^h_{f}} |e_h(\partial_{t          }\vect{u}^h(\tau))|^2\,dx\,d\tau\, dt \geq \int_0^T\int_0^t\int_{\Omega}K_{33}(x)|\partial_3 p(x,\tau)^2dx\,d\tau\,dt,
 \end{equation}
  \begin{equation} 
 	\label{odozdprva} \liminf_{h\to 0} \frac{1}{2}\displaystyle\int_0^T\int_{\Omega}\eta \kappa^{h}|\partial_{t}\vect{u}^h(t)|^2\,dx\, dt\geq \min\{\kappa^h\}\eta 
 	\liminf_{h\to 0} \frac{1}{2}\displaystyle\int_0^T\int_{\Omega} |\partial_{t}\vect{u}^h(t)|^2\,dx\, dt \geq 0.
 \end{equation}
 By combining \eqref{odozddruga}-\eqref{odozdprva} with \eqref{pred1} we obtain \eqref{strong2may16}, \eqref{lastvel} and \eqref{strong1may16pom}. It remains to prove \eqref{strong1may16}.
  We take a sequence $(\vect{u}_0^{f,n})_{n \in \N} \subset H^1(0,T;C^1(\Omega;H^1(\mathcal{Y}_f(x_3)
 ;\mathbb{R}^3)))$, such that
 $$ \varepsilon \chi_{\Omega_f^h}(x) e_x \left(\partial_t \vect{u}_0^{f,n}(x,\frac{x}{\varepsilon})\right) \xrightarrow{t,2-r,2} e_y(\partial_t\vect{u}_f^0).$$
 Notice that \eqref{strong1may16pom} is equivalent with (taking into account \eqref{cr5gamma00}):
 \begin{equation} \label{andrey} 
 	\lim_{n \to \infty}\lim_{h \to 0} \int_0^T \int_0^t\int_{\Omega_f^h} \left|\frac{\varepsilon}{h^2} e_h(\partial_{t} \vect{u}^h)- \varepsilon e_x\left(\partial_{t} \vect{u}_0^{f,n}(x,\frac{x}{\varepsilon})\right) \right|^2dx\,d \tau\,dt=0. 
 \end{equation} 
 \eqref{andrey} implies that 
 $$ \lim_{h\to 0}\frac{\ee^2}{h^4}\int_{0}^{T'}\int_{\Omega^{h}_f}|e_h\left( \partial_{t}\vect{u}^h\right)|^2dx\,dt=\int_{0}^{T'}\int_{\Omega\times\mathcal{Y}_f(x_3)}\left|e_y\left(\partial_{t} \vect{u}^f_0\right)\right|^2dx\,dy\,dt, \quad \forall T'<T.$$
 
 In order to conclude it also for $T'=T$ one can extend the forces on the interval $[0,T'']$, where $T''>T$ with constant forces $\vect{F}^h(T)$ and then use
 the uniqueness of the solution and the above result. 
 \end{proof}

We have the following strong convergence result which includes correctors. 

 \begin{teorem}\label{twoscaleconvergencesmay16}
Under the conditions of Theorem \ref{proplenka1} the following convergences hold:
\begin{equation}\label{strongtwoscale1may16}
\displaystyle\lim_{n \to \infty}\lim_{h\to 0}\int_{0}^{T}\int_{\Omega^{h}_f}\left|\frac{\ee}{h^2}e_h(\partial_\tau \vect{u}^h)(x,\tau)+e_y (\vect{q}^3_n)\left(x,\frac{\widehat{x}}{\ee},\frac{x_3}{\frac{\ee}{h}}\right)\partial_3 p(x,\tau)\right|^2\,dx\,d\tau=0,
 \end{equation}
{\allowdisplaybreaks
\begin{align}
    &\nonumber\displaystyle\lim_{n \to \infty} \lim_{h\to 0}\int_0^T\int_{\Omega^{h}_{s} } \left|\frac{1}{h}e_{h}\left(\vect{u}^h\right)-e_{h}\left(
    \begin{pmatrix}
          \mathfrak{a}_1(\widehat{x},t)-x_3\partial_1\mathfrak{b}(\widehat{x},t) \\
           \mathfrak{a}_2(\widehat{x},t)-x_3\partial_2\mathfrak{b}(\widehat{x},t) \\
           h^{-1}\mathfrak{b}(\widehat{x},t)
    \end{pmatrix} \right)\right.\\
    &\nonumber\quad-e_{h}\left(\ee \displaystyle\sum_{\alpha, \beta}\left(e_{\widehat{x}}(\mathfrak{a}_n(\widehat{x},t))\right)_{\alpha\beta}\vect{w}^{\alpha\beta}_n\left(x,\frac{\widehat{x}}{\ee},\frac{x_3}{\frac{\ee}{h}}\right) \right)\\
    &\nonumber\quad+e_{h}\left(\ee x_3\displaystyle\sum_{\alpha, \beta}\partial_{\alpha\beta}\mathfrak{b}_n(\widehat{x},t)\vect{w}^{\alpha\beta}_n\left(x,\frac{\widehat{x}}{\ee},\frac{x_3}{\frac{\ee}{h}}\right) \right)-e_{h}\left(\ee p_n(x,t)\widetilde{\vect{w}}_{2,n}\left(x,\frac{\widehat{x}}{\ee},\frac{x_3}{\frac{\ee}{h}}\right) \right)\\
    &\quad-e_{h}\left(h \displaystyle\sum_{\alpha, \beta}\left(e_{\widehat{x}}(\mathfrak{a}_n(\widehat{x},t))\right)_{\alpha\beta}\int_{0}^{x_3}\vect{g}^{\alpha\beta}_n(x)\,dx_3 \right) \label{strongtwoscale4may16}\\
    & \nonumber \left.\quad+e_{h}\left(h x_3\displaystyle\sum_{\alpha, \beta}\partial_{\alpha\beta}\mathfrak{b}_n(\widehat{x},t)\int_{0}^{x_3}\vect{g}^{\alpha\beta}_n(x)\,dx_3 \right)-e_{h}\left(h p_n(x,t)\int_{0}^{x_3}\widetilde{\vect{g}}_{2,n}(x)\,dx_3 \right)\right|^2dx\,dt=0.
 \end{align}
 }
Here $\vect{q}^3_n \in L^{\infty}(\Omega; C^1(\mathcal{Y}_f(x_3);\mathbb{R}^3)) $, $
 \mathfrak{a}_n(\cdot,t)\in L^{\infty}(0,T;C^2_{\#}(\omega;\mathbb{R}^2))$, $\mathfrak{b}_n \in L^{\infty} (0,T;C_{\#}^3(\omega))$, for $\alpha,\beta=1,2$, $\vect{w}^{\alpha\beta}_n \in C^1(\Omega;C^1(\mathcal{Y}_s(x_3);\mathbb{R}^3))$, $\widetilde{\vect{w}}_{2,n} \in C^1(\Omega;C^1(\mathcal{Y}_s(x_3);\mathbb{R}^3))$, $\vect{g}^{\alpha\beta}_n \in C^1(\Omega;\mathbb{R}^3)$, $\widetilde{\vect{g}}_{2,n} \in C^1(\Omega; \mathbb{R}^3)$ are chosen such that:
 \begin{eqnarray*} 
 	& & \|\vect{q}_n^3-\vect{q}^3\|_{L^2(\Omega;H^1(\mathcal{Y}_f(x_3);\mathbb{R}^3))} \to 0, \quad \|\vect{q}_n^3\|_{L^\infty(\Omega;H^1(\mathcal{Y}_f(x_3);\mathbb{R}^3))} \textrm{ bounded}; \\ & &
 	\|\mathfrak{a}_n-\mathfrak{a}\|_{L^2(0,T;H_{\#}^1(\omega;\mathbb{R}^2))} \to 0, \quad \min\{\frac{\eps}{h},h\}\| \mathfrak{a}_n\|_{L^{\infty}(0,T;C^2_{\#}(\omega;\mathbb{R}^2))} \to 0; 	
 	\\ & &
 	\|\mathfrak{b}_n-\mathfrak{b}\|_{L^2(0,T;H_{\#}^2(\omega))} \to 0, \quad \min\{\frac{\eps}{h},h\}\| \mathfrak{b}_n\|_{L^{\infty}(0,T;C^3_{\#}(\omega))} \to 0; \\ & &
 	\|\vect{w}^{\alpha\beta}_n- \vect{w}^{\alpha\beta}\|_{L^2(\Omega;H^1(\mathcal{Y}_s(x_3);\mathbb{R}^3))}\to 0, \quad \|\vect{w}^{\alpha\beta}_n\|_{L^{\infty}(\Omega;H^1(\mathcal{Y}_s(x_3);\mathbb{R}^3))} \textrm{ bounded; } 
 	\\ & & \|p_n-p\|_{L^2(0,T;M)}\to 0, \quad \min\{\frac{\varepsilon}{h},h\} \|p_n\|_{L^{\infty}(0,T;C^1(\omega))}\to 0; \\ & &
 	\|\widetilde{\vect{w}}_{2,n} -\widetilde{\vect{w}}_{2}\|_{L^2(\Omega;H^1(\mathcal{Y}_s(x_3);\mathbb{R}^3))}\to 0, \quad \|\widetilde{\vect{w}}_{2,n}\|_{L^\infty(\Omega;H^1(\mathcal{Y}_s(x_3);\mathbb{R}^3))} \textrm{ bounded }; 
 	\\ & & \| \vect{g}_n^{\alpha\beta}-\vect{g}^{\alpha\beta}\|_{L^2(\Omega;\mathbb{R}^3)} \to 0, \quad  \| \vect{g}_n^{\alpha\beta}\|_{L^\infty(\Omega;\mathbb{R}^3)} \textrm{ bounded}; \\ & &
 	\| \widetilde{\vect{g}}_{2,n}-\widetilde{\vect{g}}_{2}\|_{L^2(\Omega;\mathbb{R}^3)} \to 0, \quad \| \widetilde{\vect{g}}_{2,n}\|_{L^\infty(\Omega;\mathbb{R}^3)} \textrm{ bounded}.
 \end{eqnarray*} 	  
\CCC  Here for $x=(\hat{x}, x_3)$ we put $\vect{q}^3(x,\cdot):=\vect{q}^3_{x_3}(\cdot)$, $\vect{w}^{\alpha\beta}(x,\cdot):=\vect{w}^{\alpha\beta}_{x_3}(\cdot)$, $\widetilde{\vect{w}}_2(x,\cdot):=\widetilde{\vect{w}}_{2,x_3}(\cdot)$, $\widetilde{\vect{g}}_{2}(x):=\widetilde{\vect{g}}_{2,x_3}$ (recall \eqref{eq unique uationapril16}, \eqref{cellproblem1april1}, \eqref{diferformfluid2march301}) . \BBB
 \end{teorem} 

 \begin{proof}
  Using the convexity of the quadratic energy, \eqref{strongtwoscale1may16} is equivalent with
 \eqref{strong1may16} as a consequence of \eqref{diferformfluid2march302}. 
 In the similar way \eqref{strongtwoscale4may16} is equivalent with \eqref{strong2may16}, using  \eqref{gassmantensorapril01}, \eqref{limitequationhom1} and \eqref{eq:abril123}. \CCC To compare, see the proof of \cite[Proposition 6.1, Theorem 6.1]{Mikelic2012}. 
 Additional regularization through the approximation sequences is done in order to make the convergences valid. 
 \BBB

 \end{proof}
 
 \subsection{Further discussions on the limit model} 
 \label{secgen} 
 In this section we discuss the possible generalization of the Assumptions \ref{assumption on regions},  the possibility of different initial conditions and boundary conditions at the transverse boundary, possibility of imposing surface loads as well as the situation when condition \eqref{nodownup} is not satisfied. Everything discussed here can also be concluded 
for the situation with inertial term analyzed in Section \ref{secinertia}.
\CCC 
\subsubsection{Different initial and boundary conditions at transverse boundary} 
It is not a problem to assume different initial conditions (the initial position and initial velocity different than zero). However, in that case one needs to assume certain natural compatibility conditions, see \cite[Section 2.2, Section 3.1]{DUGun}.

Also, non-periodic boundary conditions at the transverse boundary  can be imposed (like e.g., Dirichlet boundary conditions, see \cite{Murat} for discussion of different possibilities for boundary conditions for Stokes equation). All the estimates obtained in this paper can be easily adapted to Dirichlet boundary conditions (see, e.g., \cite{mbukalvelcic2017,Marinvelciczubrinic2022}).  
However, one would need to do further discussion like analyze separately $C^{1,1}$ domain when using results for corrector from Appendix \ref{griso}. Periodic boundary conditions simplify this and are also assumed in \cite{clopeaumikelic2001,Mikelic2012}.

\BBB
\subsubsection{Generalization of Assumptions \ref{assumption on regions}}
By analyzing the proof of Theorem \ref{gammazeromay11} and Theorem \ref{effectiveequationsmay11} (see also Remark \ref{remgotovo}) we can see that it was not necessary to demand that the sets $U_i$ are of the form $U_i=\omega \times I_i$. Instead, they can be assumed to be generally Lipschitz. \CCC It is not difficult to adapt the property of the existence or non-existence of flow, which has to be localized for every point of the interface. 
\BBB
 This would imply that the sets $J_p$ and $J_K$ defined in \eqref{defjp} and \eqref{defjk} depend also on $\hat{x}$. We also emphasize the fact that the contact of elastic and poroelastic plate as well as poroelastic and poroelastic plate that have interface on the plane parallel to $x_3$ axis doesn't produce any problems for the analysis (the reason being that only the third derivative of the pressure appears in the limit model).

The other possible generalization might go in the direction of having infinite number of possible cell-types. This is plausible, provided that we can construct the extension operator with controllable constant (see Lemma \ref{lext} and Lemma \ref{lemma 1}). As a consequence it is possible to obtain the limit tensors which are not piecewise constants. It is also conceivable that with suitable interface conditions we would obtain the similar limit equations as we obtained here. 

 \subsubsection{Surface  loads}
From the application point of view it is important to have surface loads on the upper and lower part of the boundary, $\omega \times\{-\frac{1}{2},\frac{1}{2}\}$. The discussion here resembles \cite[Corollary 3.42]{Marinvelciczubrinic2022},  where the case of elastic high-contrast inclusions is analyzed. 

We assume that in \eqref{varformporous2gamma0may10}, on the right hand side,  we have the surface loads: 
$$\mathcal{G}^{h}(\vect G^{h})(\vect v)= \int_{\omega \times\{-1/2,1/2\}} \vect G^{h} \vect v\, d\widehat{x},\qquad \vect v \in H^1_{\#}(\Omega;\R^3),$$ 
satisfying
\begin{align}
	&\pi_h \vect{G}^{h} \subset  H^1(0,T,L^2\bigl(\omega\times \bigl\{-1/2, 1/2\bigr\};\R^3\bigr)) \textrm{ is bounded},   \\  & \pi_h \vect{G}^h \xrightharpoonup{t,2-r,2} \vect G  \in H^1\Bigl(0,T;L^2\bigl(\omega\times\bigl\{-1/2, 1/2\bigr\};\R^3\bigr)\Bigr).
\end{align} 
Instead of the condition \eqref{kirill3}, we demand that for every $t>0$ we have
$$\int_{\Omega} \vect{F}^h(t) \,dx+ \int_{\omega\times \{-1/2,1/2\}} \vect{G}^h(t) \, d\widehat{x}=0. $$
In the limit equations we then obtain the following terms on the right hand side of \eqref{diferformpen1*may11} 
$$ \int_{\omega} \vect{G}\cdot \vect{\theta} \,d\widehat{x}+\int_{\omega} \left(\vect{G}(\cdot,-1/2)-\vect{G}(\cdot,1/2)\right)\cdot \nabla_{\widehat{x}} \theta_3 \, d\widehat{x}. $$
Under the assumption on strong convergence,  the strong convergence of the solutions can also be proved in the spirit of Section \ref{secstrongcon}. Note also that the existence and uniqueness from Section \ref{analysis} can be used. 
\subsubsection{Fluid part intersects the upper and lower boundary} \label{secintersection} 

In cases where condition \eqref{nodownup} is not met, it becomes necessary to enforce boundary conditions for the fluid component at  the upper and lower boundary. While incorporating surface loads on the elastic component presents no issue, specifying boundary conditions for the fluid poses a challenge.

In \cite{Murat}, various boundary conditions and their corresponding equations in the weak formulation for the Stokes problem are discussed. It's important to note that non-homogeneous Dirichlet boundary conditions cannot be directly imposed on the upper and lower boundary for the fluid. This limitation arises due to the necessity of employing a special scaling that we used to derive the limiting equations.

Under this scaling, it follows that $\vect{u}_f$  is of the order of $h^2$. Here, $\vect{u}_f$ doesn't represent the actual position of the fluid (which is given by $\vect{u}=\widehat{\vect{u}}+\vect{u}^f$), and the scaling of the fluid's position consequently differs (being of order one in the third component). Consequently, the imposition of non-homogeneous Dirichlet boundary conditions on the upper and lower boundary is not possible.

As an alternative, one can only apply artificial Dirichlet boundary conditions on these boundaries, taking the form:
$$ \vect{u}_f= h^2 \vect{g} \textrm{ on } \partial \omega \times \{-1/2,1/2\} \cap \Omega^h_f, \quad  
\vect{g} \in H^1(0,T;V^h) \cap H^2(0,T;L^2(\Omega;\RR^3)).$$
As a consequence of the derivation \eqref{conv2may11} we would obtain the additional term of the form $$\int_{\omega} \left( g_3(\cdot,1/2)\cdot \varphi(\cdot,1/2)- g_3(\cdot,-1/2)\cdot \varphi(\cdot,-1/2) \right)\, d\widehat{x} $$
on the right hand side of \eqref{diferformpen0*may11}. 
Section \ref{analysis} also includes the existence and uniqueness result for this case. However, due to the real Dirichlet condition for velocity being $\vect{u}=\widehat{\vect{u}}+\vect{u}_f$, where $\widehat{\vect u}$ denotes the extension from the solid part, and considering that prescribing deformation for the solid part at the upper and lower boundaries isn't feasible (thus surface loads are applied), this scenario necessitates a thorough analysis concerning the solution's existence and\emph{ a priori }estimates for the microscopic solution.

Additionally, the other boundary conditions discussed in \cite{Murat} involve aspects of the fluid velocity (i.e. position) and pressure, encountering a similar issue as with the Dirichlet boundary condition. Nonetheless, if these conditions are imposed on $\vect{u}_f$, it becomes feasible to address the subsequent condition:
\begin{equation}\label{bcc}  
 (u_f)_{\alpha}=h^2 g_{\alpha}, \quad p=p_l, \textrm{ on }  \partial \omega \times \{-1/2,1/2\} \cap \Omega^h_f, 
\end{equation}
 The impact of condition \eqref{bcc} on the weak formulation is addressed in \cite{Murat}. This involves the introduction of additional terms on the right-hand side of \eqref{varformporous2gamma0may10}, consequently leading to the emergence of further terms on the right-hand side of \eqref{diferformpen1*may11}. Once more, this necessitates preliminary analysis regarding the solution's existence and\emph{ a priori }estimates for the microscopic problem.

 Generally, it is difficult to justify from the microscopic point of view, any other condition besides the condition $\partial_3 p=0$ at the upper and lower boundary. 
 \section{The problem with inertial term}
\label{secinertia} 
We consider equation \eqref{varformporous2gamma0may10} with $\eta=\eta(h)=1$. Following the structure laid out in Section \ref{secquasi-static}, we maintain a similar approach but omit the corresponding proofs.

In Section \ref{compactness2}, we present the compactness result and formulate the limit problem. Subsequently, in Section \ref{anallim2}, we establish both the existence and uniqueness of the solution, along with demonstrating an energy-type equality for the limit problem. 
\CCC In \cite[Appendix B]{gurvich} it is already noted that the semigroup approach can be helpful in obtaining the existence result. However, here again an additional effort needs to be done to define the appropriate operators, since the limit equations do not decouple, see the proof of Proposition \ref{propsemigroup} below). In the end, we will obtain the existence result by using Galerkin approximation, but we will use semigroup approach to obtain the appropriate approximation of     the solution, which we will use to prove energy-type equality in the similar way as it was done in quasi-static case. 
\BBB
Finally, in Section \ref{strong2}, we articulate the result regarding strong convergence.
 The following \CCC set \BBB will be of importance: \CCC For given $\mathfrak{b} \in H^1_{\#}(\omega)$ we denote 
$$L_{\kappa,\mathfrak{b}}:=\left\{\vect{\psi} \in H^1_{\#}(\omega;\mathbb{R}^2) : \int_{\omega} \bar{\kappa}\vect{\psi}=-\int_{\Omega} \kappa \nabla_{\hat{x}} \mathfrak{b} \, d\hat{x} \right\}.$$
Note that for any $\mathfrak{b} \in H^1_{\#}(\omega) $ we have $L_{1}\equiv L_{1,\mathfrak{b}}=\dot{H}^1_{\#}(\omega;\R^2)
$. \BBB
	\BBB
Furthermore, we introduce  $N:=H^2_{\#}(\omega)$ and recall $M:=L^2(\omega; H^1(J_K)\oplus L^2(J_p\backslash J_K))$.

\subsection{Compactness and limit equations}\label{compactness2} 
The following theorem can be proved in the same way as Theorem \ref{gammazeromay11} and Theorem \ref{effectiveequationsmay11}. 
\CCC It provides compactness statement and establishes the limit problem. \BBB 
We will state it without proof. 
\begin{teorem}\label{compinertia}  
Let assumptions \eqref{kirill3} and \eqref{forcesassumptions0} be satisfied. 
We also suppose that Assumptions \ref{assumption on regions} is satisfied and $\pi_h\vect{F}^{h}\xrightharpoonup{t,2-r,2} (\vect{F}_*,F_3) \footnote{Again as a consequence of \eqref{kirill3} we have $\int_{\omega} \overline{\langle \vect{F}\rangle_{Y}}=0.$},$
where $\vect{F}_* \in H^1(0,T; L^2(\Omega \times \mathcal{Y};\mathbb{R}^2))$, $F_3 \in L^2(0,T; L^2(\Omega \times \mathcal{Y}))$. 
The following statements hold: Let $(\vect{u}^h,p^h)$ be the solution of (\ref{varformporous2gamma0may10}) with initial condition \eqref{eq:51rescaledver}. Then there exist limits $\mathfrak{a} \in L^{\infty}(0,T;L_{\kappa,\mathfrak{b}}), \mathfrak{b}\in L^{\infty}(0,T;N)\CCC \cap H^1(0,T;L^2(\omega))\BBB,\vect{w}, \vect{g}$ such that they satisfy \eqref{rev1111}-\eqref{kirill10} and $p \in L^{\infty}(0,T; L^2(\Omega_p)) \cap L^2(0,T;M)$ such that the convergences \eqref{crgamma0}-\eqref{cr5gamma0} are satisfied and
\begin{equation} \label{kreso2} 
\partial_t \vect{u}^h \xrightharpoonup{L^\infty(0,T; L^2(\Omega))} (0,0, \partial_t \mathfrak{b}),\quad \kappa^h\partial_t \vect{u}^h \xrightharpoonup{L^\infty(0,T; L^2(\Omega))} \bar{\kappa} (0,0, \partial_t \mathfrak{b}).  
\end{equation} 	
Furthermore, the following limit equations are satisfied\footnote{Again in the second term in the second equation we can have as a test function $\varphi \in L^2(0,T;M)$ since $\mathbb{K}_{33}=0$ on $\Omega_p \backslash \Omega_K$.  }: 
\begin{equation}
\label{kreso3}  \mathfrak{b}(0)=0,
\end{equation}
\begin{equation}\label{kreso4}
\begin{split}
    	&-\int_0^T\int_{\omega} \bar{\kappa}(\widehat{x}) \partial_t \mathfrak{b}\partial_t \theta_3 \,dt+ \displaystyle\int_0^T\int_{\omega}\mathbb{A}^{\rm hom}(e_{\widehat{x}}(\mathfrak{a}), \nabla^2_{\widehat{x}}\mathfrak{b}):(e_{\widehat{x}}(\vect{\theta}_*), \nabla^2_{\widehat{x}}\theta_3) d\widehat{x}\, dt\\ &-\displaystyle\int_0^T\int_{\omega}\int_{J_p}\left(|\mathcal{Y}_f(x_3)| \mathbb{I}-\mathbb{B}^{H}(x_3)\right)p \, dx_3:\left[\iota (e_{\widehat{x}}(\vect{\theta}_{*}))\right]\,d\widehat{x}\,dt\\
	&+\displaystyle\int_0^T\int_{\omega} \int_{J_p}\left(|\mathcal{Y}_f(x_3)| \mathbb{I}-\mathbb{B}^{H}(x_3)\right)x_3 p \, dx_3:\left[\iota (\nabla^2_{\widehat{x}}\theta_3)\right]\,d\widehat{x}\,dt =\displaystyle\int_0^T\int_{\omega} \overline{\langle \vect{F}\rangle_{Y}}\cdot(\vect{\theta}_*,\theta_3)\, d\widehat{x}\, dt\\
	 &-\displaystyle\int_0^T\int_{\omega} \overline{\langle x_3\vect{F}_{*}\rangle_{Y}}\cdot\nabla_{\widehat{x}}\theta_3  \, d\widehat{x}\, dt,\quad\forall \left(\vect{\theta}_{*},{\theta}_3\right) \in L^2(0,T;L_{1})\times H^1(0,T;N),  \textrm{ s.t. } {\theta}_3(T)=0,
\end{split}
\end{equation}
\begin{equation}\label{kreso5}
\begin{split}
    &\displaystyle-\int_0^T\int_{\Omega_p}M_0(x_3) p\, \partial_t \varphi dx\,dt+\displaystyle\int_0^T\int_{\Omega_p}\mathbb{K}_{33}(x_3)\partial_{3}p\,\partial_{3}  \varphi dx\,dt\\ & -\displaystyle\int_0^T\int_{\omega}\int_{J_p}(|\mathcal{Y}_{f}(x_3)| \mathbb{I}-\mathbb{B}^{H}(x_3))\partial_t\varphi \, dx_3:\iota\left(e_{\widehat{x}}\left(\mathfrak{a}\right)\right) d\widehat{x}\,dt\\&+\displaystyle\int_0^T\int_{\omega}\int_{J_p}\left(|\mathcal{Y}_f(x_3)|\mathbb{I}-\mathbb{B}^H(x_3)\right)x_3\partial_t\varphi \, dx_3:\iota\left(\nabla^2_{\widehat{x}}\mathfrak{b}\right)\,d\widehat{x}\,dt=0,\\ &  \forall  \varphi \in L^2(0,T;M) \cap H^1(0,T;L^2(\Omega_p)) \textrm{ such that } \varphi(T)=0. 
\end{split}
\end{equation}
\end{teorem} 
\CCC 
\begin{remarkica} \label{interpret1} 
In the limit one obtains partially quasi-static model (there is no inertia term for the part of in-plane component $\mathfrak{a}$).  This also happens in the case of elastic plate (see \cite{ciarlet2}). The non-decoupling of bending (inertial) and membrane (quasi-static) equations also happens in elasticity when there are heterogeneities across thickness, see \cite{Marinvelciczubrinic2022}. As already mentioned, although the 3D evolution equations derived in \cite{clopeaumikelic2001} have memory effects, the  2D  equations obtained here do not have memory effects. This is the consequence of (implicit) scaling of time that is done to obtain the limit model (this scaling of time is also done in linearized elasticity to obtain the plate models, see \cite{ciarlet2,Marinvelciczubrinic2022}). However, in \cite{Marinvelciczubrinic2022} also the models obtained without scaling of time (they are called membrane models) are discussed.  
\end{remarkica} 	
\BBB
 \subsection{Analysis of the limit equations}\label{anallim2}
We again equip the system \eqref{kreso3}-\eqref{kreso5} with non-zero initial conditions and substitute the loads \CCC and replace set $L_{\kappa,\mathfrak{b}}$ with $L_{1}$.  \BBB
For $\mathfrak{b}_0 \in N$, $\mathfrak{b}_1 \in L^2(\omega)$, $t_0 \in L^2(\Omega_p)$,  $\vect{F}_{1} \in H^1(0,T; L'_{1})$, $\vect{F}_2 \in L^2(0,T;N')$, $\vect{G} \in L^2(0,T;M')$,  we define the following problem: Find $\mathfrak{a} \in L^2(0,T;L_1)$, $\mathfrak{b} \in L^2(0,T;N)\cap H^1(0,T;L^2(\omega))$, $p \in L^2(0,T; M)$ such that
\begin{equation}
   \label{kreso33}  \mathfrak{b}(0)=\mathfrak{b}_0, 
\end{equation}
\begin{equation}\begin{split}
        	&-\int_0^T \int_{\omega} \bar{\kappa}(\widehat{x})\partial_t \mathfrak{b}\partial_t \theta_3 \,dt+ \displaystyle\int_0^T\int_{\omega}\mathbb{A}^{\rm hom}(e_{\widehat{x}}(\mathfrak{a}), \nabla^2_{\widehat{x}}\mathfrak{b}):(e_{\widehat{x}}(\vect{\theta}_*), \nabla^2_{\widehat{x}}\theta_3) d\widehat{x}\, dt\\ &-\displaystyle\int_0^T\int_{\omega}\int_{J_p}\left(|\mathcal{Y}_f(x_3)| \mathbb{I}-\mathbb{B}^{H}(x_3)\right)p \, dx_3:\left[\iota (e_{\widehat{x}}(\vect{\theta}_{*}))\right]\,d\widehat{x}\,dt\\
	&+\displaystyle\int_0^T\int_{\omega} \int_{J_p}\left(|\mathcal{Y}_f(x_3)| \mathbb{I}-\mathbb{B}^{H}(x_3)\right)x_3 p \, dx_3:\left[\iota (\nabla^2_{\widehat{x}}\theta_3)\right]\,d\widehat{x}\,dt =\displaystyle\int_0^T{_{L'_{1}}} \langle \vect{F}_1,\vect{\theta}_*\rangle_{L_1}\, dt\, \label{kreso44} \\
	&+\int_0^T{_{N'}} \langle \vect{F}_2,\theta_3\rangle_{N} dt+\int_{\omega} \mathfrak{b}_1 \cdot \theta_3(0) \,d \widehat{x},\quad\forall \left(\vect{\theta}_{*},{\theta}_3\right) \in L^2(0,T;L_{1})\times H^1(0,T;N),  \textrm{ s.t. } {\theta}_3(T)=0, \\
\end{split}\end{equation}
\begin{equation}\label{kreso55} 
    \begin{split}
        &\displaystyle-\int_0^T\int_{\Omega_p}M_0(x_3) p\, \partial_t \varphi dx\,dt+\displaystyle\int_0^T\int_{\Omega_p}\mathbb{K}_{33}(x_3)\partial_{3}p\,\partial_{3}  \varphi dx\,dt\\ & -\displaystyle\int_0^T\int_{\omega}\int_{J_p}(|\mathcal{Y}_{f}(x_3)| \mathbb{I}-\mathbb{B}^{H}(x_3))\partial_t\varphi \, dx_3:\iota\left(e_{\widehat{x}}\left(\mathfrak{a}\right)\right) d\widehat{x}\,dt\\
&+\displaystyle\int_0^T\int_{\omega}\int_{J_p}\left(|\mathcal{Y}_f(x_3)|\mathbb{I}-\mathbb{B}^H(x_3)\right)x_3\partial_t\varphi \, dx_3:\iota\left(\nabla^2_{\widehat{x}}\mathfrak{b}\right)\,d\widehat{x}\,dt=\int_0^T{_{M'}} \langle \vect{G},\varphi\rangle_{M}\,dt\\  &+ \int_{\Omega_p} t_0 \cdot \varphi(0) ,\quad  \forall  \varphi \in L^2(0,T;M) \cap H^1(0,T;L^2(\Omega_p)) \textrm{ such that } \varphi(T)=0. 
    \end{split}
\end{equation}
\CCC Again, we can assume that $\mathfrak{a} \in L_{1}$, since the solution in $\mathfrak{a} \in L_{\kappa}$ can be obtained by adding a constant for every $t \in [0,T]$. Note that, however, we keep $\overline{\kappa}$ in the first term on the left hand side of \eqref{kreso44}. \BBB 

This section  discusses existence, uniqueness and regularity of the solution of \eqref{kreso33}-\eqref{kreso55} as well as the energy-type equality.
The following theorem gives us uniqueness of solution  of \eqref{kreso33}-\eqref{kreso55}. 
\begin{teorem} \label{uniqueinertia} 
	Let $\mathfrak{b}_0 \in N$, $\mathfrak{b}_1 \in L^2(\omega)$, $t_0 \in L^2(\Omega_p)$,  $\vect{F}_{1} \in H^1(0,T; L'_{1})$, $\vect{F}_2 \in L^2(0,T;N')$, $\vect{G} \in L^2(0,T;M')$. The solution of \eqref{kreso33}-\eqref{kreso55}, if it exists, is unique. 
\end{teorem} 	
\begin{proof} 
	We assume that $\vect{F}_{\alpha}=\vect{G}=0$, for $\alpha=1,2$, $\mathfrak{b}_0=\mathfrak{b}_1=0$, $t_0=0$ in \eqref{kreso33}-\eqref{kreso55} and for fixed $0<T'<T$  take the following test functions in \eqref{kreso44} 
	$$ \theta_3=\begin{cases} \int_{T'}^t \mathfrak{b} (s) \,ds, & \textrm{ for all } 0\leq t<T', \\ 0, & \textrm{ for all } T'\leq t \leq T.   \end{cases}  \quad \vect{\theta}_*=\begin{cases} \int_{T'}^t \mathfrak{a}(s)\,ds, & \textrm{ for all } 0 \leq t <T' \\ 0, & \textrm{ for all } T'\leq t \leq T.  \end{cases}, $$
	and in \eqref{kreso55}  
	$$   \varphi(t)=\begin{cases} \int_{T'}^{t} \int_0^s p(r)\,dr\,ds,& \textrm{ for all }  0\leq t <T', \\ 0,& \textrm{ for all } T'\leq t \leq T.  \end{cases}  $$
	We have that 
\begin{flalign*}
    & -\int_0^T \int_{\omega} \bar{\kappa}(\widehat{x}) \partial_t \mathfrak{b} \partial_t \theta_3 \,dx\,dt =- \frac{1}{2} \int_{\omega} \bar{\kappa}(\widehat{x})\int_0^{T'}\frac{d}{dt} \mathfrak{b}^2 \,dt\,d\widehat{x}=-\frac{1}{2}\int_{\omega}\bar{\kappa}(\widehat{x})\mathfrak{b}^2(T')\, d\widehat{x}, &&
\end{flalign*}\vspace{-0.75\baselineskip}
\begin{flalign*}
    & \int_0^{T}\int_{\omega}\mathbb{A}^{\rm hom}(e_{\widehat{x}}(\mathfrak{a}), \nabla^2_{\widehat{x}}\mathfrak{b}):(e_{\widehat{x}}(\vect{\theta}_*), \nabla^2_{\widehat{x}}\theta_3)\, d\widehat{x}\, dt &&
    \\& = \frac{1}{2} \int_0^{T'} \frac{d}{dt} \int_{\omega} \mathbb{A}^{\rm hom}\left(e_{\widehat{x}}(\int_{T'}^t\mathfrak{a}(s)\,ds), \nabla^2_{\widehat{x}}(\int_{T'}^t\mathfrak{b}(s) \,ds)\right):\left(e_{\widehat{x}}(\int_{T'}^t\mathfrak{a}(s)\,ds), \nabla^2_{\widehat{x}}(\int_{T'}^t\mathfrak{b}(s) \,ds)\right) d\widehat{x}\,dt &&
    \\& = -\frac{1}{2}  \int_{\omega} \mathbb{A}^{\rm hom}\left(e_{\widehat{x}}(\int_0^{T'}\mathfrak{a}(s)\,ds), \nabla^2_{\widehat{x}}(\int_0^{T'}\mathfrak{b}(s) \,ds)\right):\left(e_{\widehat{x}}(\int_0^{T'}\mathfrak{a}(s)\,ds), \nabla^2_{\widehat{x}}(\int_0^{T'}\mathfrak{b}(s) \,ds)\right) d\widehat{x}, &&
\end{flalign*}\vspace{-0.75\baselineskip}
\begin{flalign*}
    \int_0^{T}\int_{\omega}\int_{J_p}(|\mathcal{Y}_f(x_3)| & \mathbb{I}-\mathbb{B}^{H}(x_3))p \, dx_3:\left[\iota (e_{\widehat{x}}(\vect{\theta}_{*}))\right]\,d\widehat{x}\,dt &&
    \\& = \int_0^{T'}\int_{\omega}\int_{J_p}\left(|\mathcal{Y}_f(x_3)| \mathbb{I}-\mathbb{B}^{H}(x_3)\right)\frac{d}{dt} \int_0^tp(s) \,ds\, dx_3:\left[\iota (e_{\widehat{x}}(\int_{T'}^t \mathfrak{a}(s) \,ds))\right]\,d\widehat{x}\,dt &&
    \\& = -\int_0^{T'}\int_{\omega}\int_{J_p}\left(|\mathcal{Y}_f(x_3)| \mathbb{I}-\mathbb{B}^{H}(x_3)\right) \int_0^tp \, dx_3:\left[\iota (e_{\widehat{x}}( \mathfrak{a}))\right]\,d\widehat{x}\,dt, &&
\end{flalign*}\vspace{-\baselineskip}
\begin{flalign*}     
    \int_0^T\int_{\omega} \int_{J_p}\left(|\mathcal{Y}_f(x_3)| \mathbb{I}-\mathbb{B}^{H}(x_3)\right)&x_3 p \, dx_3:\left[\iota (\nabla^2_{\widehat{x}}\theta_3)\right]\,d\widehat{x}\,dt &&
    \\& = - \int_0^{T'}\int_{\omega}\int_{J_p}\left(|\mathcal{Y}_f(x_3)| \mathbb{I}-\mathbb{B}^{H}(x_3)\right) \int_0^tp \, dx_3:\left[\iota (\nabla^2_{\widehat{x}}\mathfrak{b})\right]\,d\widehat{x}\,dt, &&
\end{flalign*}\vspace{-\baselineskip}
\begin{flalign*}
    \int_0^T\int_{\Omega_p}M_0(x_3) p\, \partial_t \varphi dx\,dt &= \int_0^{T'} \int_{\Omega} M_0(x_3) p\int_0^t p dx \,dt &&
    \\&  = \frac{1}{2}\int_{0}^{T'} \frac{d}{dt} \int_{\Omega_p}M_0(x_3) \left( \int_0^t p\right)^2dx \,dt=\frac{1}{2}\int_{\Omega}M_0(x_3)  \left( \int_0^{T'} p\right)^2\,dx, &&
\end{flalign*}\vspace{-\baselineskip}
\begin{flalign*}
    \int_0^T\int_{\Omega_p}\mathbb{K}_{33}(x_3)\partial_{3}p\,\partial_{3}  \varphi dx\,dt&= \int_0^{T'} \int_{\Omega} \mathbb{K}_{33}(x_3) \partial_t  \int_0^t \partial_3p\, \int_{T'}^t \int_0^s \partial_3  p(r)\,dr \,ds \,dx \,dt &&
    \\& =-\int_0^{T'} \int_{\Omega} \mathbb{K}_{33}(x_3) \partial_3 \int_0^t p\, \partial_3 \int_0^t p \,dx\,dt, &&
\end{flalign*}\vspace{-\baselineskip}
\begin{flalign*}
	\int_0^T\int_{\omega}\int_{J_p}(|\mathcal{Y}_{f}(x_3)| \mathbb{I}-\mathbb{B}^{H}(x_3))&\partial_t\varphi \, dx_3:\iota\left(e_{\widehat{x}}\left(\mathfrak{a}\right)\right) d\widehat{x}\,dt &&
    \\& = \int_0^{T'}\int_{\omega}\int_{J_p}\left(|\mathcal{Y}_f(x_3)| \mathbb{I}-\mathbb{B}^{H}(x_3)\right) \int_0^tp \, dx_3:\left[\iota (e_{\widehat{x}}( \mathfrak{a}))\right]\,d\widehat{x}\,dt, &&
\end{flalign*}\vspace{-\baselineskip}
\begin{flalign*}
	\int_0^T\int_{\omega}\int_{J_p}\left(|\mathcal{Y}_f(x_3)|\mathbb{I}-\mathbb{B}^H(x_3)\right)&x_3\partial_t\varphi \, dx_3:\iota\left(\nabla^2_{\widehat{x}}\mathfrak{b}\right)\,d\widehat{x}\,dt &&
    \\& = \int_0^{T'}\int_{\omega}\int_{J_p}\left(|\mathcal{Y}_f(x_3)| \mathbb{I}-\mathbb{B}^{H}(x_3)\right) \int_0^tp \, dx_3:\left[\iota (\nabla^2_{\widehat{x}}\mathfrak{b})\right]\,d\widehat{x}\,dt. &&
\end{flalign*}
	Multiplying both \eqref{kreso44} and \eqref{kreso55} with $-1$ and adding them we obtain 
	\begin{eqnarray*} 
		& &	\frac{1}{2} \int_{\omega} \bar{\kappa}(\widehat{x})\mathfrak{b}^2(T')\, d \widehat{x}\\ & &+\frac{1}{2}  \int_{\omega} \mathbb{A}^{\rm hom}\left(e_{\widehat{x}}(\int_0^{T'}\mathfrak{a}(s)\,ds), \nabla^2_{\widehat{x}}(\int_0^{T'}\mathfrak{b}(s) \,ds)\right):\left(e_{\widehat{x}}(\int_0^{T'}\mathfrak{a}(s)\,ds), \nabla^2_{\widehat{x}}(\int_0^{T'}\mathfrak{b}(s) \,ds)\right) d\widehat{x} \\ & &
		+\frac{1}{2}\int_{\Omega}M_0(x_3)  \left( \int_0^{T'} p\right)^2\,dx+\int_0^{T'} \int_{\Omega} \mathbb{K}_{33}(x_3) \partial_3 \int_0^t p\, \partial_3 \int_0^t p \,dx\,dt=0.
	\end{eqnarray*} 	
	Since this is valid for arbitrary $0\leq T'\leq T$ we obtain $\mathfrak{b}=p=0$ on $[0,T]$. To conclude uniqueness of $\mathfrak{a}$ we use \eqref{kreso44} with $\theta_3=0$ and $\vect{\theta}_*$ arbitrary and the fact that $\mathfrak{a} \in L_1$. This proves uniqueness. 
\end{proof}

Before proving the existence result, we first recall some results from operator theory. Let  $(H, \langle \cdot,\cdot \rangle_H)$ be a separable  Hilbert space.  Let $\mathcal{A}:H \to H$ be a generator of the strongly continuous  semigroup $e^{t\mathcal{A}}$, $\|e^{t\mathcal{A}}\|_{H \to H}\leq M_{\mathcal{A}} e^{\omega_{\mathcal{A}} t}$, $t\geq 0$, where $M_{\mathcal{A}}>0$ is a positive constant and $\omega_\mathcal{A} \in \R$.  Here $\|\cdot \|_{H\to H}$ denotes the operator norm for the bounded operators from $H$ to $H$.
For $\vect{f} \in L^1(0,T;H)$, consider the following evolution problem:
\begin{equation}
	\label{abstractevolutionproblem}
     \begin{cases}
    \partial_{t} \vect u(t) =\mathcal{A} \vect u(t) + \vect f(t), &\\
    \vect u(0) = \vect u_0, \quad
    \vect u_0 \in H.  &
    \end{cases}
\end{equation}
\begin{definition}\label{sol1} 
We call $\vect{u}\in C(0,T;H)$ a mild solution of \eqref{abstractevolutionproblem} if 
	\begin{equation} \label{ivan903}
		{\vect u}(t) =  e^{t\mathcal{A}}{\vect u}_0 + \int_0^T  e^{(t-s)\mathcal{A}}{ \vect f}(s)ds. 
	\end{equation}
\end{definition} 		   	   
\begin{remarkica}\label{pedro1}
Under the condition that $\vect{u}_0 \in \mathcal{D}(\mathcal{A})$ and $\vect{f} \in L^p(0,T;\mathcal{D}(\mathcal{A}))$ we have that the mild solution is the strong solution of \eqref{abstractevolutionproblem} in $W^{1,p}(0,T;H)\cap C (0,T;\mathcal{D}(\mathcal{A}))$ (recall that since $\mathcal{A}$ is closed operator, $\mathcal{D}(\mathcal{A})$ is a Hilbert space with graph norm).  	
\end{remarkica} 		   
	   \CCC The following lemma gives us a stability estimate.\BBB
\begin{lema} 	   
We have that the mild solution of \eqref{abstractevolutionproblem} satisfies:
	\begin{equation} 
		\| \vect u\|_{L^{\infty}(0,T;H)} \leq M_{\mathcal{A}} e^{\omega_{\mathcal{A}} T}\left(\|{\vect u}_0\|_H+\|\vect f\|_{L^1(0,T;H)}\right).
	\end{equation}
\end{lema} 
\begin{proof} 
The claims follows directly from \eqref{ivan903}.
\end{proof}

Before proving the existence result we will show that the system \eqref{kreso33}-\eqref{kreso55} can be put formally in the form of \eqref{abstractevolutionproblem} (with $\mathcal{A}$ generator of the contraction semigroup on appropriate Hilbert space $H$). Although from this we will not be able to conclude the existence of the mild solution for \eqref{kreso33}-\eqref{kreso55} (since the loads are not regular enough) it will give us some information on the system which we will use to prove energy-type equality. Also this  will enable us to obtain the existence of the solution for each step of Galerkin approximation.

\begin{propozicija} \label{propsemigroup}
There exists a Hilbert space $(H,\langle\cdot\rangle_H)$ and a Hilbert space $(V,\langle \cdot\rangle_V)$ densely embedded in $H$ such that the system \eqref{kreso33}-\eqref{kreso55} can be put in the form \eqref{abstractevolutionproblem} for $\vect{f} \in L^2(0,T; V')$ and where $\mathcal{A}$ is a generator of a strongly continuous semigroup on $H$.  
\end{propozicija} 	
\begin{proof} 
{\bf Step 1}.{\it  Definition of $H$.} 	 We define 
	\begin{equation} \label{defH11} H:=N \times L^2(\omega) \times L^2(\Omega_p),\end{equation} 
	endowed with the scalar product
	$$\langle (\mathfrak{b}_1,\mathfrak{v}_1, p_1)^T,(\mathfrak{b}_2,\mathfrak{v}_2, p_2)^T  \rangle_{\tilde b,H} =\fint_{\omega}\mathfrak{b}_1\cdot \fint_{\omega}\mathfrak{b}_2+a_1(\mathfrak{b}_1,\mathfrak{b}_2) +\langle \bar{\kappa} \mathfrak{v}_1,\mathfrak{v}_2 \rangle_{L^2(\omega)}+\widetilde{b}(p_1,p_2),  $$
	where $a_1$ is a bilinear form on $H^2_{\#}(\omega)$ is  defined by: 
	\begin{equation}  \label{utoest1} 
		a_1(\mathfrak{b}_1, \mathfrak{b}_2)=\int_{\omega}\mathbb{A}^{\rm hom}(e_{\widehat{x}}(\mathfrak{a}_1^{\mathfrak{b}_1}), \nabla^2_{\widehat{x}}\mathfrak{b}_1):(e_{\widehat{x}}(\mathfrak{a}_1^{\mathfrak{b}_2}), \nabla^2_{\widehat{x}}\mathfrak{b}_2) d\widehat{x}.
	\end{equation} 	
	Here $\mathfrak{a}_1^{\widehat{\mathfrak{b}}} \in L_1$ for $\widehat{\mathfrak{b}} \in N$ is a unique  solution of  
	\begin{equation}\label{poneq1} 
		\int_{\omega}\mathbb{A}^{\rm hom}(e_{\widehat{x}}(\mathfrak{a}_1^{\widehat{\mathfrak{b}}}), \nabla^2_{\widehat{x}}\widehat{\mathfrak{b}}):(e_{\widehat{x}}(\vect{\theta}_*), 0) d\widehat{x}=0, \quad \forall \vect{\theta}_* \in L_1.
	\end{equation} 	
	 $\widetilde{b}$ is a bilinear form on $L^2(\Omega_p)$ defined by 
	$
		\widetilde{b}(p_1,p_2)=\langle \widetilde{\mathcal{B}}p_1,p_2\rangle_{L^2(\Omega_p)},	
	$
	where $\widetilde{\mathcal{B}}:L^2(\Omega) \to L^2(\Omega)$ is a positive definite bounded operator defined by bilinear form
	\begin{equation} \label{utoest3}  
		\langle \widetilde{\mathcal{B}}p_1,p_2 \rangle_{L^2(\Omega_p)} =\int_{\Omega_p} M_0(x_3) p_1p_2\, dx+ \int_{\omega}\int_{J_p}(|\mathcal{Y}_{f}(x_3)| \mathbb{I}-\mathbb{B}^{H}(x_3))p_2 \, dx_3:\iota\left(e_{\widehat{x}}\left(\mathfrak{a}_2^{p_1}\right)\right) d\widehat{x},
	\end{equation} 	
	and $\mathfrak{a}_2^{\widehat{p}} \in L_1$ for $\widehat{p} \in L^2(\Omega)$ is the unique solution of 
 \begin{equation}\label{utoest2}\begin{split}
      &\int_{\omega}\mathbb{A}^{\rm hom}(e_{\widehat{x}}(\mathfrak{a}_2^{\widehat{p}}), 0):(e_{\widehat{x}}(\vect{\theta}_*), 0) d\widehat{x}= \int_{\omega}\int_{J_p}\left(|\mathcal{Y}_f(x_3)| \mathbb{I}-\mathbb{B}^{H}(x_3)\right)\widehat{p}\, dx_3:\left[\iota (e_{\widehat{x}}(\vect{\theta}_{*}))\right]\,d\widehat{x}, \\
      &\quad \forall \vect{\theta}_* \in \dot{H}^1_{\#}(\omega;\mathbb{R}^2).
 \end{split}\end{equation} 
	From \eqref{utoest3} and \eqref{utoest2} it follows that
	\begin{equation} 
		\langle \widetilde{\mathcal{B}}p_1,p_2 \rangle_{L^2(\Omega_p)} =\int_{\Omega_p} M_0(x_3) p_1 p_2\,dx+ \int_{\omega}\mathbb{A}^{\rm hom}(e_{\widehat{x}}(\mathfrak{a}_2^{p_1}), 0):(e_{\widehat{x}}(\mathfrak{a}_2^{p_2}), 0) d\widehat{x},
	\end{equation} 	
	from which we have positive definiteness. We denote the corresponding norm by $\|\cdot\|_{\tilde b,H}$.  
	We also endow $H$ with the scalar product
	which is given by
	$$\langle (\mathfrak{b}_1,\mathfrak{v}_1, p_1)^T,(\mathfrak{b}_2,\mathfrak{v}_2, p_2)^T  \rangle_{H} =\fint_{\omega}\mathfrak{b}_1\cdot \fint_{\omega}\mathfrak{b}_2+a_1(\mathfrak{b}_1,\mathfrak{b}_2) +\langle \mathfrak{v}_1,\mathfrak{v}_2 \rangle_{L^2(\omega)}+\langle p_1,p_2\rangle_{L^2(\Omega_p)}.  $$
	The corresponding norm, denoted by $\|\cdot\|_H$, is equivalent to the norm $\|\cdot\|_{\tilde b,H}$. Note that if we denote by $\mathcal{B}:H\to H$ the self adjoint continuous operator operator given by 
	\begin{equation} \label{defb11} 
	 \mathcal{B} (\mathfrak{b},\mathfrak{v},p)^T=(\mathfrak{b},\bar{\kappa} \mathfrak{v},\widetilde{\mathcal{B}}p)^T, \quad  \forall (\mathfrak{b},\mathfrak{v},p)^T \in H,
	 \end{equation} 
	then we have 
	\begin{equation} \label{krnic2} 
		\langle (\mathfrak{b}_1,\mathfrak{v}_1,p_1)^T, (\mathfrak{b}_2,\mathfrak{v}_2,p_2)^T \rangle_{\tilde b,H}= \langle \mathcal{B}  (\mathfrak{b}_1,\mathfrak{v}_1,p_1)^T, (\mathfrak{b}_2,\mathfrak{v}_2,p_2)^T\rangle_H, \quad \forall (\mathfrak{b}_{\alpha},\mathfrak{v}_{\alpha},p_{\alpha})^T \in H, \ \alpha=1,2,
	\end{equation} 	
	from which it follows that 
	$
		\|(\mathfrak{b},\mathfrak{v},p)^T\|_{\tilde b,H}= \|\mathcal{B}^{1/2}(\mathfrak{b},\mathfrak{v},p)^T\|_H. 
	$
 
	{\bf Step 2.} {\it Definition of $\mathcal{A}$}. We define formally the operator 
	$\widetilde{\mathcal{A}}$ on $H$ by
	\begin{equation} \label{bukal1} 
		\widetilde{\mathcal{A}}\left(   \mathfrak{b} , \mathfrak{v}, p  \right)^T=\left( \mathfrak{v}, \widetilde{\mathcal{A}_1}(\mathfrak{b}, p)^T , \widetilde{\mathcal{A}}_2 (\mathfrak{v},p)^T   \right)^T, 
	\end{equation}  
	where we formally put\footnote{At this point it is only formally since we did not specify the domain of $\widetilde{\mathcal{A}}$.}
  \begin{equation}\label{bukal2} 
     \begin{split}
         \langle \widetilde{\mathcal{A}}_1(\mathfrak{b},p)^T, \widetilde{\theta}_3\rangle_{L^2(\omega)}=&- \int_{\omega}\mathbb{A}^{\rm hom}(e_{\widehat{x}}(\mathfrak{a}_1^{\mathfrak{b}}+\mathfrak{a}_2^{p}), \nabla^2_{\widehat{x}}\mathfrak{b}):(0, \nabla^2_{\widehat{x}}\widetilde{\theta}_3)\, d\widehat{x} \\ 
		&  -\int_{\omega} \int_{J_p}\left(|\mathcal{Y}_f(x_3)| \mathbb{I}-\mathbb{B}^{H}(x_3)\right)x_3 p \, dx_3:\left[\iota (\nabla^2_{\widehat{x}}\widetilde{\theta}_3)\right]\,d\widehat{x}, \quad \forall \widetilde{\theta}_3 \in N. 
     \end{split}
 \end{equation}
	Here $\mathfrak{a}_1^{\mathfrak{b}} ,\mathfrak{a}_2^{\mathfrak{p}} \in L_1$ are defined by \eqref{poneq1} and \eqref{utoest2} respectively. 
	For the definition of $\widetilde{\mathcal{A}}_2$ we formally put 
 \begin{equation}\label{bukal3} 
     \begin{split}
         &\langle \widetilde{\mathcal{A}}_2 (\mathfrak{v},p)^T,  \varphi \rangle_{L^2(\Omega_p)}= -
		\int_{\omega}\int_{J_p}(|\mathcal{Y}_{f}(x_3)| \mathbb{I}-\mathbb{B}^{H}(x_3))\varphi \, dx_3:\iota\left(e_{\widehat{x}}\left(\mathfrak{a}_1^{\mathfrak{v}}\right)\right) d\widehat{x}\\ 
		& +\displaystyle  \int_{\omega}\int_{J_p}\left(|\mathcal{Y}_f(x_3)|\mathbb{I}-\mathbb{B}^H(x_3)\right)x_3\varphi \, dx_3:\iota\left(\nabla^2_{\widehat{x}}\mathfrak{v}\right)\,d\widehat{x}- \int_{\Omega_p}\mathbb{K}_{33}(x_3)\partial_{3}p\,\partial_{3}  \varphi \,dx,  \quad \forall \varphi \in M,  
     \end{split}
 \end{equation}
	where $\mathfrak{a}_1^{\mathfrak{v}} \in L_1$ is defined by \eqref{poneq1}. Next we want to show that (properly defined) modification of  operator $\widetilde{\mathcal{A}}$ is a generator of contraction semigroup on $H$. Firstly we want to show that the operator $I-\widetilde{\mathcal{A}}$ is bijective, i.e.
	for $\vect{x} \in H$ we want to show that there exists a unique $(\mathfrak{b},\mathfrak{v},p)^T \in H $ such that 
	\begin{equation} \label{utoest10} 
		(I-\widetilde{\mathcal{A}})\left(  \mathfrak{b}, \mathfrak{v}, p \right)^T=\vect{x}. 
	\end{equation}
	From the first equation in the system \eqref{utoest10} we have that 
	\begin{equation} \label{utoest11}
		\mathfrak{b}-\mathfrak{v} =\vect{x}_1 \Rightarrow \mathfrak{b}= \mathfrak{v}+\vect{x}_1. 
	\end{equation} 	
	Next we plug this into the second and third equation of \eqref{utoest10} and want to solve the second and third equation of \eqref{utoest10}.  
	We define the bilinear form on $N \times M $ in the following way: 
 \begin{equation}\label{tvuk2} 
     \begin{split}
         a((\mathfrak{v},p)^T, (\widetilde{\theta}_3,\varphi)^T)=&
		\int_{\omega}\mathbb{A}^{\rm hom}(e_{\widehat{x}}(\mathfrak{a}_1^{\mathfrak{v}}+\mathfrak{a}_2^{p}), \nabla^2_{\widehat{x}}\mathfrak{v}):(0, \nabla^2_{\widehat{x}}\widetilde{\theta}_3)\, d\widehat{x}\\  &+
		\int_{\omega} \int_{J_p}\left(|\mathcal{Y}_f(x_3)| \mathbb{I}-\mathbb{B}^{H}(x_3)\right)x_3 p \, dx_3:\left[\iota (\nabla^2_{\widehat{x}}\widetilde{\theta}_3)\right]\,d\widehat{x}\\  &+\int_{\omega}\int_{J_p}(|\mathcal{Y}_{f}(x_3)| \mathbb{I}-\mathbb{B}^{H}(x_3))\varphi \, dx_3:\iota\left(e_{\widehat{x}}\left(\mathfrak{a}_1^{\mathfrak{v}}\right)\right) d\widehat{x}\\  &-\int_{\omega}\int_{J_p}\left(|\mathcal{Y}_f(x_3)|\mathbb{I}-\mathbb{B}^H(x_3)\right)x_3\varphi \, dx_3:\iota\left(\nabla^2_{\widehat{x}}\mathfrak{v}\right)\,d\widehat{x}\\  &+\int_{\Omega_p}\mathbb{K}_{33}(x_3)\partial_{3}p\,\partial_{3}  \varphi \,dx. 
     \end{split}
 \end{equation}	
	Note that the second and third equation of \eqref{utoest10}, taking into account \eqref{utoest11},  can be interpreted as solving the equation: 
  \begin{equation}\label{utoest20} 
     \begin{split}
         \langle \mathfrak{v},\widetilde{\theta}_3\rangle_{L^2(\omega)}&+\langle p,\varphi \rangle_{L^2(\Omega_p)}+a((\mathfrak{v},p)^T, (\widetilde{\theta}_3,\varphi)^T)= \langle \vect{x}_2,\widetilde{\theta}_3 \rangle_{L^2(\omega)}+\langle \vect{x}_3,\varphi \rangle_{L^2(\Omega_p)} \\ &  
		  -\int_{\omega}\mathbb{A}^{\rm hom}(e_{\widehat{x}}(\mathfrak{a}_1^{\vect{x}_1}), \nabla^2_{\widehat{x}}\vect{x}_1):(0, \nabla^2_{\widehat{x}}\widetilde{\theta}_3) d\widehat{x}, \ \ \forall (\widetilde{\theta}_3,\varphi)^T \in N \times M.
     \end{split}
 \end{equation}
	As a consequence of the definition of  $\mathfrak{a}_1$, $\mathfrak{a}_2$ (see \eqref{poneq1}, \eqref{utoest2}) we have that 
  \begin{equation}\label{nedest4}
     \begin{split}
         	\int_{\omega}\mathbb{A}^{\rm hom}(e_{\widehat{x}}(\mathfrak{a}_1^{\mathfrak{v}}), \nabla^2_{\widehat{x}}\mathfrak{v}):(0, \nabla^2_{\widehat{x}}\widetilde{\theta}_3) d\widehat{x}=	\int_{\omega}\mathbb{A}^{\rm hom}(e_{\widehat{x}}(\mathfrak{a}_1^{\mathfrak{v}}), \nabla^2_{\widehat{x}}\mathfrak{v}):(e_{\widehat{x}}(\mathfrak{a}_1^{\widetilde{\theta}_3}), \nabla^2_{\widehat{x}}\widetilde{\theta}_3) d\widehat{x}, \,\,  \forall \mathfrak{v},\widetilde{\theta}_3 \in N. 
     \end{split}
 \end{equation}
	and also
	 \begin{equation}\label{subest2}
     \begin{split}
         &\int_{\omega}\mathbb{A}^{\rm hom}(e_{\widehat{x}}(\mathfrak{a}_2^{p}), 0):(0, \nabla^2_{\widehat{x}}\widetilde{\theta}_3) d\widehat{x}
		= -\int_{\omega}\mathbb{A}^{\rm hom}(e_{\widehat{x}}(\mathfrak{a}_2^{p}), 0):(e_{\widehat{x}}(\mathfrak{a}_1^{\widetilde{\theta}_3}), 0) d\widehat{x}\\ &=-\int_{\omega}\int_I(|\mathcal{Y}_{f}(x_3)| \mathbb{I}-\mathbb{B}^{H}(x_3))p \, dx_3:\iota(e_{\widehat{x}}(\mathfrak{a}_1^{\mathfrak{\widetilde{\theta}_3}})) d\widehat{x}, \quad  \forall p\in L^2(\Omega), \forall\widetilde{\theta}_3 \in N. 
     \end{split}
 \end{equation}
 \CCC
 Taking into account \eqref{nedest4} and \eqref{subest2} we have from \eqref{tvuk2} 
   \begin{equation}\label{nedest3}
 	\begin{split}
 	 a((\mathfrak{v},p)^T, (\widetilde{\theta}_3,\varphi)^T)	 =&  \int_{\omega}\mathbb{A}^{\rm hom}(e_{\widehat{x}}(\mathfrak{a}_1^{\mathfrak{v}}, \nabla^2_{\widehat{x}}\mathfrak{v}):(e_{\widehat{x}}(\mathfrak{a}_1^{\widetilde{\theta}_3}), \nabla^2_{\widehat{x}}\widetilde{\theta}_3) d\widehat{x}\\  & +
 		\int_{\omega}\int_{J_p}(|\mathcal{Y}_{f}(x_3)| \mathbb{I}-\mathbb{B}^{H}(x_3))\varphi \, dx_3:\iota\left(e_{\widehat{x}}\left(\mathfrak{a}_1^{\mathfrak{v}}\right)\right) d\widehat{x}\\
 		&-
 		\int_{\omega} \int_{J_p}\left(|\mathcal{Y}_f(x_3)| \mathbb{I}-\mathbb{B}^{H}(x_3)\right)x_3 \varphi \, dx_3:\left[\iota (\nabla^2_{\widehat{x}}\mathfrak{v})\right]\,d\widehat{x}\\  &-\int_{\omega}\int_{J_p}(|\mathcal{Y}_{f}(x_3)| \mathbb{I}-\mathbb{B}^{H}(x_3))p \, dx_3:\iota\left(e_{\widehat{x}}\left(\mathfrak{a}_1^{\widetilde{\theta}_3}\right)\right) d\widehat{x}\\  &+\int_{\omega}\int_{J_p}\left(|\mathcal{Y}_f(x_3)|\mathbb{I}-\mathbb{B}^H(x_3)\right)x_3p \, dx_3:\iota\left(\nabla^2_{\widehat{x}}\widetilde{\theta}_3\right)\,d\widehat{x}+\int_{\Omega_p}\mathbb{K}_{33}(x_3)\partial_{3}p\,\partial_{3}  \varphi \,dx,
 	\end{split}
 \end{equation}	
 \BBB 
From this  and Proposition \ref{tensorprop},  one easily obtains that the form on the left hand side of \eqref{utoest20} is coercive on $N \times M$. Since it is obviously continuous, it follows from Lax-Milgram \CCC(non-symmetric case) \BBB that \eqref{utoest20} has a solution for $\vect{x} \in H$. Taking into account \eqref{utoest11} we have that we can solve uniquely \eqref{utoest10} for $\vect{x} \in H$. Next we define the operator $\widetilde{\mathcal{A}}$ in the following way: 
	\begin{itemize} 
		\item the domain of the operator $\mathcal{D} (\widetilde{\mathcal{A}})$ is given as the subset of $H$ consisting of the solutions of \eqref{utoest10} as $\vect{x}$ exhausts $H$. Obviously 
		\begin{equation} \label{defV11}
			\mathcal{D}(\widetilde{\mathcal{A}}) \subset N \times N \times M=:V. 	
		\end{equation} 	
		\item for given $\vect{y} \in \mathcal{D}(\widetilde{\mathcal{A}})	$ we define 
		\begin{equation} \label{deftildea11}  \widetilde{\mathcal{A}}\vect{y}:=-\vect{x}(\vect{y})+\vect{y},
		\end{equation} 
		where $\vect{x}(\vect{y})$ is an element of $H$ that satisfies  
		$(I-\widetilde{\mathcal{A}})\vect{y}=\vect{x}(\vect{y}).  $
	\end{itemize} 	
	Note that $( I-\widetilde{\mathcal{A}})^{-1}:H \to H$ is a bounded operator. Furthermore, as a consequence of \eqref{utoest11} and \eqref{utoest20}, we have that \eqref{bukal2},\eqref{bukal3} are valid for every $(\mathfrak{b},\mathfrak{v},p)^T \in \mathcal{D}(\widetilde{\mathcal{A}})$, $\widetilde{\theta}_3 \in N$, $\varphi \in M$. 
 
 Next, we define the operator $\mathcal{C}$,   $H \ni \vect{y}\mapsto \mathcal{C} \vect{y}= (\mathfrak{b}_1,\widetilde{\theta}_3,\varphi)^T \in H$, where\footnote{Below we will show that $\mathcal{C}=\left((I-\widetilde{\mathcal{A}})^{-1}\right)^T$.} 
	\begin{equation} \label{nedest1}  \mathfrak{b}_1=\vect{y}_1-\widetilde{\theta}_3
		+\fint_{\omega} \widetilde{\theta}_3\, d\hat{x}, 
	\end{equation} 
	and $(\widetilde{\theta}_3,\varphi)^T \in N \times M$ are solutions of:
\begin{equation}\label{subest1} 
     \begin{split}
         & \langle \widetilde{\theta}_3,\mathfrak{v}\rangle_{L^2(\omega)}+\langle\varphi, p \rangle_{L^2(\Omega)}+\widetilde{a}((\widetilde{\theta}_3,\varphi)^T,(\mathfrak{v},p)^T)= \langle \vect{y}_2,\mathfrak{v} \rangle_{L^2(\omega)}+\langle \vect{y}_3,p \rangle_{L^2(\Omega)}\\   
		&   \hspace{+1ex}+\fint_{\omega} \vect{y}_1 \cdot \fint_{\omega} \mathfrak{v} +\int_{\omega}\mathbb{A}^{\rm hom}(e_{\widehat{x}}(\mathfrak{a}_1^{\vect{y}_1}), \nabla^2_{\widehat{x}}\vect{y}_1):(0, \nabla^2_{\widehat{x}}\mathfrak{v}) d\widehat{x}, \quad  \forall (\mathfrak{v},p)^T \in N \times M. 
     \end{split}
 \end{equation}
	Here $\widetilde{a}$ is the bilinear form on $N \times M$ is defined in the following way
	\CCC
 \begin{equation} \label{tvuk1} 
         \widetilde{a}( (\widetilde{\theta}_3,\varphi)^T, (\mathfrak{v},p)^T):=a((\mathfrak{v},p)^T, (\widetilde{\theta}_3,\varphi)^T), \quad \forall
         (\widetilde{\theta}_3,\varphi)^T, (\mathfrak{v},p)^T \in N \times M.   
    \end{equation}
    \BBB
 \BBB
	As a consequence of \eqref{nedest1} and \eqref{subest1} we easily see that the operator $\mathcal{C}$ has trivial kernel. To see this take $(\mathfrak{b}_1, \widetilde{\theta}_3,p)^T=0$ in \eqref{nedest1}, \eqref{subest1} and conclude $\vect{y}=0$.  
	As a consequence of \eqref{utoest11}, \eqref{utoest20},\eqref{nedest1},\eqref{subest1} and \eqref{tvuk1} we have that
	$$
		\langle (I-\widetilde{\mathcal{A}})^{-1} \vect{x}, \vect{y} \rangle_H=\langle \vect{x}, \mathcal{C} \vect{y} \rangle_H, \quad \forall \vect{x},\vect{y} \in H. 	
	$$	
	This means that $\mathcal{C}=\left((I-\widetilde{\mathcal{A}})^{-1}\right)^T$. 
	The density of $\mathcal{D}(\widetilde{\mathcal{A}})$ in $H$ follows from the fact that the operator $\mathcal{C}$ has trivial kernel, which implies that the range of the operator $(I-\widetilde{\mathcal{A}})^{-1}$ is dense in $H$. 
	Directly from the definition we see that the operator $\widetilde{\mathcal{A}}$ is closed. Namely, if $\vect{x}_n \xrightarrow{H} \vect{x}$, $\widetilde{\mathcal{A}} \vect{x}_n  \xrightarrow{H} \vect{y}$, we have the following:
	$$ \vect{x} \xleftarrow{H} \vect{x}_n=(I-\widetilde{\mathcal{A}})^{-1}(I-\widetilde{\mathcal{A}})\vect{x}_n\xrightarrow{H} (I-\widetilde{\mathcal{A}})^{-1}(\vect{x}-\vect{y}).  $$
	From this it follows that $\vect{x} \in \mathcal{D}(\widetilde{\mathcal{A}})$ and $(I-\widetilde{\mathcal{A}}) \vect{x}=\vect{x}-\vect{y}$, which implies that $\widetilde{\mathcal{A}} \vect{x}=\vect{y}$. Furthermore, from \eqref{bukal2}, \eqref{bukal3}, \eqref{utoest11} we obtain that for every $(\mathfrak{b},\mathfrak{v},p)^T \in \mathcal{D}(\widetilde{\mathcal{A}})$ we have
	\begin{equation} 
		\label{krnic1} \langle \widetilde{\mathcal{A}} (\mathfrak{b},\mathfrak{v},p)^T, (\mathfrak{b},\mathfrak{v},p)^T\rangle_H =\fint_{\omega} \mathfrak{b} \, d\widehat{x}  \fint_{\omega} \mathfrak{v} \, d\widehat{x}- \int_{\Omega_p}\mathbb{K}_{33}(x_3)\partial_{3}p\,\partial_{3}  p \,dx.
	\end{equation} 	 
	where we used \eqref{nedest4} for $\widetilde{\theta}_3=\mathfrak{b}$ and \eqref{subest2} for $\widetilde{\theta}_3=\mathfrak{v}$. 
	From \eqref{defb11} we conclude that there exists $c>0$ such that  
		\begin{equation} 
		\label{krnic100} \langle (\widetilde{\mathcal{A}}-c\mathcal{B}) (\mathfrak{b},\mathfrak{v},p)^T, (\mathfrak{b},\mathfrak{v},p)^T\rangle_H <0.
	\end{equation}
Obviously  the operator 
	\begin{equation} \label{defA11} 
	\mathcal{A}_{c}:=\mathcal{B}^{-1} \widetilde{\mathcal{A}}-c I  :H \to H
	\end{equation} 
	 (where $H$ is equipped with either scalar product $\langle \cdot,\cdot \rangle_H$ or $\langle \cdot,\cdot \rangle_{\widetilde{b},H})$ is a closed, densely defined operator. From \eqref{krnic2} and \eqref{krnic100}  we obtain for every $(\mathfrak{b},\mathfrak{v},p)^T \in \mathcal{D}(\widetilde{\mathcal{A}})=\mathcal{D}( \mathcal{A}_{c})$  
	\begin{equation} 
		 \langle (\mathcal{B}^{-1}\widetilde{\mathcal{A}}-c I) (\mathfrak{b},\mathfrak{v},p)^T, (\mathfrak{b},\mathfrak{v},p)^T\rangle_{\tilde b,H} = \langle (\widetilde{\mathcal{A}}-c \mathcal{B}) (\mathfrak{b},\mathfrak{v},p)^T, (\mathfrak{b},\mathfrak{v},p)^T\rangle_{H}\leq0. \label{krnic555}
	\end{equation}
	From \eqref{krnic555} we conclude that for every $(\mathfrak{b},\mathfrak{v},p)^T \in \mathcal{D}(\widetilde{\mathcal{A}})$ and $\lambda>0$ we have
	\begin{equation} \label{krnic11} 
		\langle (\lambda I- \mathcal{A}_{c})(\mathfrak{b},\mathfrak{v},p)^T, (\mathfrak{b},\mathfrak{v},p)^T \rangle_{\tilde b,H}\geq \lambda \|(\mathfrak{b},\mathfrak{v},p)^T\|^2_{\tilde b,H}. 
	\end{equation} 
	From \eqref{krnic11} we easily conclude that for every $\lambda>0$: 
	\begin{equation} \label{krnic13} 
		\textrm{Ker} (\lambda I- \mathcal{A}_{c})=\{0\}. 
	\end{equation} 
By repeating the same analysis for $\widetilde{\mathcal{A}}^T$ (note that the identity \eqref{krnic1} is also valid for $\widetilde{\mathcal{A}}^T$) we analogously conclude that $ \label{krnic1333} 
		\textrm{Ker}(\lambda I-\mathcal{A}_{c})^T=\{0\},$ 
	where $(\cdot)^T$ is for $\widetilde{A}$ taken with respect to $\langle \cdot, \cdot\rangle_H$ while for $\lambda I-\mathcal{A}_{c}$ is taken with respect to $\langle \cdot,\cdot\rangle_{\tilde b,H}$ scalar product. This implies that the range of the operator $\lambda I- \mathcal{A}_{c} $ is dense. From \eqref{krnic11} we again conclude that the inverse of the operator $\lambda I-\mathcal{A}_{c}$ can be extended to continuous operator and that 
	\begin{equation} \label{krnic12} 
		\|(\lambda I-\mathcal{A}_{c})^{-1}\|_{(\tilde b,H)\to (\tilde b,H)}\leq \frac{1}{\lambda}, 
	\end{equation} 
	where $\|\cdot\|_{(\tilde b,H)\to(\tilde b,H)}$ denotes the operator norm taken with respect to $\|\cdot\|_{\tilde b,H}$ norm on $H$. Using the closedness of the operator $\lambda I- \mathcal{A}_{c}$ we conclude from \eqref{krnic12} that the operator $\lambda I- \mathcal{A}_{c}$ is surjective. Since by \eqref{krnic13} it is also injective, it is bijective and \eqref{krnic12} is valid. From Hille-Yosida theorem it follows that $ \mathcal{A}_{c}$ is a generator of a contraction semigroup.
	Thus the operator $\mathcal{A }:= \mathcal{A}_{c} +c I=\mathcal{B}^{-1}\widetilde{\mathcal{A}}$ is a generator of a strongly continuous semigroup. 
 
{\bf Step 3.} {\it Definition of the loads and the initial pressure.}	
	Obviously, from the construction the system \eqref{kreso33}-\eqref{kreso55} has the form \eqref{abstractevolutionproblem} where  
\begin{equation}\label{defloads} 
    \begin{split}
        		_{V'}\langle \vect{f} (t), (\theta_3,\widetilde{\theta}_3,\varphi)^T \rangle_{V} &= {_{N'}}\langle \vect{F}_2,\theta_3\rangle_{N}-\int_{\omega}\mathbb{A}^{\rm hom}(e_{\widehat{x}}(\mathfrak{a}_3^{\vect{F}_1}), 0):(0, \nabla^2_{\widehat{x}}\theta_3) d\widehat{x}+  {_{M'}} \langle\vect{G},\varphi\rangle_{M}  \\ 
		&  - \int_{\omega}\int_{J_p}(|\mathcal{Y}_{f}(x_3)| \mathbb{I}-\mathbb{B}^{H}(x_3))\varphi \, dx_3:\iota\left(e_{\widehat{x}}\left(\mathfrak{a}_3^{\partial_t \vect{F}_1(t)}\right)\right) d\widehat{x},
    \end{split}
\end{equation}
	and for $\widehat{\vect{F}} \in L'_{1} $, $\mathfrak{a}_3^{\widehat{\vect{F}}} \in L_1$ is a unique  solution of 
	\begin{equation} 
		\label{loads1} \displaystyle\int_{\omega}\mathbb{A}^{\rm hom}(e_{\widehat{x}}(\mathfrak{a}_3^{\widehat{\vect{F}}}),0):(e_{\widehat{x}}(\vect{\theta}_*), 0)d\widehat{x}= \displaystyle {_{L_1'}} \langle  \widehat{\vect{F}},\vect{\theta}_*\rangle_{L_1},
		\quad\forall \vect{\theta}_{*} \in L_1.
	\end{equation}
	The initial condition is given by 
	\begin{equation} \label{defin} 
	(\mathfrak{b}(0),\mathfrak{v}(0),p(0) )^{\top}=(\mathfrak{b}_0,\mathfrak{b}_1, p_0)^{\top},  
	\end{equation} 
	where $p_0 \in L^2(\Omega_p)$ is given by 
\begin{equation}\label{krnj111} 
         \widetilde{b} (p_0,h)  =\int_{\omega}\int_{J_p}(|\mathcal{Y}_{f}(x_3)| \mathbb{I}-\mathbb{B}^{H}(x_3))h \, dx_3:\iota\left(e_{\widehat{x}}\left(\mathfrak{a}_3^{ \vect{F}_1(0)}\right)\right) d\widehat{x} +\int_{\Omega_p} t_0\cdot h\,dx, \quad \forall h \in L^2(\Omega_p). 	
 \end{equation}
\CCC To see this one needs to conclude from \eqref{kreso44} that  
$$\mathfrak{a}=\mathfrak{a}_1^{\mathfrak{b}}+\mathfrak{a}_2^p+\mathfrak{a}_3^{\vect{F}_1}, $$
plug this in the equations \eqref{kreso44} and \eqref{kreso55}, do integration by parts in time to move the time derivative from the test functions $\varphi$, $\theta_3$ and introduce $\mathfrak{v}:=\partial_t \mathfrak{b}$ (this is standard when we want to write the second order problems as first order) and use \eqref{utoest3}, \eqref{defb11}, \eqref{bukal1}-\eqref{bukal3} to write \eqref{kreso33}-\eqref{kreso55} in the form:
\begin{eqnarray*}
 \langle \partial_t (\mathfrak{b},\mathfrak{v},p)^T,(\theta_3,\widetilde{\theta}_3,\varphi)^T\rangle_{\widetilde{b},H}&=&\langle \mathcal{B}\partial_t  (\mathfrak{b},\mathfrak{v},p)^T,(\theta_3,\widetilde{\theta}_3,\varphi)^T\rangle_{H} \\
 &=& \langle \widetilde{\mathcal{A}}(\mathfrak{b},\mathfrak{v},p)^T,(\theta_3,\widetilde{\theta}_3,\varphi)^T \rangle_H+_{V'}\langle \vect{f} (t), (\theta_3,\widetilde{\theta}_3,\varphi)^T \rangle_{V}\\ &=&\langle \mathcal{A}(\mathfrak{b},\mathfrak{v},p)^T,(\theta_3,\widetilde{\theta}_3,\varphi)^T \rangle_{\widetilde{b},H}+_{V'}\langle \vect{f} (t), (\theta_3,\widetilde{\theta}_3,\varphi)^T \rangle_{V},
 \end{eqnarray*} 
 with \eqref{defin}. This finishes the proof. 
   \BBB
\end{proof} 	

 The following theorem gives us existence of the solution of \eqref{kreso33}-\eqref{kreso55}. 
\begin{teorem} \label{existinertia} 
	If  $\mathfrak{b}_0 \in H^2_{\#}(\omega)$, $\mathfrak{b}_1 \in L^2(\omega)$, $t_0 \in L^2(\Omega_p)$,  $\vect{F}_{1} \in H^1(0,T; L'_{1})$, $\vect{F}_2 \in L^2(0,T;N')$, $\vect{G} \in L^2(0,T;M')$,  then there exists a solution $\mathfrak{a} \in L^2(0,T;L_1)   $, $\mathfrak{b} \in L^2(0,T;N)\cap H^1(0,T;L^2(\omega))$, $p \in L^2(0,T; M)$  of \eqref{kreso3}-\eqref{kreso5}. In addition we have that $\mathfrak{a} \in C([0,T];L_1)$, $\mathfrak{b} \in C([0,T]; N) \cap C^1([0,T];L^2(\omega)) $, $p \in C([0,T];L^2(\Omega))$.  Moreover there exists
	$C>0$ such that 
 \begin{equation}\label{kreso56}
     \begin{split} &\|\mathfrak{a}\|_{C(0,T;L_1)}+\|\mathfrak{b}\|_{C(0,T;N)}+\|\mathfrak{b}\|_{H^1(0,T;L^2(\omega))}+\|p\|_{C(0,T; L^2(\Omega))}+\|p\|_{L^2(0,T; M)}\\   &  \leq C  \Big(\|\vect{F}_1\|_{H^1(0,T; L'_{1})}+\|\vect{F}_2\|_{L^2(0,T; N')}+\|\vect{G}\|_{L^2(0,T; M')}  +\|\mathfrak{b}_0\|_{N}+\|\mathfrak{b}_1\|_{L^2(\omega)}+\|t_0\|_{L^2(\Omega_p)} \Big).
     \end{split}
 \end{equation}	
\end{teorem} 	 
\begin{proof} 
 We use Galerkin approximation. We take a sequence of increasing finite dimensional subspaces $(G_n)_{n \in \mathbb{N}}$,  $G_n=G_n^1 \times G_n^2 \times G_n^3$ where $G_n^1\subset L_1$,  $G_n^2 \subset  N$, $G_n^3 \subset M$, for every $n \in \mathbb{N}$, and $\overline{\cup_{n \in \mathbb{N}}G_n^1}=L_1$,  $\overline{\cup_{n \in \mathbb{N}}G_n^2}= N$,  $\overline{\cup_{n \in \mathbb{N}}G_n^3}=M$.  We approximate the initial condition $\mathfrak{b}_0$ with the sequence  $(\mathfrak{b}_0^n)_{n \in \mathbb{N}}$ such that $ G_n^2\ni \mathfrak{b}_0^n  \to \mathfrak{b}_0$ in  N. We consider the solution of the problem: For  $\mathfrak{b}_1 \in L^2(\omega)$, $t_0 \in L^2(\Omega_p)$,   $\vect{F}_{1} \in H^1(0,T; L'_{1})$, $\vect{F}_2 \in L^2(0,T;N')$, $\vect{G}_1 \in H^1(0,T;M')$, $\vect{G}_2 \in L^2(0,T;L^2(\Omega))$ we find $\mathfrak{a}_n \in H^1(0,T;G_n^1)$, $\mathfrak{b}_n \in H^2(0,T;G_n^2)$, $p_n \in H^1(0,T; G_n^3)$ which satisfy \eqref{kreso33} with $\mathfrak{b}_0$  replaced with $\mathfrak{b}_0^n$, \eqref{kreso44} for every $(\vect{\theta}_*,\theta_3) \in H^1(0,T;G_n^1 \times G_n^2)$, and \eqref{kreso55} for every $\varphi \in H^1(0,T;G_n^3)$. The existence of such problem is guaranteed by the fact that the problem has the form \eqref{abstractevolutionproblem} (see Definition \ref{sol1} and Remark \ref{pedro1}), where $H=G_n^2 \times G_n^2 \times G_n^3$, $V=H$, and the operator $\mathcal{A}$ is given by (matrix) $\mathcal{B}_n^{-1} \widetilde{\mathcal{A}}_n$, where $\widetilde{\mathcal{A}}_n$ and $\mathcal{B}_n$ are defined with \eqref{bukal1}, \eqref{bukal2}, \eqref{bukal3} and \eqref{utoest1}, \eqref{defb11}, but restricted and with test functions on $G_n^2$ and $G_n^3$ respectively  (we take $H$ equipped with the scalar product $\langle \cdot,\cdot\rangle_{\widetilde{b},H}$ defined in the proof of Proposition \ref{propsemigroup}). 
 Since we are on a finite dimensional space the definition of $\widetilde{\mathcal{A}}_n$ with 
 \eqref{bukal1}, \eqref{bukal2} and \eqref{bukal3} is not only formal. 
 Also  $\mathcal{D}(\mathcal{B}_n^{-1} \widetilde{\mathcal{A}}_n)=G_n^2 \times G_n^2 \times G_n^3$ and thus we have additional regularity in time of $\mathfrak{a}_n,\mathfrak{b}_n,p_n$. This enables us to modify \eqref{kreso33}-\eqref{kreso55} by doing integration by parts and conclude
\begin{equation}\label{kreso333}
    \begin{split}
        \mathfrak{b}_n(0)=\mathfrak{b}_0^n, \quad \partial_t \mathfrak{b}_n(0)= \mathfrak{b}_1^n, \quad p_n(0)=p_0^n,
    \end{split}    
\end{equation}
\begin{equation}\label{kreso444}
    \begin{split}
        & \int_{\omega} \bar{\kappa}(\widehat{x})\partial_{tt} \mathfrak{b}_n \theta_3 \,d\widehat{x}+ \displaystyle\int_{\omega}\mathbb{A}^{\rm hom}(e_{\widehat{x}}(\mathfrak{a}_n), \nabla^2_{\widehat{x}}\mathfrak{b}_n):(e_{\widehat{x}}(\vect{\theta}_*), \nabla^2_{\widehat{x}}\theta_3) d\widehat{x}\\ &-\displaystyle\int_{\omega}\int_{J_p}\left(|\mathcal{Y}_f(x_3)| \mathbb{I}-\mathbb{B}^{H}(x_3)\right)p_n\, dx_3:\left[\iota (e_{\widehat{x}}(\vect{\theta}_{*}))\right]\,d\widehat{x}+\displaystyle\int_{\omega} \int_{J_p}\left(|\mathcal{Y}_f(x_3)| \mathbb{I}-\mathbb{B}^{H}(x_3)\right)x_3 p \, dx_3:\left[\iota (\nabla^2_{\widehat{x}}\theta_3)\right]\,d\widehat{x} \\
	&=\displaystyle{_{L'_{1}}} \langle \vect{F}_1,\vect{\theta}_*\rangle_{L_1}+{_{N'}} \langle \vect{F}_2,\theta_3\rangle_{N},\quad\forall \left(\vect{\theta}_{*},{\theta}_3\right) \in G_n^1 \times G_n^2, \textrm{ for a.e. t} \in (0,T); 
    \end{split}    
\end{equation}
\begin{equation}\label{kreso555} 
    \begin{split}
      &\displaystyle\int_{\Omega_p}M_0(x_3)\partial_t p_n \, \varphi dx+\displaystyle\int_{\Omega_p}\mathbb{K}_{33}(x_3)\partial_{3}p_n\,\partial_{3}  \varphi dx +\displaystyle\int_{\omega}\int_{J_p}(|\mathcal{Y}_{f}(x_3)| \mathbb{I}-\mathbb{B}^{H}(x_3))\varphi \, dx_3:\iota\left(e_{\widehat{x}}\left(\partial_t \mathfrak{a}_n\right)\right) d\widehat{x}\,dt\\
	&-\displaystyle\int_{\omega}\int_{J_p}\left(|\mathcal{Y}_f(x_3)|\mathbb{I}-\mathbb{B}^H(x_3)\right)x_3\varphi \, dx_3:\iota\left(\nabla^2_{\widehat{x}}\partial_t \mathfrak{b}_n\right)\,d\widehat{x} ={_{M'}} \langle \vect{G},\varphi\rangle_{M}  \quad \forall  \varphi \in G_n^3 \textrm{ for a.e. } t \in (0,T). 
    \end{split}    
\end{equation}
Here $\mathfrak{b}_1^n$ is the orthogonal projection of $\mathfrak{b}_1$ on $G_n^1$ with respect to the $\langle \cdot,\cdot\rangle_{L^2(\omega)}$ scalar product, while $p_0^n \in G_n^3$ is given by the expression (cf. \eqref{krnj111}):
\begin{equation}\label{krnj1111}
    \begin{split}
        \widetilde{b} (p_0^n,h)  =\int_{\omega}\int_{J_p}(|\mathcal{Y}_{f}(x_3)| \mathbb{I}-\mathbb{B}^{H}(x_3))h \, dx_3:\iota\left(e_{\widehat{x}}\left(\mathfrak{a}_3^{ \vect{F}_1(0)}\right)\right) d\widehat{x} +\int_{\Omega_p} t_0\cdot h\,dx,   \quad \forall h \in G_n^3. 
    \end{split}
\end{equation}
By testing then \eqref{kreso444} with  $(\partial_{t}\mathfrak{a}_n,\partial_t\mathfrak{b}_n)$ and  \eqref{kreso555} with $\varphi=p_n$ and integrating over $(0,t)$ for $t<T$ we obtain 
\begin{equation}\label{vesna1}
    \begin{split}
       & \|\bar{\kappa}\partial_t \mathfrak{b}_n(t)\|^2_{L^2(\omega)}+\int_{\omega}\mathbb{A}^{\rm hom}(e_{\widehat{x}}(\mathfrak{a}_n(t)), \nabla^2_{\widehat{x}}\mathfrak{b}_n(t)):(e_{\widehat{x}}(\mathfrak{a}_n(t)), \nabla^2_{\widehat{x}}\mathfrak{b}_n(t)) d\widehat{x}\\   &\quad+\int_{\Omega_p} M_0(x_3) p_n^2(t)\,dx+\int_0^{t}\int_{\Omega_p} \mathbb{K}_{33}(x_3)(\partial_{3}p_n)^2  \, dx\, d\tau=\|\bar{\kappa}\mathfrak{b}_1^n\|_{L^2(\omega)}\\  &\quad+\int_{\omega} \mathbb{A}^{\rm hom}(e_{\widehat{x}}(\mathfrak{a}_n(0)), \nabla^2_{\widehat{x}}\mathfrak{b}_n(0)):(e_{\widehat{x}}(\mathfrak{a}_n(0)), \nabla^2_{\widehat{x}}\mathfrak{b}_n(0)) d\widehat{x} + \int_0^{t} {_{L'_{1}}}\langle \partial_t \vect{F}_1, \mathfrak{a}_n \rangle_{L_1}\,d\tau\\  & \quad -{_{L'_{1}}}\langle  \vect{F}_1(0), \mathfrak{a}_n (0)\rangle_{L_1}+\int_0^{t}{_{N'}} \langle \vect{F}_2,\mathfrak{b}_n\rangle_{N}\,d\tau+\int_{\Omega_p} M_0(x_3) (p_0^n)^2 \,dx+\int_0^{t}{_{M'}} \langle \vect{G},p_n\rangle_{M}\,d\tau. \\ 
    \end{split}
\end{equation}	 
To obtain the equation for $\mathfrak{a}_n(0)$ from $\mathfrak{b}_n(0)$ and $p_n(0)$ one can test \eqref{kreso444} with $\theta_3=0$ for $t=0$. After doing Young's inequality on the terms of the right side of equality sign of \eqref{vesna1} and by using Proposition \ref{tensorprop} one obtains the bound: 
\begin{equation}\label{andrei2}
    \begin{split}
        &\|\partial_t \mathfrak{b}_n\|_{L^\infty(0,T;L^2(\omega))}+\|\mathfrak{a}_n\|_{L^\infty(0,T;L_1)}
+\|\mathfrak{b}_n\|_{L^{\infty}(0,T;N)}+\|p_n\|_{L^{\infty}(0,T;L^2(\Omega_p))} +\|p_n\|_{L^2(0,T;M)} \\  &  \leq C\big( \|\vect{F}_1\|_{H^1(0,T;L'_{1})}+\|\vect{F}_2\|_{L^2(0,T;N')}+\|\vect{G}\|_{L^2(0,T;M')}   +\|\mathfrak{b}_0^n\|_{N}+\|\mathfrak{b}_1^n\|_{L^2(\omega)}+\|p_0^n\|_{L^2(\Omega_p)}\big),
    \end{split}
\end{equation}	
for some $C>0$. 
From \eqref{andrei2} we conclude that there exists weak limit $$(\mathfrak{a},\mathfrak{b},p) \in \left(L^\infty(0,T;L) \cap W^{1,\infty}(0,T;L^2(\omega))\right)\times L^{\infty}(0,T;N)\times \left(L^{\infty}(0,T;L^2(\Omega))\cap L^2(0,T;M)\right)$$ of $(\mathfrak{a}_n, \mathfrak{b}_n,p_n)$. By letting $n$ to infinity in \eqref{kreso333}-\eqref{kreso555} (after writing them in the form of \eqref{kreso44} and \eqref{kreso55}) we obtain that $(\mathfrak{a},\mathfrak{b},p)$ satisfies \eqref{kreso33}-\eqref{kreso55}.  By letting $n$ to infinity in \eqref{andrei2} and by using \eqref{krnj1111} we obtain
\begin{equation}\label{andrei3} 
    \begin{split}
        &\|\partial_t \mathfrak{b}\|_{L^\infty(0,T;L^2(\omega))}+\|\mathfrak{a}\|_{L^\infty(0,T;L_1)}
	+\|\mathfrak{b}\|_{L^{\infty}(0,T;N)}+\|p\|_{L^{\infty}(0,T;L^2(\Omega_p))} +\|p\|_{L^2(0,T;M)} \\  & \leq C\big( \|\vect{F}_1\|_{H^1(0,T;L'_{1})}+\|\vect{F}_2\|_{L^2(0,T;N')}+\|\vect{G}\|_{L^2(0,T;M')}  +\|\mathfrak{b}_0\|_{N}+\|\mathfrak{b}_1\|_{L^2(\omega)}+\|t_0\|_{L^2(\Omega_p)}\big). 
    \end{split}
\end{equation}	
The additional regularity properties follow again from   \cite[Remark 1, Chapter XVIII.5]{dautraylions} (cf. Remark \ref{remkrnj1}) and thus we can replace $L^{\infty}([0,T])$ norms  by $C([0,T]$ on the left hand side of \eqref{andrei3}.   
\end{proof}
Next we prove the energy-type identity, \CCC which is again needed to prove the strong convergence. \BBB

\begin{propozicija} \label{propen2}
	If $\vect{F}_1 \in H^1(0,T;L'_{1})$, $\vect{F}_2 \in L^2(0,T;N')$, $\vect{G} \in L^2(0,T;M')$, $\mathfrak{b}_0 \in N$, $\mathfrak{b}_1 \in L^2(\Omega)$, $t_0 \in L^2(\Omega)$ and $\mathfrak{a} \in C(0,T;L_{1})$, $\mathfrak{b} \in C(0,T;N)\cap H^1(0,T;L^2(\omega))$, $p \in L^2(0,T; M)\cap C(0,T;L^2(\Omega))$ is the (unique) solution of \eqref{kreso333}-\eqref{kreso555}, then we have  for every $t \in [0,T]$:
\begin{equation}\label{vesna2}
    \begin{split}
        & 
	\|	\bar{\kappa}\partial_t \mathfrak{b}(t)\|^2_{L^2(\omega)}+\int_{\omega}\mathbb{A}^{\rm hom}(e_{\widehat{x}}(\mathfrak{a}(t)), \nabla^2_{\widehat{x}}\mathfrak{b}(t)):(e_{\widehat{x}}(\mathfrak{a}(t)), \nabla^2_{\widehat{x}}\mathfrak{b}(t)) d\widehat{x}\\  &\quad+\int_{\Omega_p} M_0(x_3) p^2(t)\,dx+\int_0^{t}\int_{\Omega_p} \mathbb{K}_{33}(x_3)(\partial_{3}p)^2  \, dx\, d\tau=\|\bar{\kappa}\mathfrak{b}_1\|_{L^2(\omega)}\\  &\quad+\int_{\omega} \mathbb{A}^{\rm hom}(e_{\widehat{x}}(\mathfrak{a}(0)), \nabla^2_{\widehat{x}}\mathfrak{b}(0)):(e_{\widehat{x}}(\mathfrak{a}(0)), \nabla^2_{\widehat{x}}\mathfrak{b}(0)) d\widehat{x} + \int_0^{t} {_{L_{1}'}}\langle \partial_t \vect{F}_1, \mathfrak{a} \rangle_{L_{1}}\,d\tau\\  & \quad -{_{L_1'}}\langle  \vect{F}_1(0), \mathfrak{a} (0)\rangle_{L_{1}}+\int_0^{t}{_{N'}} \langle \vect{F}_2,\mathfrak{b}\rangle_{N}\,d\tau+\int_{\Omega_p} M_0(x_3) (p_0)^2 \,dx+\int_0^{t}{_{M'}} \langle \vect{G},p\rangle_{M}\,d\tau. 
    \end{split}
\end{equation}
Here $p_0 \in L^2(\Omega_p)$ is defined by \eqref{krnj111}. 
\end{propozicija} 	
\begin{proof} 
If $\mathfrak{a} \in H^1(0,T;L_{1})$, $\mathfrak{b} \in H^1(0,T;N) \cap H^2(0,T;L^2(\omega))$ and $p \in L^2(0,T;M) \cap H^1(0,T;L^2(\Omega))$ we would obtain \eqref{vesna2} in the same way as we obtained \eqref{vesna1}. For the general case we use approximation. We know that $\mathcal{D}(\mathcal{A})$, where $\mathcal{A}$ is defined in the proof of Proposition \ref{propsemigroup}, is  dense in $H$ defined by \eqref{defH11} and $H$ is dense in $V'$, where $V$ is defined by \eqref{defV11}. From this we also have that $L^2(0,T;\mathcal{D}(\mathcal{A}))$ is dense in $L^2(0,T;V')$. The claim for general $(\mathfrak{a},\mathfrak{b},p)$ follows by approximation argument by approximating  with  $\vect{F}_2$, $\vect{G}$ such that $\vect{f}$ defined by \eqref{defloads} is in $L^2(0,T;\mathcal{D}(\mathcal{A}))$ and approximating $\mathfrak{b}_0,\mathfrak{b}_1,t_0$ such that the initial conditions \eqref{defin} is in $\mathcal{D}(\mathcal{A})$.  
By using Remark \ref{pedro1} we know that the strong solution of \eqref{abstractevolutionproblem} has the above required regularity (with this approximated loads and initial condition) and solves \eqref{kreso333}-\eqref{kreso555} and thus satisfies \eqref{vesna2}. The claim for the original loads follows then by using stability estimate \eqref{kreso56}.  
\end{proof} 	
 \subsection{Strong convergence} \label{strong2}

We will state the convergence result which is analogous to Theorem \ref{proplenka1} and Theorem \ref{twoscaleconvergencesmay16}. The proof goes in an analogous way and we will skip it.
\begin{teorem} \label{stronginertia} 
 Let assumptions \eqref{kirill3} and \eqref{forcesassumptions0}  be satisfied, together with the Assumptions \ref{assumption on regions}. Furthermore, we assume that $
		\pi_h\vect{F}^{h}\xrightarrow{t,2-r,2} \vect{F}, \quad h\partial_t F_{\alpha}^{h}\xrightarrow{t,2-r,2} \partial_t F_{\alpha}, \textrm{ for } \alpha=1,2, $ 
	where $\vect{F} \in H^1(0,T; L^2(\Omega \times \mathcal{Y};\mathbb{R}^2)) \times  L^2(0,T; L^2(\Omega \times \mathcal{Y})) $ and 
	$F_{\alpha}(0)=0$, for $\alpha=1,2$. 
	Let $(\vect{u}^h,p^h)$ be the solution of \eqref{varformporous2gamma0may10} with initial conditions \eqref{eq:51rescaledver}. Than, the convergences \eqref{strong2may16} and \eqref{strong1may16} are valid as well as
$$	\lim_{h\to 0}\int_0^T\int_{\Omega}\kappa^{h}(x)|\partial_{t}\vect{u}^h(t)|^2\,dx\, dt =\int_0^T\int_{\omega}\bar{\kappa}(\widehat{x})\left|\partial_t\mathfrak{b}\right|^2\,d\widehat{x}\, dt. $$
Furthermore, convergences \eqref{strongtwoscale1may16} and \eqref{strongtwoscale4may16} are valid as well as: 
$$ \lim_{h\to 0}\int_0^T\int_{\Omega}|\partial_{t}\vect{u}^h(t)-(0,0,\partial_t \mathfrak{b}(t))^T|^2\,dx\, dt=0. $$
Here $(\mathfrak{a},\mathfrak{b}, p) \in L^2(0,T;L_{\kappa,\mathfrak{b}})  \times\left(L^2(0,T;N)\cap H^1(0,T;L^2(\omega))\right)\times L^2 (0,T;M)$ are solutions of \eqref{kreso3}-\eqref{kreso5} and $\vect{u}^f_0 \in  H^1(0,T;$ $L^2(\Omega; H^1(\mathcal{Y}_f(x_3);\RR^3)))$ is defined with \eqref{diferformfluid1march301}.
\end{teorem}

\begin{appendices}
\section*{Appendix}

\section{Griso's decomposition and its consequences}
\label{griso} 
In this section we use Griso's decomposition to characterize the sequences with bounded symmetrized scaled gradients. \CCC These results present fundamental tools to deal with dimension reduction problems, since in these problems, in the {\it a priori} estimates, one only has the boundedness of symmetrized scaled gradients, which has to be combined with boundary conditions (and we show how they can be combined with Dirichlet or periodic boundary conditions).
The analogue results were shown in \cite{Vel14a,mbukalvelcic2017,Marinvelciczubrinic2022}. Here we need some additional claims where we deal with periodic boundary conditions. 
 
Griso's decomposition tells us how on $h$ level we can decompose sequences with bounded elastic energy (on thin domain) and gives estimate of the reminder. This result has to be combined with certain compactness statements with respect two-scale convergence to fully characterize possible  two-scale limits of sequences of symmetrized scaled gradients. 

 \BBB

Let $\omega \subset \R^2$ be an open set and let us denote $\Omega = \omega \times I$.
Let $\gamma_D\subset \partial \omega$ be of positive measure and let $\Gamma_D=\gamma_D \times I$. We denote by $H^1_{\Gamma_D}(\Omega)$ the subspace of $H^1(\Omega)$ containing functions  with zero trace on $\Gamma_D$. Analogously we define $H^1_{\gamma_D}(\omega)$. Similarly, $H^2_{\gamma_D}(\omega)$ denotes the subspace of $H^2(\omega)$ such that the trace of the functions and its first derivatives on $\gamma_D$ is zero.

First we state Griso's decomposition. \CCC  Although in \cite{Grisosdecomp2005} it is stated on thin domain, we state the result here on the canonical domain. \BBB
\begin{teorem}[Griso's decomposition, \cite{Grisosdecomp2005}, Theorem 2.3]\label{aux:thm.griso}
	Let $ \omega\subset\R^2$ be a set with Lipschitz boundary and $\mat{\psi}\in H^1(\Omega,\R^3)$. Then for arbitrary $h \in (0,1)$ 
	the following identity holds
	\begin{equation}\label{griso1}
		\mat{\psi} = \widehat{\mat{\psi}}(\widehat{x}) + \vect{r}(\widehat{x}) \wedge x_3\vect{e}_3 + \overline{\mat{\psi}}(x)
		= \left\{   \begin{array}{l}
			\widehat{\psi}_1(\widehat{x}) + r_2(\widehat{x})x_3 + \overline{\psi}_1(x)\\
			\widehat{\psi}_2(\widehat{x}) - r_1(\widehat{x})x_3 + \overline{\psi}_2(x)\\
			\widehat{\psi}_3(\widehat{x}) + \overline{\psi}_3(x)
		\end{array}
		\right.\,,
	\end{equation}
	where
	\begin{equation} \label{griso2}
		\widehat{\mat\psi}(\widehat{x}) = \int_I \mat\psi(\widehat{x},x_3)\dd x_3\,,\quad \vect{r}(x') 
		= \CCC 12 \BBB\int_Ix_3\vect{e}_3\wedge\mat\psi(\widehat{x},x_3)\dd x_3\,.
	\end{equation}
	Moreover, the following inequality holds
	\begin{equation}\label{prvaKorn}
		\|e_h(\widehat{\mat\psi} + \vect{r}\wedge x_3\vect{e}_3)\|_{L^2(\Omega;\RR^{3 \times 3})}^2 
		+ \|\nabla_h\overline{\mat\psi}\|_{L^2(\Omega;\RR^{3 \times 3})}^2 + \frac{1}{h^2}\|\overline{\mat\psi}\|_{L^2(\Omega;\RR^3)}^2
		\leq C\|e_h(\mat\psi)\|_{L^2(\Omega;\RR^{3 \times 3})}^2\,,
	\end{equation} 
	with constant $C>0$ depending only on $\omega$.
\end{teorem} 
\begin{remarkica} 
We have that 
\begin{equation} \label{remnak10}
    \begin{split}
        \|e_h(\widehat{\mat\psi} + \vect{r}\wedge x_3\vect{e}_3)&\|_{L^2(\Omega;\RR^{3 \times 3})}^2
        =\|e_{\widehat{x}}(\widehat{\psi}_1,\widehat \psi_2)\|_{L^2(\omega;\RR^{2 \times 2})}^2 + \frac{1}{12}\|e_{\widehat{x}}(r_2,-r_1)\|_{L^2(\omega;\RR^{2 \times 2})}^2
        \\& + \frac{1}{h^2}\|\pa_1(h\widehat{\psi}_3) + r_2\|_{L^2(\omega)}^2 + \frac{1}{h^2}\|\pa_2(h\widehat{\psi}_3) - r_1\|_{L^2(\omega)}^2. 
    \end{split}
\end{equation}
\end{remarkica} 	
The following lemma is a direct consequence of Griso's decomposition and is proven in \CCC \cite[Lemma A4]{mbukalvelcic2017}. 
It decomposes the sequence which has bounded symmetrized scaled gradients in the "limit deformation" that follows  Kirchoff-Love ansatz and the remainder, which is further characterized in Lemma  \ref{lemmaA10}. 
\BBB

\begin{lema}\label{lemmaA8}
Consider a bounded set $\omega\subset\RR^2$ with Lipschitz boundary. Suppose that $(\vect{\psi}^h)_{h \in (0,1)}\subset H^1_{\Gamma_D}(\Omega;\RR^3)$ is such that
$
\displaystyle\limsup_{h\to 0}\|e_h (\vect{\psi}^{h})\|_{L^2(\Omega;\RR^{3\times 3})}<\infty. 
$
Then, there exists a subsequence (still labeled by $h>0$) such that
\begin{equation}
\vect{\psi}^{h}=(
      \mathfrak{a}_1 -x_3\partial_1\mathfrak{b}, 
       \mathfrak{a}_2-x_3\partial_2\mathfrak{b}, 
       h^{-1}\mathfrak{b})^T+\widetilde{\vect{\psi}}^{h},
\end{equation}
for some $\mathfrak{a}\in H^1_{\gamma_D}(\omega;\RR^2)$, $\mathfrak{b}\in H^2_{\gamma_D}(\omega)$ and a sequence $(\widetilde{\vect{\psi}}^{h})_{h \in (0,1)}\subset H^1_{\Gamma_D}(\Omega;\RR^3)$, which satisfies $h\pi_{1/h}\widetilde{\vect{\psi}}^{h}\overset{L^2}{\to} 0$. 
In particular,
$
e_h(\vect{\psi}^h)=\iota(e_{\widehat{x}}(\mathfrak{a})-x_3\nabla^2_{\widehat{x}}\mathfrak{b})+e_h(\widetilde{\vect{\psi}}^{h}).
$

\end{lema}
\begin{remarkica} \label{remnak1}
For the case when $\omega$ is rectangle and $(\vect{\psi}^h)_{h \in (0,1)}\subset H^1_{\#}(\Omega;\RR^3)$ \CCC satisfies that $h \int_{\Omega} \pi_{1/h}(\kappa^h \vect{\psi}^h) $ is bounded for some $(\kappa^h)_{h \in (0,1)} \subset L^{\infty}(\Omega)$, which is uniformly bounded (independently of $h$) from below and above by a positive constant\BBB, we have that   $\mathfrak{b}\in H^2_{\#}(\omega)$, $\mathfrak{a}\in H^1_{\#}(\omega;\RR^2)$ and $(\widetilde{\vect{\psi}}^{h})_{h \in (0,1)}\subset H^1_{\#}(\Omega;\RR^3)$. This can be shown following the proof of Lemma \ref{lemmaA8} \CCC and using the result of Proposition \ref{propapp1} \BBB. 
\end{remarkica} 	

The following is a part of \cite[Lemma A10]{Marinvelciczubrinic2022}. \CCC Together with Lemma \ref{lemmaA8} it gives the desired compactness result for the sequence with bounded symmetrized scaled gradients, which is stronger than the one that would follow from Theorem \ref{aux:thm.griso}.
 However, it requires $C^{1,1}$ regularity of the domain (or periodic boundary conditions, see Remark \ref{remnak2} below). That is the reason why we still needed Theorem \ref{aux:thm.griso} in \cite{Vel13,mbukalvelcic2017,Marinvelciczubrinic2022}
and combined it with Lemma \ref{lemmaA10} on $C^{1,1}$ regular subsets of Lipschitz regular domain.  
 \BBB
\begin{lema}\label{lemmaA10}
Let $\omega\subset\RR^2$ be  a connected set with $C^{1,1}$ boundary.  If $(\vect{\psi}^h)_{h>0}\subset H^1(\Omega;\RR^3)$ is such that
\begin{equation}
h\pi_{1/h}\vect{\psi}^h\overset{L^2}{\to}0,\quad \displaystyle \limsup_{h\to 0}\|e_h (\vect{\psi}^{h})\|_{L^2(\Omega;\RR^{3\times 3})}\leq M<\infty,
\end{equation}
then there exist $(\varphi^h)_{h>0}\subset H^2(\omega)$, $(\widetilde{\vect{\psi}}^h)_{h>0}\subset H^1(\Omega;\RR^3)$ such that
\begin{equation}
e_h(\vect{\psi}^h)=-x_3\iota(\nabla^2_{\widehat{x}}\varphi^h)+e_h(\widetilde{\vect{\psi}}^h)+o^h,
\end{equation}
where $o^h\in L^2(\Omega;\RR^{3\times 3})$ is such that $o^h\overset{L^2}{\to}0$. In addition, the following properties hold:
\begin{equation}
\displaystyle\lim_{h\to 0}\left(\|\varphi^h\|_{H^1(\omega)}+\|\widetilde{\vect{\psi}}^h\|_{L^2(\Omega;\RR^3)}\right)=0,\quad \displaystyle\limsup_{h\to 0}\left(\|\varphi^h\|_{H^2(\omega)}+\|\nabla_h\widetilde{\vect{\psi}}^h\|_{L^2(\Omega;\RR^{3 \times 3})}\right)\leq CM,
\end{equation}
where $C$ depends only on $\omega$.
\end{lema}
\begin{remarkica} \label{remnak2}
	In the case when $\omega$ is a rectangle and $(\vect{\psi}^h)_{h>0}\subset H_{\#}^1(\Omega;\RR^3)$, the statement of Lemma \ref{lemmaA10} holds with  $(\varphi^h)_{h \in (0,1)}\subset H_{\#}^2(\omega)$, $ (\widetilde{\vect{\psi}}^h)_{h \in (0,1)}\subset H_{\#}^1(\Omega;\RR^3)$. This can be proved following the proof of \cite[Lemma A10]{Marinvelciczubrinic2022}.
\end{remarkica}

Using Theorem \ref{aux:thm.griso} we  prove the following statement. \CCC It tells us that the components of the sequence with bounded symmetrized scaled gradients that satisfy periodic boundary condition, after appropriately scaled, are bounded in $H^1$ norm.  
This scaling is standard in the context of dimension reduction problems in linearized elasticity. 
Moreover, in order to obtain $H^1$ bondedness of the scaled components, in addition one only needs to have the control on the scaled weighted mean value of the components.  
\BBB
\begin{propozicija} \label{propapp1} 
Let $\omega \subset \R^2$ be a rectangle.  Let $(\kappa^h)_{h \in (0,1)} \subset L^\infty(\Omega)$, be such that $c_0 < \kappa^h < c_1$ almost everywhere, for some $c_0,\,c_1>0$ independent of $h$. Then there exists $C>0$ such that for every $(\vect{\psi}^h)_{h \in (0,1)}\subset H^1_{\#}(\Omega;\RR^3)$ we have 
\begin{equation} \label{appnalk} 
	\|\pi_{1/h} \vect{\psi}^h\|_{H^1(\Omega;\RR^3)} \leq C\left(\frac{1}{h} \|e_h(\vect{\psi}^h)\|_{L^2(\Omega;\RR^{3 \times 3})}+\left|\int_{\Omega} \pi_{1/h} (\kappa^h\vect{\psi}^h) \right| \right). 
\end{equation} 	
\end{propozicija}
\begin{proof} 
We assume the opposite, that for every $n \in \mathbb{N}$, there exists a sequence $(h_n)_{n \in \mathbb{N}} \subset (0,1) $ and $(\vect{\psi}^{h_n})_{n \in \mathbb{N}}\subset H^1_{\#}(\Omega;\RR^3)$ such that
\begin{equation} \label{appnalk11} 
	\|\pi_{1/h_n} \vect{\psi}^{h_n}\|_{H^1(\Omega;\RR^3)} \geq n\left(\frac{1}{h_n} \|e_h(\vect{\psi}^{h_n})\|_{L^2(\Omega;\RR^{3 \times 3})}+\left|\int_{\Omega} \pi_{1/h_n} (\kappa^{h_n}\vect{\psi}^{h_n}) \right| \right). 
\end{equation} 	
Without a loss of generality, we can assume that for every $n \in \mathbb{N}$ we have 	
\begin{equation} \label{ivo1}
\|\pi_{1/h_n} \vect{\psi}^{h_n}\|_{H^1(\Omega;\RR^3)}=1,\quad \frac{1}{h_n} \|e_h(\vect{\psi}^{h_n})\|_{L^2(\Omega;\RR^{3 \times 3})}+\left|\int_{\Omega} \pi_{1/h_n} (\kappa^{h_n}\vect{\psi}^{h_n}) \right| \leq \frac{1}{n}.  
\end{equation} 
We decompose $\vect{\psi}^h$ according to  \eqref{griso1}, noting that the decomposition satisfies the properties \eqref{griso2} and \eqref{prvaKorn}. Using \eqref{remnak10}, periodicity and \eqref{ivo1} we conclude that 
\begin{equation} \label{ivo2} 
 \left|\int_{\omega} \vect{r}^{h_n} \right| \leq \frac{Ch_n^2}{n}.  \end{equation} 

By using \eqref{ivo2},  Korn's inequality for periodic functions (\CCC see Lemma \ref{kornper} below)\BBB, we conclude from \eqref{remnak10} that \CCC
\begin{equation} \label{ivo3} 
	\frac{1}{h_n}\| \vect{r}^{h_n}\|_{H^1(\omega;\RR^2)} \leq  \frac{C}{n}. 
\end{equation} 
  From \eqref{ivo3} we conclude that $\frac{\vect{r}^{h_n}}{h_n} \to 0$ strongly in $H^1(\omega)$.

From the decomposition \eqref{griso1}, using \eqref{prvaKorn},  \eqref{ivo1} and \eqref{ivo3} we conclude that 
\begin{equation} \label{revoz1} 
	\left| \int_{\omega} \pi_{1/h_n}(\overline{\kappa}^{h_n} \widehat{\vect{\psi}}^{h_n})   \right| \leq \frac{C}{n},
\end{equation} 	
where $\overline{\kappa}^{h_n}=\int_I \kappa^{h_n} \, dx_3$. 
Using \eqref{revoz1} and again \eqref{prvaKorn} and \eqref{remnak10}  and Korn's inequality for periodic functions (Lemma \ref{kornper} below)    for $(\widehat{\psi}_1^{h_n},\widehat{\psi}_2^{h_n})$ and Poincare inequality for $\widehat{\psi}_3^{h_n}$  we conclude
\begin{equation} \label{ivo33} 
\left\|\pi_{1/h_n}\widehat{\vect{\psi}}^{h_n}-\fint_{\omega}\pi_{1/h_n}\widehat{\vect{\psi}}^{h_n}\right\|_{H^1(\omega;\R^3)} \leq  \frac{C}{n}, 
\end{equation}

We denote by $C_n=\fint_{\omega}\pi_{1/h_n}\widehat{\vect{\psi}}^{h_n} $. From \eqref{ivo33} it follows that $\pi_{1/h_n}\widehat{\vect{\psi}}^{h_n} - C_n \to 0$, strongly in $L^2(\omega)$. From \eqref{revoz1}, using the boundedness of $\overline{\kappa}^{h_n}$  we conclude that
 that $C_n \to 0$ as $n\to \infty$. 

It follows that $\pi_{1/h_n}\widehat{\vect{\psi}}^{h_n}\to 0$, strongly in $H^1(\omega)$. Taking into account \eqref{griso1}, \eqref{prvaKorn}, \eqref{ivo1} and \eqref{ivo3} this contradicts the fact that $\|\pi_{1/h_n} \vect{\psi}^{h_n}\|_{H^1(\Omega;\RR^3)}=1$. 
\end{proof} 	
\CCC The following lemma is the Korn's inequality with periodic boundary condition. It is a direct consequence of \cite[Theorem 2.5]{Oleinik}.  We will state it only for dimension two, although it is valid in arbitrary dimension. 
\begin{lema} \label{kornper} 
Let $\omega \subset \R^2$ be a rectangle. Then we have that there exists $C=C(\omega)$ such that 
$$\|\vect{u} \|_{H^1(\omega;\R^2)} \leq C \left(\|e(\vect{u}) \|_{L^2(\omega;\R^{2 \times 2})} +\left|\int_{\omega} \vect{u}   \right|\right), \quad \forall \vect{u} \in H^1_{\#} (\omega;\R^2).    $$
\end{lema} 	
\begin{remarkica}\label{korrem} 
As a trivial consequence of Lemma 	\ref{kornper} we see that, if for $\vect{u} \in H^1_{\#} (\omega;\R^2)$ we have that $e(\vect{u})=0$, then we have $\vect{u}=C \in \R^2$. 
\end{remarkica} 	
\BBB

\section{Auxiliary claims on two-scale convergence}
\label{secaptwoscale}
In this section we prove some auxiliary claims on two-scale convergence. We assume $h>0$ and $\eps=\eps(h) > 0$ are small parameters such that $\lim_{h \to 0} \frac{\eps(h)}{h}=0$. Unless otherwise stated we assume $\Omega \subset \mathbb{R}^3$ with Lipschitz boundary. We take as before $Y=[0,1]^3$ and 
$\mathcal{Y}=\mathbb{R}^3/\mathbb{Z}^3$ a unit flat torus \CCC and $\widehat{Y}=[0,1]^2$, $\widehat{\mathcal{Y}}=\RR^2/\ZZ^2$ (recall Section \ref{notation}).\BBB We give the following two definitions.

\begin{definition}\label{twoscaleconvergencedef1}
Let $(u^h)_{h>0}$ be a bounded sequence in $L^2(\Omega)$. We say that $u^h$ weakly two-scale rescaled converges to $\vect{u}\in L^2(\Omega\times Y)$ if
\begin{equation}
\displaystyle\int_{\Omega}u^h(x)\phi\left(x,\frac{\widehat{x}}{\ee}, \frac{x_3}{\frac{\ee}{h}}\right)\to \displaystyle\int_{\Omega}\int_{Y}u(x,y)\phi(x,y)\,dy\,dx\quad \forall \phi\in C_{c}(\Omega; C(\mathcal{Y})).
\end{equation}
We denote this by $
u^h\xrightharpoonup{2-r}u(x,y). $
Furthermore, we say that $(u^h)_{h>0}$ strongly two-scale converges to $u\in L^2(\Omega\times Y)$ if
$u^h\xrightharpoonup{2-r}u(x,y)$ and $\|u^h\|_{L^2} \to \|u\|_{L^2(\Omega \times Y)}$.
We denote this by $
u^h\xrightarrow{2-r}u(x,y).$
\end{definition}
 \CCC The main motivation behind rescaled two-scale convergence in problem under consideration is that we always use the change of variables to transform the problem on the thin domain to the problem on the domain of unit thickness. In that way the cubes with $\eps$ side length become rectangular cuboid with length and width of size $\eps$ and height of size $\eps/h$.     \BBB

The standard two-scale weak and strong convergence we denote by $\xrightharpoonup{2}$  and $\xrightarrow{2}$ (for these notions see \cite{Allaire-92}). Note that if $\Omega=\omega \times I$, where $\omega \subset \RR^2$ a set with Lipschitz boundary and if $(u^h)_{h>0} \subset L^2(\Omega)$, but depends only on $\hat{x}$, then the notion of weak two-scale convergence and weak two-scale rescaled convergence coincide.  \CCC It is not difficult to check that the standard compactness statement for a sequence bounded in $L^2$ is also valid in this case.  \BBB 
The definition of two-scale rescaled convergence naturally extends to time dependent spaces.

\begin{definition}\label{twoscaleconvergencedef2}
Let $(u^h)_{h>0}$ be a bounded sequence in $L^p(0,T; L^2(\Omega))$, for $p \in (1,+\infty]$. We say that $(u^h)_{h>0}$ weakly two-scale rescaled converges to $u\in L^p(0,T; L^2(\Omega\times Y))$ if
\begin{equation}
\displaystyle\int_0^T\int_{\Omega}u^h(x,t)\phi\left(x,\frac{\widehat{x}}{\ee}, \frac{x_3}{\frac{\ee}{h}}\right)\varphi (t)\,dx\,dt\to \displaystyle\int_0^T\int_{\Omega}\int_{Y}u(t, x,y)\phi(x,y)\varphi (t)\,dy\,dx\,dt,
\end{equation}
$\forall \phi\in C_{c}(\Omega; C(\mathcal{Y})), \varphi\in C(0,T)$. We denote this by 
$u^h\xrightharpoonup{t, 2-r,p}u(x,y,t).$
\CCC Again the standard compactness statement for the sequence bounded in $L^p(0,T;L^2(\Omega))$ is also valid in this case. \BBB

For $p<\infty$, if in addition we have that $u^h(t)\xrightharpoonup{t, 2-r,p} u(x,y,t)$ for almost every $t\in [0,T]$ and
\begin{equation}
\displaystyle\int_{0}^{T}\|u^h\|^p_{L^2(\Omega)}dt\to \displaystyle\int_{0}^{T}\|u\|^p_{L^2(\Omega \times Y)}\,dt,
\end{equation}
then we will say that $(u^h)_{\ee>0}$ strongly two-scale rescaled convergences to $u$ and denote this by $
u^h\xrightarrow{t, 2-r,p}u(x,y,t).$
\end{definition}
We firstly state the following lemma which characterizes the two scale limits of Hessians. \CCC It is well known statement and its proof is e.g. the special case of \cite[Lemma 3]{Vel13}. \BBB
\begin{lema}\label{lemmaA17}
Let $(\varphi^h)_{h>0}\subset H^2(\omega)$ be a bounded sequence. Assume that $\varphi^h\to \varphi_0$ strongly in $L^2(\omega)$.  Then there exists $\varphi_1\in L^2(\omega; \dot{H}^2(\CCC \hat{\mathcal{Y}}\BBB))$ such that on a subsequence we have
\begin{equation}
\nabla^2 \varphi^h\xrightharpoonup{2} \nabla^2\varphi_0 (x)+\nabla^2_{\hat y} \varphi_1(x,\hat{y}).
\end{equation}
\end{lema}

We prove the following lemma and corollary which also characterizes two-scale limits of sequence that satisfies certain properties \CCC (the boundedness of the scaled gradients). \BBB The heuristic argument behind it goes via two-scale expansion, \CCC but the proof goes via duality argument. \BBB

\begin{lema}\label{newlemmamarch21}
Let $\Omega=\omega\times I$, where $\omega\subset\RR^2$ a bounded set with Lipschitz boundary and let $(\vect{V}^h)_{h>0}\subset H^1(\Omega;\RR^3)$ be such that there exists $C>0$ with
\begin{equation}\label{boundmarch231} 
\|\vect{V}^h\|_{L^2(\Omega;\RR^3)}+\ee\| \nabla_h \vect{V}^h\|_{L^2(\Omega;\RR^{3 \times 3})}\leq C.
\end{equation}
Then there exists $\vect{\Psi}\in L^2(\Omega,H^1(\mathcal{Y};\RR^3))$ such that \CCC(up to a subsequence)\BBB
\begin{equation}\label{convermarch231}
\vect{V}^h \xrightharpoonup{2-r} \vect{\Psi}, \quad \ee \nabla_h \vect{V}^h\xrightharpoonup{2-r} \nabla_y \vect{\Psi}.
\end{equation}
\end{lema}

\begin{proof}
We know by two-scale compactness that there exist $\vect{\Psi}\in L^2(\Omega\times\mathcal{Y};\RR^{3})$ and $\mathcal{V}\in L^2(\Omega\times\mathcal{Y};\RR^{3\times 3})$ such that $\vect{V}^h\xrightharpoonup{2-r} \vect{\Psi}$ and $\ee \nabla_h \vect{V}^h\xrightharpoonup{2-r} \mathcal{V}.$
 Let $v\in C_c^{1}(\Omega;\RR)$  and $\vect{\psi}\in C^{1}(\mathcal{Y};\RR^3)$ be arbitrary. We have
 {\allowdisplaybreaks
\begin{equation}
    \begin{split}
        \int_{\Omega}\int_{\mathcal{Y}}\mathcal{V}(x,y)\cdot v(x)\cdot \vect{\psi}(y)dydx
        &=\displaystyle\lim_{h\to 0}\int_{\Omega}\ee\nabla_h\vect{V}^h(x)\cdot v(x)\cdot \vect{\psi}\left(\frac{\widehat{x}}{\ee}, \frac{x_3}{\frac{\ee}{h}}\right)dx
        \\&=-\displaystyle\lim_{h\to 0}\int_{\Omega}\vect{V}^h(x)\cdot v(x)\cdot\left(\partial_{y_1}\psi_1+\partial_{y_2}\psi_2+\partial_{y_3}\psi_3\right)dx
        \\*&\quad-\displaystyle\lim_{h\to 0}\int_{\Omega}\vect{V}^h(x)\cdot \left(\ee \partial_{1}v\psi_1+\ee \partial_{2}v\psi_2+\frac{\ee}{h}\partial_{3}v\psi_3\right)dx
        \\&=-\displaystyle\int_{\Omega}\int_{\mathcal{Y}}\vect{\Psi}(x,y)\cdot v(x)\cdot\left(\partial_{y_1}\psi_1+\partial_{y_2}\psi_2+\partial_{y_3}\psi_3\right)\,dy\,dx
        \\*&=\displaystyle\int_{\Omega}\int_{\mathcal{Y}}\nabla_y\vect{\Psi}(x,y)\cdot v(x)\cdot \vect{\psi}(y)\,dy\,dx.
    \end{split}
\end{equation}
}
From this we have the claim. 
\end{proof}
\CCC
\begin{corollary}\label{cornewlemmamarch21}
	Let $p\in (1,+\infty]$,  $\Omega=\omega\times I$, where $\omega\subset\RR^2$ a bounded set with Lipschitz boundary and let $(\vect{V}^h)_{h>0}\subset L^p(0,T;H^1(\Omega;\RR^3))$ be such that there exists $C>0$ with
	\begin{equation}\label{corboundmarch231} 
		\|\vect{V}^h\|_{L^p(0,T;L^2(\Omega;\RR^3))}+\ee\| \nabla_h \vect{V}^h\|_{L^p(0,T;L^2(\Omega;\RR^{3 \times 3}))}\leq C.
	\end{equation}
	Then there exists $\vect{\Psi}\in L^p(0,T;L^2(\Omega,H^1(\mathcal{Y};\RR^3)))$ such that \CCC(up to a subsequence)
	\begin{equation}\label{corconvermarch231}
		\vect{V}^h \xrightharpoonup{t,2-r,p} \vect{\Psi}, \quad \ee \nabla_h \vect{V}^h\xrightharpoonup{t,2-r,p} \nabla_y \vect{\Psi}.
	\end{equation}
\end{corollary}
\begin{proof} 
	The proof goes in an analogous way as the proof of Lemma \ref{newlemmamarch21}. 
\end{proof} 	
\BBB
The following lemma is in the same spirit as Lemma \ref{newlemmamarch21}. \CCC However, here one imposes the boundedness of the scaled gradients. \BBB

\begin{lema}

\label{newlemmamarch212}
Let $\Omega=\omega\times I$, where $\omega\subset\RR^2$ a bounded set with Lipschitz boundary and let $(\vect{V}^h)_{h>0}\subset H^1(\Omega;\RR^3)$ be such that there exists $C>0$ with
\begin{equation}\label{boundmarch23}
\|\vect{V}^h\|_{L^2(\Omega;\RR^{3})}+\| \nabla_h \vect{V}^h\|_{L^2(\Omega;\RR^{3 \times 3})}\leq C.
\end{equation}
There exist  $\widehat{\vect{V}}\in L^2(\omega, \RR^3)$, $\vect{V}_1\in L^2(\Omega,\RR^3)$,  and $\widetilde{\vect{V}}\in L^2(\Omega,\dot{H}^1(\mathcal{Y},\RR^3))$ such that (\CCC up to a subsequence) \BBB
\begin{equation}
\label{convermarch23}
\nabla_h \vect{V}^h\xrightharpoonup{2-r}\vect{V}:=\left(\nabla_{\widehat{x}}\widehat{\vect{V}}|0\right)+ \nabla_y \widetilde{\vect{V}}+(0|0|\vect{V}_1).
\end{equation}
Here $\widehat{\vect{V}}$ is the \CCC strong limit of $\vect{V}^h$ in $L^2$\BBB \CCC(again on a subsequence). \BBB
\end{lema}

\begin{proof}

 We know by two-scale compactness that there exists $\vect{V}\in L^2(\Omega\times\mathcal{Y};\RR^{3\times 3})$ such that $
\nabla_h \vect{V}^h\xrightharpoonup{2-r} \vect{V}$, \CCC on a subsequence. \BBB
 Let $v\in C^{\infty}_c(\Omega)$  and $\vect{\psi}=(\vect{\psi}^1|\vect{\psi}^2|\vect{\psi}^3)\in C^{\infty}(\mathcal{Y};\RR^{3\times 3})$ be such that $\diver_y \vect{\psi}=0$\footnote{\CCC The divergence is taken per row \BBB}.
 We compute \CCC
\begin{equation}
\label{averagemarch22}
    \begin{split}
         &\int_{\Omega}\int_{\mathcal{Y}}\vect{V}(x,y): \left( v(x) \vect{\psi}(y)\right)dydx =\displaystyle\lim_{h\to 0}\int_{\Omega}\nabla_h\vect{V}^h(x): \left(v(x) \vect{\psi}\left(\frac{\widehat{x}}{\ee}, \frac{x_3}{\frac{\ee}{h}}\right)\right)dx\\
        &=\displaystyle\lim_{h\to 0}\int_{\Omega}\nabla_h \left(\vect{V}^h-\int_{I}\vect{V}^hdx_3\right)(x):\left(v(x) \vect{\psi}\left(\frac{\widehat{x}}{\ee}, \frac{x_3}{\frac{\ee}{h}}\right)\right)\,dx\\
        &\quad+\displaystyle\lim_{h\to 0}\int_{\Omega}\left(\partial_1\left(\int_{I}\vect{V}^hdx_3\right)(x)\Big| \partial_2\left(\int_{I}\vect{V}^hdx_3\right)(x)\Big| 0\right):  \left(v(x) \vect{\psi}\left(\frac{\widehat{x}}{\ee}, \frac{x_3}{\frac{\ee}{h}}\right)\right)dx.
    \end{split}
\end{equation}
\BBB
From \eqref{boundmarch23} it follows that
 \begin{equation}
 \label{eqap1}
\left\|\int_I\vect{V}^h\right\|_{L^2(\omega;\R^3)}\leq C, \quad \left\|\nabla_{\hat{x}} \int_{I} \vect{V}^h\right\|_{L^2(\omega;\R^{3 \times 2})}\leq C,
 \end{equation}
and
\begin{equation}
\label{eqap2}
\displaystyle\int_{\omega}\left|\vect{V}^h-\int_{I}\vect{V}^hdx_3\right|^2\leq \int_{\omega}\int_{I}|\partial_{3}\vect{V}^h|^2\leq Ch^2.
\end{equation}
Note that for $\vect{\Psi}^h:=\vect{V}^h-\int_{I}\vect{V}^hdx_3$, we have that $h^{-1}\vect{\Psi}^h$ is bounded in $L^2(\Omega;\RR^3)$ and $\nabla_h \vect{\Psi}^h$ is bounded in $L^2(\Omega;\RR^{3 \times 3})$. By using Lemma \ref{lemamarch214} below we conclude that on a subsequence we have that\footnote{\CCC The additional regularity of $\vect{\Psi}$ follows from the fact that $h^{-1} \vect{\Psi}^h$ and $h^{-1} \partial_3 \vect{\Psi}^h$ are bounded in $L^2(\Omega)$.  \BBB}  
\begin{equation}
	h^{-1}\vect{\Psi}^h \xrightharpoonup{2-r} \vect{\Psi}(x), \quad \vect{\Psi} \in L^2(\omega;H^1(I;\RR^3)). 
\end{equation}
We can write the first integral on the right-hand side in (\ref{averagemarch22}) as follows \CCC
{\allowdisplaybreaks
\begin{equation}
\label{eq:2303}
    \begin{split}
   &\displaystyle\lim_{h\to 0} \int_{\Omega}\nabla_h \vect{\Psi}^h(x):  \left(v(x) \vect{\psi}\left(\frac{\widehat{x}}{\ee}, \frac{x_3}{\frac{\ee}{h}}\right)\right)dx
    \\&=\displaystyle\lim_{h\to 0}\frac{1}{\eps} \int_{\Omega}\vect{\Psi}^h(x)  v(x)\underbrace{\diver_y    \vect{\psi}\left(\frac{\widehat{x}}{\ee}, \frac{x_3}{\frac{\ee}{h}}\right)}_{=0}dx+\displaystyle\lim_{h\to 0}\int_{\Omega}\frac{\vect{\Psi}^h}{h}(x)\partial_3 v(x)\vect{\psi}^3\left(\frac{\widehat{x}}{\ee}, \frac{x_3}{\frac{\ee}{h}}\right)\,dx\\ &+ \lim_{h \to 0}\int_{\Omega}\vect{\Psi}^h(x) \left(\partial_{1}v(x)\vect{\psi}^1\left(\frac{\widehat{x}}{\ee}, \frac{x_3}{\frac{\ee}{h}}\right)+\partial_{2}v(x)\vect{\psi}^2\left(\frac{\widehat{x}}{\ee}, \frac{x_3}{\frac{\ee}{h}}\right)\right)dx
    \\&=\displaystyle\int_{\Omega}\vect{\Psi}(x)\partial_{3} v(x)\vect{\psi}^3(y)dydx=\displaystyle\int_{\Omega}\vect{\Psi}(x)\partial_{3}v(x)\int_{\mathcal{Y}}\vect{\psi}^3(y)\,dy\,dx
    \\*&=\displaystyle\int_{\Omega}\partial_{3}\vect{\Psi}(x) v(x)\int_{\mathcal{Y}}\vect{\psi}^3(y)\,dy\,dx.
    \end{split}
\end{equation}
}
 \BBB
\CCC From \eqref{eqap1} it follows that (on a subsequence) \BBB  $\int_{I}\vect{V}^hdx_3\xrightarrow{L^2}\widehat{ \vect{V}},$
for some $\widehat{\vect{V}} \in H^1(\omega;\RR^3)$, which is as a consequence of \eqref{eqap2} equivalent to $\vect{V}^h\xrightarrow{L^2} \widehat{\vect{V}}$.
By the basic result of two-scale convergence (see \cite{Allaire-92}), there exists $\vect{V}_2 \in L^2(\omega; H^1(\widehat{\mathcal{Y}};\RR^3))$ such that
\begin{equation}
\label{twoscale234}
\nabla_{\hat{x}}\int_{I}\vect{V}^hdx_3\xrightharpoonup{2-r} \nabla_{\widehat{x}}\widehat{\vect{V}}+\nabla_{\widehat{y}}\vect{V}_2.
\end{equation}
 Using the convergences \eqref{eq:2303} and \eqref{twoscale234} in \eqref{averagemarch22}, we obtain
 \begin{equation}
     \begin{split}
         	&\int_{\Omega}\int_{\mathcal{Y}}\vect{V}(x,y): v(x) \vect{\psi}(y)\,dy\,dx\\
	&=\displaystyle\int_{\Omega}\partial_{3}\vect{\Psi}(x) v(x)dx\,\int_{\mathcal{Y}}\vect{\psi}^3(y)\, dy+\int_{\Omega}\int_{\mathcal{Y}}\left(\nabla_{\widehat{x}}\widehat{\vect{V}}+\nabla_{\widehat{y}}\vect{V}_2|0\right): v(\vect{\psi}^1 |\vect{\psi}^2| \vect{\psi}^3)\,dy\,dx\\
	&=\displaystyle\int_{\Omega}\int_{\mathcal{Y}}\left[\left(\nabla_{\widehat{x}}\widehat{\vect{V}}|\partial_{3}\vect{\Psi}\right)+\left(\nabla_{\widehat{y}}\vect{V}_2|0\right)\right]: v \vect{\psi}\,dy\,dx.
     \end{split}
 \end{equation}
Therefore
\begin{eqnarray*}
& &	\int_{\Omega}\int_{\mathcal{Y}}\left\{\vect{V}(x,y)-\left[\left(\nabla_{\widehat{x}}\widehat{\vect{V}}|\partial_{3}\vect{\Psi}\right)+\left(\nabla_{\widehat{y}}\vect{V}_2|0\right)\right]\right\}: v  \vect{\psi}\,dy\,dx=0,\\ & &\quad \quad \forall v \in C_c (\Omega),\  \forall \vect{\psi}\in  C^{\infty}(\mathcal{Y};\RR^{3\times 3})\,\,\text{such that}\,\,\diver_y \vect{\psi}=0,
\end{eqnarray*}
which implies that there exists $\vect{V}_3\in L^2(\Omega;\dot{H}^1(\mathcal{Y};\RR^3))$ such that
\begin{equation}\label{claim2}
	\vect{V}(x,y)-\left[\left(\nabla_{\widehat{x}}\widehat{\vect{V}}|\partial_{3}\vect{\Psi}\right)+\left(\nabla_{\widehat{y}}\vect{V}_2|0\right)\right]=\nabla_y \vect{V}_3.
\end{equation}
Finally, the lemma follows directly from \eqref{claim2} by setting $\widetilde{\vect{V}}:=\vect{V}_2+\vect{V}_3$, $\vect{V}_1:=\partial_3 \vect{\Psi}$.
\end{proof}
The following lemma was necessary for the proof of Lemma \ref{newlemmamarch212}.  \CCC It states that under certain condition we can guarantee that two-scale limit of the sequence doesn't depend on the fast variable. \BBB
\begin{lema}\label{lemamarch214}
Let $({\varphi}^{h})_{h>0}\subset H^1(\Omega)$ be a bounded sequence in $L^2(\Omega)$ such that $\varphi^{h}\overset{2-r}{\rightharpoonup} \varphi(x,y)$. Assume that \CCC $\ee\nabla_h\varphi^{h}\rightarrow 0$ \BBB strongly in $L^2(\Omega)$. Then we have that $\varphi(x,y)$ depends only on $x$.
\end{lema}

\begin{proof}
Let $k\neq 0$, $k \in \ZZ^3$, and $v\in C^{\infty}_{c}(\Omega)$. Without loss of generality we assume $k_1 \neq 0$. We define $$p(y):=\frac{\sin  2\pi k y}{k_1}-i\frac{\cos 2\pi k y}{k_1}.$$
Next we have:
\begin{align}
&\displaystyle\int_{\Omega}\int_{Y}\varphi(x,y)\cdot e^{i 2\pi\langle k,y\rangle}v(x)\,dy\,dx\\
&=\int_{\Omega}\int_{\mathcal{Y}}\varphi(x,y)\cdot\partial_{y_1}p(y)\cdot v(x)dxdy=\displaystyle \lim_{h\to 0}\int_{\Omega}\varphi^{h}(x)\cdot \partial_{y_1}p\left(\frac{\widehat{x}}{\ee}, \frac{x_3}{\frac{\ee}{h}}\right)\cdot v(x)\,dx\\
&=\displaystyle\lim_{h \to 0}\int_{\Omega}\varphi^{h}(x)\cdot \partial_{1}\left(\ee p\left(\frac{\widehat{x}}{\ee}, \frac{x_3}{\frac{\ee}{h}}\right)\cdot v(x)\right)dx-\displaystyle\lim_{h\to 0}\int_{\Omega}\varphi^{h}(x)\cdot \ee p\left(\frac{\widehat{x}}{\ee}, \frac{x_3}{\frac{\ee}{h}}\right)\cdot \partial_{1}v(x)\,dx\\
&=\displaystyle\lim_{h\to 0}\int_{\Omega}\varphi^{h}(x)\cdot \partial_{1}\left(\ee p\left(\frac{\widehat{x}}{\ee}, \frac{x_3}{\frac{\ee}{h}}\right)\cdot v(x)\right)dx=\displaystyle\lim_{h \to 0}\int_{\Omega}\ee \partial_{1}\varphi^{h}(x)\cdot p\left(\frac{\widehat{x}}{\ee}, \frac{x_3}{\frac{\ee}{h}}\right)\cdot v(x)\,dx=0.
\end{align}
In the analogous way we treat the case when $k_2 \neq 0$ or $k_3 \neq 0$. Thus we conclude that $\varphi$ depends only on $x$. 
\end{proof}
 \end{appendices}

\section*{Acknowledgements}
The authors were supported by Croatian Science Foundation under Grant Agreement no. IP-2018-01-8904 (Homdirestroptcm). P. Hern\'andez-Llanos was also supported by Agencia Nacional de Investigaci\'on y Desarrollo, FONDECYT Postdoctorado 2023 under Grant no. 3230202. 
 The warm hospitality at the University of Zagreb is gratefully acknowledged by P Hern\'andez-Llanos. M. Bu\v{z}an\v{c}i\'{c}, I. Vel\v{c}i\'{c} and Josip \v{Z}ubrini\'{c} were also supported by  Croatian Science Foundation under Grant agreement No. IP-2022-10-5181 (HOMeOS).
 \CCC We are grateful to anonymous referee for his suggestions that helped in improving the paper. \BBB
 \\
 
\vspace{+1ex} 
\noindent \textbf{Data Availability:} Due to the nature of the research, there is no supporting data.\\

\vspace{+1ex} 
\noindent \textbf{Conflict of Interest:} The authors declare to have no conflict of interest.




\end{document}